\documentstyle{amsppt}
\newcount\mgnf\newcount\tipi\newcount\tipoformule\newcount\greco 
\tipi=2          
\tipoformule=0   

\global\newcount\numsec\global\newcount\numfor
\global\newcount\numapp\global\newcount\numcap
\global\newcount\numfig\global\newcount\numpag
\global\newcount\numnf

\def\SIA #1,#2,#3 {\senondefinito{#1#2}%
\expandafter\xdef\csname #1#2\endcsname{#3}\else
\write16{???? ma #1,#2 e' gia' stato definito !!!!} \fi}

\def \FU(#1)#2{\SIA fu,#1,#2 }

\def\etichetta(#1){(\veroparagrafo.\veraformula)%
\SIA e,#1,(\veroparagrafo.\veraformula) %
\global\advance\numfor by 1%
\write15{\string\FU (#1){\equ(#1)}}%
\write16{ EQ #1 ==> \equ(#1)  }}
\def\etichettaa(#1){(A\veraappendice.\veraformula)
 \SIA e,#1,(A\veraappendice.\veraformula)
 \global\advance\numfor by 1
 \write15{\string\FU (#1){\equ(#1)}}
 \write16{ EQ #1 ==> \equ(#1) }}
\def\getichetta(#1){Fig. \verafigura
 \SIA g,#1,{\verafigura}
 \global\advance\numfig by 1
 \write15{\string\FU (#1){\graf(#1)}}
 \write16{ Fig. #1 ==> \graf(#1) }}
\def\retichetta(#1){\numpag=\pgn\SIA r,#1,{\verapagina}
 \write15{\string\FU (#1){\rif(#1)}}
 \write16{\rif(#1) ha simbolo  #1  }}
\def\etichettan(#1){(n\verocapitolo.\veranformula)
 \SIA e,#1,(n\verocapitolo.\veranformula)
 \global\advance\numnf by 1
\write16{\equ(#1) <= #1  }}

\newdimen\gwidth
\gdef\profonditastruttura{\dp\strutbox}
\def\senondefinito#1{\expandafter\ifx\csname#1\endcsname\relax}
\def\BOZZA{
\def\alato(##1){
 {\vtop to \profonditastruttura{\baselineskip
 \profonditastruttura\vss
 \rlap{\kern-\hsize\kern-1.2truecm{$\scriptstyle##1$}}}}}
\def\galato(##1){ \gwidth=\hsize \divide\gwidth by 2
 {\vtop to \profonditastruttura{\baselineskip
 \profonditastruttura\vss
 \rlap{\kern-\gwidth\kern-1.2truecm{$\scriptstyle##1$}}}}}
\def\verapagina{
{\romannumeral\number\numcap}.\number\numsec.\number\numpag}}

\def\alato(#1){}
\def\galato(#1){}
\def\veroparagrafo{\number\numsec}\def\veraformula{\number\numfor}
\def\veraappendice{\number\numapp}
\def\verapagina{\number\pageno}\def\veranformula{\number\numnf}
\def\verafigura{{\romannumeral\number\numcap}.\number\numfig}
\def\verocapitolo{\number\numcap}\def\veranformula{\number\numnf}
\def\Eqn(#1){\eqno{\etichettan(#1)\alato(#1)}}
\def\eqn(#1){\etichettan(#1)\alato(#1)}
\def\ver{\veroparagrafo}
\def\Eq(#1){\eqno{\etichetta(#1)\alato(#1)}}
\def\eq(#1){\etichetta(#1)\alato(#1)}
\def\Eqa(#1){\eqno{\etichettaa(#1)\alato(#1)}}
\def\eqa(#1){\etichettaa(#1)\alato(#1)}
\def\dgraf(#1){\getichetta(#1)\galato(#1)}
\def\drif(#1){\retichetta(#1)}

\def\eqv(#1){\senondefinito{fu#1}$\clubsuit$#1\else\csname fu#1\endcsname\fi}
\def\equ(#1){\senondefinito{e#1}\eqv(#1)\else\csname e#1\endcsname\fi}
\def\graf(#1){\senondefinito{g#1}\eqv(#1)\else\csname g#1\endcsname\fi}
\def\rif(#1){\senondefinito{r#1}\eqv(#1)\else\csname r#1\endcsname\fi}
\def\bib[#1]{[#1]\numpag=\pgn
\write13{\string[#1],\verapagina}}

\def\include#1{
\openin13=#1.aux \ifeof13 \relax \else
\input #1.aux \closein13 \fi}

\openin14=\jobname.aux \ifeof14 \relax \else
\input \jobname.aux \closein14 \fi
\openout15=\jobname.aux
\openout13=\jobname.bib


\ifnum\tipoformule=1\let\Eq=\eqno\def\eq{}\let\Eqa=\eqno\def\eqa{}
\def\equ{}\fi


{\count255=\time\divide\count255 by 60 \xdef\hourmin{\number\count255}
        \multiply\count255 by-60\advance\count255 by\time
   \xdef\hourmin{\hourmin:\ifnum\count255<10 0\fi\the\count255}}

\def\oramin{\hourmin }

\def\data{\number\day/\ifcase\month\or january \or february \or march \or
april \or may \or june \or july \or august \or september
\or october \or november \or december \fi/\number\year;\ \oramin}

\setbox200\hbox{$\scriptscriptstyle \data $}

\newcount\pgn \pgn=1
\def\foglio{\number\numsec:\number\pgn
\global\advance\pgn by 1}
\def\foglioa{A\number\numsec:\number\pgn
\global\advance\pgn by 1}

\footline={\rlap{\hbox{\copy200}}\hss\tenrm\folio\hss}


\global\newcount\numpunt

\magnification=\magstephalf
\baselineskip=16pt
\parskip=8pt

\voffset=2.5truepc
\hoffset=0.5truepc
\hsize=6.1truein
\vsize=8.4truein 
{\headline={\ifodd\pageno\rightheadline \else \leftheadline \fi}}
\def\rightheadline{\it  {tralala}\hfil\tenrm\folio}
\def\leftheadline{\tenrm \folio \hfil\it  {Section $\ver$}}

\def\a{\alpha}

\def\d{\delta}
\def\e{\epsilon}

\def\f{\phi}
\def\g{\gamma}

\def\l{\lambda}

\def\s{\sigma}
\def\t{\tau}

\def\z{\zeta}

\def\D{\Delta}
\def\L{\Lambda}
\def\G{\Gamma}

\def\1{{1\kern-.25em\roman{I}}}
\def\eu{{1\kern-.25em\roman{I}}}
\def\f1{{1\kern-.25em\roman{I}}}

\def\R{{\Bbb R}}  
\def\N{{\Bbb N}}  
\def\Z{{\Bbb Z}}  
\def\E{{\Bbb E}}  

\def\dist{\,\roman{dist}}

\let\cal=\Cal
\def\AA{{\cal A}}
\def\BB{{\cal B}}

\def\DD{{\cal D}}
\def\EE{{\cal E}}
\def\FF{{\cal F}}
\def\GG{{\cal G}}

\def\LL{{\cal L}}
\def\MM{{\cal M}}

\def\PP{{\cal P}}

\def\SS{{\cal S}}

\def\A{{\cal A}}

\def\chap #1#2{\line{\ch #1\hfill}\numsec=#2\numfor=1}

\def\ov #1{\overline{#1}}
\def\ba{{\backslash}}

\def\wt{\widetilde}

\def\ov#1{\overline{#1}}


\def\note#1{\footnote{#1}}

\def\frac#1#2{{#1\over #2}}
\def\sfrac#1#2{{\textstyle{#1\over #2}}}

\def\text#1{\quad{\hbox{#1}}\quad}
\def\newpage{\vfill\eject}
\def\proposition #1{\noindent{\thbf Proposition #1:}}

\def\theo #1{\noindent{\thbf Theorem #1: }}

\def\lemma #1{\noindent{\thbf Lemma #1: }}
\def\definition #1{\noindent{\thbf Definition #1: }}

\def\proof{{\noindent\pr Proof: }}
\def\proofof #1{{\noindent\pr Proof of #1: }}
\def\endproof{$\diamondsuit$}
\def\remark{\noindent{\bf Remark: }}
\def\thanks{\noindent{\bf Acknowledgements: }}
\def\lg{\buildrel {\textstyle <}\over  {>}}
\font\pr=cmbxsl10
\font\thbf=cmbxsl10 scaled\magstephalf

\font\ch=cmbx12

\font\it=cmti10
\font\bf=cmbx10


 \def\lg{\buildrel {\textstyle <}\over  {>}}
\def\dom{{\roman {dom}}}
\def\diam{\,{\roman {diam}}}
\def\conv{\,{\roman {conv}}}
\def\dist{\,{\roman {dist}}}
\def\inte{\,{\roman {int}}}
\def\ri{\,{\roman {ri}}}
\def\cl{\,{\roman {cl}}}
\def\bd{\,{\roman {bd}}}
\def\rbd{\,{\roman {rbd}}}

\def\supp{\,{\roman {supp}}}
\def\st{\,\,{\roman {s.t.}}\,}
\overfullrule=0pt
\nopagenumbers
\font\tit=cmbx12
\font\aut=cmbx12

\def\s{\char'31}
\centerline{\tit SAMPLE PATH LARGE DEVIATIONS FOR A CLASS OF}
\vskip.2truecm
\centerline{\tit MARKOV CHAINS RELATED TO }
\vskip0.2cm
\centerline{\tit  DISORDERED MEAN FIELD MODELS}
\vskip.2truecm 
\vskip1.5truecm

\centerline{\aut Anton Bovier 
\note{ Weierstrass-Institut f\"ur Angewandte Analysis und Stochastik,
Mohrenstrasse 39, D-10117 Berlin,\hfill\break Germany.
 e-mail: bovier\@wias-berlin.de},
 and  V\'eronique Gayrard\note{D\'epartement de math\'ematiques,
E.P.F.L., CH-1015 Lausanne, Switzerland; on leave from
Centre de Physique Th\'eorique, CNRS,
Luminy, Case 907, F-13288 Marseille, Cedex 9, France.\hfill\break
email: gayrard\@masg1.epfl.ch}}

\vskip2truecm\rm
\def\s{\sigma}
\noindent {\bf Abstract:} We prove a large deviation principle 
on path space for a class of discrete time Markov processes 
whose state space is the intersection of a regular domain $\L\subset \R^d$
with some lattice of spacing $\e$. Transitions from $x$ to $y$ are allowed if
$\e^{-1}(x-y)\in \D$ for some fixed set of vectors $\D$.
The transition 
probabilities $p_\e(t,x,y)$, which themselves depend on $\e$,
 are allowed to depend on the starting 
point $x$ and the time $t$ in a sufficiently regular way, 
except near the boundaries, where some singular behaviour 
is allowed. The rate function is 
identified as an action functional which is given as the integral of 
a Lagrange function.
Markov processes of this type arise in the study of mean field dynamics
of disordered mean field models.


 \vfill
$ {} $
\newpage
{\headline={\ifodd\pageno\rightheadline \else \leftheadline \fi}}
\def\rightheadline{\it  {Sample path LDP}\hfil\tenrm\folio}
\def\leftheadline{\tenrm \folio \hfil\it  {Section $\ver$}}

\chap{1. Introduction.}1

 In this paper we study a class of Markov processes with discrete state
space which have the property that their transition probabilities 
vary slowly with time as the processes progresses (we will give a precise
meaning to this later). Such processes occur 
in many applications and have been studied both in the physical 
and mathematical literature. For an extensive discussion, we refer e.g. 
to van Kampen's book [vK], Chapter IX. It has been shown 
by Kurtz [Ku], under suitable
conditions, that these processes can be scaled in such a way
that a law of large numbers holds that states that the rescaled process 
converges, almost surely, to the solution of a certain differential equation.
He also established a central limit theorem showing that the deviations 
from the solution under proper scaling converges to a
generalized Ornstein-Uhlenbeck process [Ku2].
The simplest example of such Markov processes are of course symmetric 
random walks
(in $\Z^d$, say). In this case one
the LLN scaling consists in considering the process
(for $t\in \R_+$)
 $Z_n(t)=\frac1n\sum_{i=1}^{[nt]}X_i$, and one has the obvious result that as 
$n$ tends to infinity, $Z_n(t)$ converges to $0$, which solves the
differential equation
is $X'(t)=0$. The corresponding central limit theorem is then nothing but 
Donsker's invariance principle [Do] which asserts that $\sqrt n Z_n(t)$
converges to Brownian motion. In this simple situation, the LLN and the CLT 
are accompanied by a large deviation principle, due to Mogulskii [Mo]
that states that the family of laws of the processes  $Z_n(t), t\in [0,T]$ 
satisfies
a large deviation principle with some rate function of the form
$S(x)=\int_0^T dt \LL(\dot x(t))$. This LDP is the analog of 
Schilder's theorem for Brownian motion (in which case the function 
$\LL$ is just the square). 
Generalizations of Mogulskii's theorem were studied in a series of paper by 
Wentzell [W1-4]. A partial account of this work is given in Section 5 of the
book by Wentzell and Freidlin [WF]. The class of locally infinitely divisible
processes studied there include Markov jump processes. Wentzell 
proved large deviation principles under some spatial 
regularity assumptions on the moment generating functions of the local 
jump-distributions and its Legendre transforms. 
The particular case of pure Markov jump processes is worked out in [SW].
This theory has been developed considerably in a large number of works
principlly by 
Dupois,  Ellis, and Weiss and co-workers (see e.g. [DEW,DE,DE1,DE2,DR,AD,SW] 
and references therein). The main thrust of this line of research was to weaken
the spatial regularity hypothesis on the transition rates to
include situations with boundaries and discontinuities. The main motivation
was furnished by applications to queing systems. Given the variety of 
possibls situations, is not surprising that there is no complete theory 
availble, but rather a large set of examples satisfying particular 
hypothesis. Among the rare general results is 
an upper large deviation bound proven in 
  [DEW] that holds  under measurability assumptions only; the question under 
which conditions
these bounds are sharp remain  open in general. The upper bounds in [DEW] 
are also stated for discrete time Markov processes. 
Needless to say, the bulk of the literature is concerned with the diffusion
 case, i.e. large deviations for solutions of stochastic differential 
equations driven by Wiener processes [WF,Az]. Questions of discontinuous 
statistics have been considered in this context in [BDE,CS]. 
For other related 
large deviation principles, see also [Ki1,Ki2].

In the present case we consider discrete time Markov chains depending on 
a small  parameter $\e$ defined on a state space $\L_\e\subset \R^{d}$
that have transition rates $p_\e(x,y,t)$ of the form 
$p_\e(x,x+\e\d,t)=\exp(f_\e(x,\d,t))$, for $\d\in\D$ where
$\D$ is some finite   set and $f_\e$
is required to satisfy some regularity conditions to be specified
in detail later. The new feature of our results are
\item{(i)} The functions $f_\e$ themselves 
are allowed to depend (in a controlled way)
on the small parameter $\e$.
\item{(ii)} Regularity conditions are required in the interiors of 
the domains, but some singular behaviour near the boundary is allowed.
\item{(iii)} The transition rates are time-dependent.

Features (i) and (ii) are motivated from applications to stochastic dynamics in
disordered mean-field models of statistical mechanics which we will not 
discuss here. See e.g. [BEGK,BG]. Let us mention that the 
large deviations results obtained in the present paper were needed
(in the particular setting of time-homogeneous and reversible processes)
in [BEGK] to show that a general transition between metastable states
proceeds along a (asymptotically) deterministic sequence of 
so-called admissible transitions. The necessity to consider (i) arises mainly
from the
fact that in such systems, rather strong finite size effect due to the disorder
are present and these effect the transition probabilities. Control of this
dependence requires a certain amount of extra work.

 The problem at boundaries (ii) is also intrinsic for most of the systems 
we are interested in. While for many application it would  be sufficient 
to have an large deviation estimates for sets of paths that stay away from the 
boundary, we feel that it is more satisfactory to have a full LDP under 
conditions that are generally met in the systems we are interested in.
The types of singularities we must deal with differ from those treated 
in the queing motivated literature cited above. 

(iii) is motivated by our interest to 
study the behaviour of such systems under time dependent external
variations of parameters, and in particular to study  hysteresis phenomena.
This causes no particular additional technical difficulties.


We have chosen to give complete and elementary proves of our results, even 
though the basic ideas are now standard in large deviation theory and 
any technical lemmata (mainly from convex analysis) are also served in 
similar situations in the past. But there are some subtle points, mainly in 
the dealing with boundary effects, and we feel
that it is easier and more instructive to follow a complete 
line of argument using only the minimal amount of technical tools. 

The remainder of the paper is organized as follows. In Section 2 we 
give precise formulation of our results. Section 3 states the basic large 
deviation upper and lower bounds and shows why they imply our theorems, 
Section 4 establishes some 
elementary fact from convex analysis that will be needed later, 
and in Section 5 the upper and lower bounds are proven.

\thanks We thank J.-D. Deuschel and O. Zeitouni for pointing out some 
interesting recent references. This paper was written during 
visits of the authors at the  Centre de Physique Th\'eorique,
Marseille, the Weierstrass Institut f\"ur Angewandte Analysis, Berlin, the 
D\'epartment de math\'ematiques de l'EPFL, Lausanne, and EURANDOM, Eindhoven.
We thank these institutions for their hospitality and financial support.

\newpage


\chap{2. Statements of results}2

Let $\G$ denote some lattice in $\R^d$ and let
 $\L\subset\R^d$ be a convex set (finite or infinite)
that is complete w.r.t. the 
Euclidean metric. Define, for $\e>0$, the rescaled lattice $\e \G$ and its 
intersection with the set $\L$, $\G_{\e}\equiv\L\cap (\e\G)$. We consider 
discrete time Markov chains with state space $\G_\e$. $\e$ will play the 
r\^ole of a small parameter\note{In applications to dynamics of
mean field models $\e$ will enter as the the inverse of the 
system size $N$, hence only take discrete values. 
This will not be important here.}. Let $\D\subset \G$ 
denote a finite subset of lattice vectors.

The time $t$-to-$(t+1)$ transition probabilities, 
$(t,x,y)\in\N\times\G_{\e}\times\G_{\e}\mapsto p_{\e}(t,x,y)\in [0,1]$
are of the form
$$
p_{\e}(t,x,y)
\equiv
\cases
g_{\e}\left(\e t,x,{\e}^{-1}(y-x)\right)
& \hbox{if}\,\, {\e}^{-1}(y-x)\in\D, x\in\G_{\e}, y\in\G_{\e}\cr
&\cr
0 &\hbox{ otherwise}
\endcases
\Eq(2.1)
$$
where the functions $\{g_{\e}, \e>0\}$, 
$g_{\e}:\R^+\times\R^d\times\D\rightarrow[-\infty, 0]$, 
are obviously required to meet 
the condition
$$
\sum_{\d\in\D}g_{\e}({s},x,\d)=1
,\,\,\,\,\,\,\,\,\forall {s}\in\R^+, \forall x\in\L
\Eq(2.2)
$$
We will set 
$$
f_\e(t,x,\d)\equiv \cases \ln g_\e(t,x,\d), &\hbox{\rm if} \,\, g_\e(t,x,\d)>0
\cr -\infty,&\hbox{\rm if} \,\,  g_\e(t,x,\d)=0
\endcases
\Eq(2.2bis)
$$
These functions  will be assumed to verify a number of 
additional hypothesis; in order to 
state them we need some notation:
For any set ${\cal S}$ in $\R^d$ the {\it convex hull} of ${\cal S}$ is denoted
by $\conv{\cal S}$; the {\it closure, interior} and {\it boundary} of 
${\cal S}$
are denoted by $\cl\, {\cal S}$, $\inte\, {\cal S}$ and 
$\bd\, {\cal S}=(\cl\, {\cal S})\setminus(\inte\, {\cal S})$. For each $\e>0$ 
we define the $\e$-{\it interior} of ${\cal S}$, denoted by 
$\inte_\e {\cal S}$, to be:
$$
\inte_\e {\cal S}=\{x\in{\cal S}\mid \forall\d\in\D, x+\e\d\in{\cal S}\}
\Eq(2.3)
$$
Note that $\inte_\e {\cal S}$ is not necessarily open. 
The $\e$-{\it boundary} of ${\cal S}$ is then defined by 
$\bd_{\e}{\cal S}=(\cl\, {\cal S})\setminus(\inte_\e {\cal S})$.
For each $\e>0$ we set:
$$
\eqalign{
&\L^{(\d,\e)}=\left\{x\in\L\mid x+\e\d\in\L\right\},
\,\,\,\,\,\d\in\D\cr
&\L^{(\d)}=\left\{x\in\L\mid \exists \e>0 \st x+\e\d\in\L\right\},
\,\,\,\,\,\d\in\D\cr
}
\Eq(2.4)
$$
Obviously
$$
\eqalign{
\bigcup_{\d\in\D}\L^{(\d,\e)}=\L
&{\text{and}}
\bigcap_{\d\in\D}\L^{(\d,\e)}=\inte_{\e}\L
\cr
\bigcup_{\d\in\D}\L^{(\d)}=\L
&{\text{and}}
\bigcap_{\d\in\D}\L^{(\d)}=\inte\L
\cr
}
\Eq(2.5)
$$
Moreover we have:

\lemma{\ver.1}{\it
\item{(i)} $\inte_{\e'}\L\subset \inte_{\e}\L$ for all $\e'>\e>0$;
\item{(ii)} $\inte_{\e}\L\subset \inte\L$ for all $\e>0$;
\item{(iii)}
$
\inte\L=
\left\{x\in\L\mid \exists \e>0 \st x\in\inte_{\e}\L\right\}
$.
}

\proof (i) is immediate. Given $x\in\inte_{\e}\L$ each of the points
$x+\e\d$, $\d\in\D$, belongs to $\L$. Forming the convex hull of this 
set of points we have, by convexity of $\L$: 
$\conv\{x+\e\d\mid \d\in\D\}=x+\e\,\conv\D\subset\L$. Let $B$ be the closed 
unit ball in $\R^d$ centered at the origin. Since by assumption
$\conv\D$ is a $d$-dimensional set, there exists $r\equiv r(\diam\D)>0$ such 
that $r B\subset\conv\D$. Hence $x+\e r B\subset\L$ and
$\inte_{\e}\L\subset \{x\in\L\mid x+\e r B\subset\L\}$, proving (ii). Similarly
we obtain that for any 
$x\in\bigcup_{\e>0}\inte_{\e}\L
=\{x\in\L\mid \exists \e>0 \st \forall \d\in\D,
x+\e\d\in\L\}$ there exists $\e'>0$
such that $x+\e' B\subset\L$, which yields (iii). The lemma is proven. 
\endproof.

\noindent{\bf Hypothesis  {\ver.2}:}\note{The statements ``for each $\e>0$''
should in fact be replaced by ``for each $\e>0$ sufficiently small''.}
{\it  For each $\e>0$ and each $\d\in\D$,
$$
\eqalign{
g_{\e}({s},x,\d)&>0,\,\,\,\,\,\forall 
({s},x)\in\R^+\times\L^{(\d,\e)}
\cr
g_{\e}({s},x,\d)&=0,\,\,\,\,\,\forall 
({s},x)\in\R^+\times\L\setminus\L^{(\d,\e)}
\cr
}
\Eq(2.6)
$$
and 
$$
g_{\e}({s},x,\d)=0,\,\,\,\,\,\forall 
({s},x,\d)\in\R^+\times(\R^d\setminus\L)\times\D
\Eq(2.7)
$$
Moreover,
\item{}  $\forall x\in\inte\L, \exists \e'>0$ and $c>0$
 such that $\forall 0<\e<\e'$,
$$
g_{\e}({s},x,\d)>c,\,\,\,\,\,\forall ({s},\d)\in\R^+\times\D
\Eq(2.8bis)
$$
\item{} $\forall x\in\bd\L,\exists \e'>0$ and $c>0$ such that 
$\forall 0<\e<\e'$,
$$
g_{\e}({s},x,\d)>c,\,\,\,\,\,\forall 
{s}\in\R^+, \forall \d\in\{\d'\in\D\mid \L^{(\d')}\ni x\}
\Eq(2.8ter)
$$
and 
$$
g_{\e}({s},x,\d)=0,\,\,\,\,\,\forall 
{s}\in\R^+, \forall \d\notin\{\d'\in\D\mid \L^{(\d')}\ni x\}
\Eq(2.8qua)
$$

}

\remark Hypothesis \ver.2 implies in particular that for each $\e>0$,
$$
f_{\e}({s},x,\d)>-\infty,\,\,\,\,\,\forall 
({s},x,\d)\in\R^+\times\inte_{\e}\L\times\D
\Eq(2.8)
$$
and 
$$
\forall x\in\bd_{\e}\L, \exists \d\in\D {\text{s.t.}} f_{\e}({s},x,\d)>-\infty
\Eq(2.9)
$$

\remark  Lemma \ver.1 and Hypothesis \ver.2 also imply that for 
any $x\in\L$, $\exists \e'>0$ s.t. $\forall 0<\e<\e'$
$$
\left\{\d\in\D\,\big|\,\L^{(\d,\e)}\ni x\right\}
=\left\{\d\in\D\,\big|\,\L^{(\d)}\ni x\right\}
\Eq(2.9bis)
$$

\noindent{\bf Hypothesis  {\ver.3}:}
{\it 
There exist functions, $f^{(0)}_{\e}$ and $f^{(1)}_{\e}$ such that
$$
f_{\e}=f^{(0)}_{\e}+\e f^{(1)}_{\e},
\Eq(2.10)
$$
satisfying:

\item{(H0)} $f^{(0)}_{\e}({s},x,\d)=-\infty$ if and only if 
$f_{\e}({s},x,\d)=-\infty$.
\item{(H1)} 
For any closed bounded subset ${\cal S}\subset\inte\L$ there exists 
a positive constant $K\equiv K({\cal S})<\infty$ such that, for each $\e>0$,
$$
\sup_{x\in{\cal S}}
\sup_{\scriptstyle{\d\in\D:}\atop\scriptstyle{{\cal S}\cap\L^{(\d,\e)}\ni x}}
\left|f^{(1)}_{\e}({s},x,\d)\right|\leq K,\,\,\,\,\forall s\in\R^+
\Eq(2.11)
$$

\item{(H2)} 
There exists a constant $0<\theta<\infty$ such that, for each $\e>0$,
$$
\sup_{x\in\L}
\sup_{\scriptstyle{\d\in\D:}\atop\scriptstyle{\L^{(\d,\e)}\ni x}}
\left|f^{(0)}_{\e}({s},x,\d)-f^{(0)}_{\e}({s}',x,\d)\right|
\leq\theta|{s}-{s}'|
,\,\,\,\,\,\,\,\,\forall {s}\in\R^+, \forall {s}'\in\R^+\,.
\Eq(2.12)
$$
\item{(H3)} 
For any closed bounded subset ${\cal S}\subset\inte\L$ there exists 
a positive constant $\vartheta\equiv \vartheta({\cal S})<\infty$ such that,
for each $\e>0$,
$$
\sup_{{s}\in\R^+}
\sup_{\scriptstyle{\d\in\D:}\atop\scriptstyle
{{\cal S}\cap\L^{(\d,\e)}\ni \{x,x'\}}}
\left|f^{(0)}_{\e}({s},x,\d)-f^{(0)}_{\e}({s},x',\d)\right|
\leq\vartheta|x-x'|
,\,\,\,\,\,\,\,\,\forall x\in{\cal S}, \forall x'\in{\cal S}
\Eq(2.13)
$$
}

\noindent{\bf Hypothesis  {\ver.4}:} {\it 
The functions $g_\e$ converge uniformly to a function $g$ on the set
$\R^+\times\L\times \D$. 
Moreover, for any  $({s},x,\d)\in \R^+\times\L\times\D$, 
$$
\lim_{\e\rightarrow 0}g_{\e}({s},x,\d)
=\lim_{\e\rightarrow 0}e^{f^{(0)}_{\e}({s},x,\d)}
\Eq(2.a)
$$
}

\remark Note that Hypothesis \ver.4 together with Hypothesis \ver.2
implies that the limits
$$
\lim_{\e\rightarrow 0}f_{\e}({s},x,\d)
=\lim_{\e\rightarrow 0}f^{(0)}_{\e}({s},x,\d)=f({s},x,\d)
\Eq(2.14)
$$
exist and are finite at every $({s},x,\d)$ in the set defined by:
$$
s\in\R^+, x\in\L, \d\in\{\d'\in\D\mid \L^{(\d')}\ni x\}
\Eq(2.14bis)
$$
We put $f({s},x,\d)=-\infty$
on the complement of \eqv(2.14bis).

\remark For $x\in\inte\L$ then $\{\d'\in\D\mid \L^{(\d')}\ni x\}=\D$.

\remark The limiting function $f$ of course inherits the properties $(H2)$ and 
$(H3)$ of Hypothesis \ver.3 with $\L^{(\d,\e)}$ replaced by $\L^{(\d)}$.

 As a consequence of Hypothesis \ver.3  and \ver.4 we have:

\lemma{\ver.5} {\it 
\item{(i)} For each $\e>0$ and each $\d\in\D$, the function 
$(s,x)\mapsto f^{(0)}_{\e}({s},x,\d)$ is jointly continuous in $s$ and $x$
relative to $\R^+\times \inte(\inte_{\e}\L)$. 

\item{(ii)}For each $\d\in\D$, the function 
$(s,x)\mapsto f({s},x,\d)$ is jointly continuous in $s$ and $x$
relative to $\R^+\times\inte\L$.
}

\proof It follows from $(H2)$ of Hypothesis \ver.3 that the collection of 
functions 
$\{f^{(0)}_{\e}(\,.\,,x,\d)\mid x\in\inte_\e\L, \d\in\D\}$ is equi-Lipshitzian
on $\R^+$, implying that the function
$s\mapsto f^{(0)}_{\e}({s},x,\d)$ is continuous relative to $\R^+$
for each $x\in\inte_{\e}\L$ and $\d\in\D$.
Using Lemma \ver.1, (ii), it follows from
$(H3)$ of Hypothesis \ver.3 that the collection of functions 
$\{f^{(0)}_{\e}({s},\,.\,,\d)\mid s\in\R^+, \d\in\D\}$ is equi-Lipshitzian
on all closed bounded subsets ${\cal S}\subset\inte(\inte_{\e}\L)$ and hence,
in particular,
the function $x\mapsto f^{(0)}_{\e}({s},x,\d)$ is continuous relative to 
$\inte(\inte_{\e}\L)$ for each $s\in\R^+$ and 
$\d\in\D$.
The joint continuity of $f^{(0)}_{\e}({s},x,\d)$ in $s$ and $x$
simply results from the fact that
$\R^+$ and $\inte(\inte_{\e}\L)$ are locally compact topological space. This
proves (i). In view of the remark following Hypothesis \ver.4, the proof of 
(ii) is identical to that of (i). The lemma is proven.\endproof


Each of the following functions are mapping 
$\R^+\times\R^d\times\R^d$ into $[-\infty,+\infty]$: 
$$
{\LL}(t,u,v)=\log\sum_{\d\in\D}e^{(v,\d)+f(t,u,\d)}
\Eq(2.15)
$$
$$
{\LL}^*(t,u,v^*)=
\sup_{v\in\R^d}\left\{(v,v^*)-{\LL}(t,u,v)\right\}
\Eq(2.17)
$$
$$
L_\e(t,u,v)=\log\sum_{\d\in\D}e^{(v,\d)+f_\e(t,u,\d)}
\Eq(2.17bis)
$$
$$
L_\e^*(t,u,v^*)=\sup_{v\in\R^d}\left\{(v,v^*)-L_\e(t,u,v)\right\}
\Eq(2.17ter)
$$
We set 
$$
\LL_\e^{(r)*}(t,u,v^*)\equiv \inf_{t':|t'-s|\leq r}\inf_{u':|u'-u|\leq r}
L_\e^{*}(t',u',v^*),\,\,\,\,r>0 
\Eq(2.18)
$$
Finally, we set
$$
{\LL}^{(r)*}(t,u,v^*)\equiv \inf_{t':|t'-s|\leq r}\inf_{u':|u'-u|\leq r}
{\LL}^{*}(t,u,v^*),\,\,\,\,r>0 
\Eq(2.19bis)
$$
and
$$
\ov{\LL}^*(t,u,v^*)\equiv \lim_{r\downarrow 0} {\LL}^{(r)*}(t,u,v^*)
\Eq(2.21)
$$

The main function spaces appearing in the text are listed hereafter.
All of them are spaces of  $\R^d$-valued functions on some finite
interval $[0,T]$.
By $C([0,T])$ we denote the space of continuous functions 
equipped with the supremum norm:
$\|\phi(.)\|_C=\max_{0\leq t\leq T} |\phi(t)|$, 
where $|\,.\,|$ denotes the Euclidean norm on $\R^d$ (i.e. $|x|=\sqrt{(x,x)}$).
$L^p([0,T])$, $1\leq p<\infty$, is the familiar  
space of Lebesgue measurable functions for which
$\int_{0}^{T}|\phi(t)|^pdt$ is finite and is equipped with the norm
$\|\phi(.)\|_p=\left(\int_{0}^{T}|\phi(t)|^pdt\right)^{1/p}$.
$W([0,T])$ denotes the Banach space of absolutely continuous functions and
can be equipped, e.g., with the norm,
$\|\phi(.)\|_W=|\phi(0)|+\|\dot\phi(.)\|_1$.
Recall that
$$
W([0,T])=\Bigl\{
\phi\in C([0,T])\Big|
\forall \e>0\, \exists \d>0\, s.t. \,
\sum_{l=1}^k\left|t_l-t_{l-1}\right|<\d
\Rightarrow
\sum_{l=1}^k\left|\phi(t_l)-\phi(t_{l-1})\right|<\e
\Bigr\}
\Eq(2.15bis)
$$
or, equivalently,
$$
W([0,T])=\left\{
\phi\in C([0,T])\Big|
\forall t'\in[0,T],\, \forall t\in[t',T],\,\,
\phi(t)-\phi(t')=\int_{t'}^t\dot\phi(s)ds
,\,\, \dot\phi\in L^1([0,T])
\right\}
\Eq(2.16bis)
$$
As a rule all spaces above are metrized with the norm-induced metric and
are considered in the  metric topology 
(i.e., the topology of uniform convergence).

We need to introduce some subsets of this space. Recall that the effective
domain of a an extended-real-valued function $g$ on $X$ is the set 
$\dom g\equiv\{x\in X\mid g(x)<\infty\}$.
For each $(t,u)\in\R^+\times\L$ define the extended-real-valued function
$\overline{\Phi}^*_{t,u}$ through:
$$
\overline{\Phi}^*_{t,u}(v^*)=\overline{\cal L}^*(t,u,v^*)
\Eq(2.30)
$$
Setting
$$
D_u=\dom\overline{\Phi}^*_{t,u},\,\,\,\,D=\conv\D
\Eq(2.31)
$$ 
we define,
$$
\overline{\DD}([0,T])\equiv
\left\{
\phi\in W([0,T])\Big|
\phi(t)\in\L \, \hbox{and}\,
\dot\phi(t)\in D_{\phi(t)}\,\, \hbox{for Lebesgue a.e.}\,\, t\in[0,T]
\right\}
\Eq(2.33)
$$
$$
\DD^{\circ}([0,T])\equiv
\left\{
\phi\in W([0,T])\Big|
\phi(t)\in\inte\L \, \hbox{and}\,
\dot\phi(t)\in D\,\, \hbox{for Lebesgue a.e.}\,\, t\in[0,T]
\right\}
\Eq(2.34)
$$

Our prime interest will be in the large deviation behaviour of a family
of continuous time processes constructed from
the Markov chains $\{X_{\e},\e>0\}$ by linear interpolation on the 
coordinate variables and rescaling of the time. More precisely, 
let $[0,T]$ be an arbitrary but finite interval and define the process
$Y_{\e}$ on sample path space 
$\left(C([0,T]), {\cal B}(C([0,T]))\right)$ by setting, for each $t\in [0,T]$,
$$
Y_{\e}(t)=X_{\e}\left(\left[\sfrac{t}{\e}\right]\right)
+\left(\sfrac{t}{\e}-\left[\sfrac{t}{\e}\right]\right)
\left(X_{\e}\left(\left[\sfrac{t}{\e}\right]+1\right)
-X_{\e}\left(\left[\sfrac{t}{\e}\right]\right)\right)
\Eq(2.17quater)
$$
Let $\widetilde{\cal P}_{\e,\phi_0}\equiv 
{\cal P}_{\e,\phi_0}\circ Y_{\e}^{-1}$ denote it's law. 
We are now in a position to state our main result.

\theo{\ver.6} {\it Assume that the Hypothesis \ver.2,  \ver.3, \ver.4
are satisfied. If moreover
\item{(H4)} For any convex set $\AA\subset W([0,T])$ 
$$
\inf_{\phi\subset\AA\cap \ov\DD([0,T])}\displaystyle{\int_0^T 
\overline{\cal L}^*(t,\phi(t),\dot\phi(t))dt}
=\inf_{\phi\subset\AA\cap \DD^\circ([0,T])}\displaystyle{\int_0^T 
\overline{\cal L}^*(t,\phi(t),\dot\phi(t))dt}
\Eq(2.36)
$$
then the family of measures 
$\left\{\widetilde{\cal P}_{\e,\phi_0}, \e>0\right\}$ on 
$\left(C([0,T]), {\cal B}(C([0,T]))\right)$ obeys a full large 
deviation principle with good rate function 
${\cal I}:C([0,T])\rightarrow\R^+$ 
given by
$$
{\cal I}(\phi(.))=
\cases
\displaystyle{\int_0^T 
\overline{\cal L}^*(t,\phi(t),\dot\phi(t))dt}
& \hbox{if}\,\, \phi(.)\in\overline\DD([0,T])\,\hbox{ and}\,\,\phi(0)=\phi_0
\cr
&\cr
+\infty
& \hbox{otherwise}\cr
\endcases
\Eq(2.19)
$$
}
\proposition {\ver.7} {\it 
Condition (H4) is satisfied if the following two  conditions hold:
\item{(i)} 
At each $(t,u,v^*)\in\R^+\times\L\times\R^d$ 
$$
\lim_{i\rightarrow\infty}\overline{\LL}^*(t,u_i,v^*)\leq 
\overline{\LL}^*(t,u,v^*)
\Eq(2.22)
$$
for every sequence $u_1,u_2,\dots$ in $\inte\L$ converging to
$u\in\L$. 
\item{(ii)}  For some function 
$g:\R^+\rightarrow \R^+$ satisfying
$\lim_{\a\downarrow 0} \a g(\a)=0$, for all 
$(s,u,v^*)\in\R^+\times\inte \L\times D$, 
$$
\overline\LL^*(s,u,v^*)\leq g\left(\dist(u,\L^c)\right)
\Eq(2.37)
$$
}

\remark Since $\ov\LL^*\leq \LL^*$, it is of course enough to verify 
\eqv(2.37) for the more explicitly given function $\LL^*$. This 
condition is realized in most examples of interest. Condition 
(H4) is of course always realized in situations where the process cannot 
reach the boundary of $\L$ in finite time, and in particular if $\L=\R^d$.

Proposition \ver.7 will be proven in Section 4.
 
For later reference the properties of ${\cal I}$ are given explicitly in the 
proposition below. 

\proposition{\ver.8} {\it The function ${\cal I}$ defined in
\eqv(2.19) verifies:
\item{(i)} $0\leq{\cal I}(\phi(.))\leq\infty$ and $\dom{\cal I}=\ov\DD([0,T])$
\item{(ii)} ${\cal I}(\phi(.))$ is lower semi continuous.
\item{(iii)} For each $l<\infty$, the set 
$\{\phi(.)\mid{\cal I}(\phi(.))<l\}$ is compact in $C([0,T])$.
}

\proof The proof of this proposition is in fact a more or less 
identical rerun of the proof given Section 9.1 of  Ioffe and Tihomirov [IT]
and we will not repeat it here. \endproof

\noindent By definition (i) and (ii) are the standard properties of a
rate function while  goodness is imparted to it by property $(iii)$.


\remark The LDP of Theorem \ver.6 can easily be extended beyond the continuous
setting arising from the definition of $Y_{\e}$ in that, instead of 
$Y_{\e}$, we could consider the process $Z_{\e}$ defined by,  
$$
Z_{\e}(t)=X_{\e}\left(\left[\sfrac{t}{\e}\right]\right),\,\,\,\,
\hbox{for each}\,\,
t\in [0,T]
\Eq(2.20)
$$
Naturally the path space of $Z_{\e}$ is now the space $D([0,T])$ of 
functions that are right 
continuous and have left limits which, equipped with the Skorohod topology,
${\cal S}$, is rendered Polish (we refer to the beautiful book by [Bi] for 
questions related to this space). It can then be shown that the family of 
measures 
$\left\{\widehat{\cal P}_{\e,\phi_0}, \e>0\right\}$ on 
$\left(D([0,T]),{\cal S}\right)$ obeys a full large 
deviation principle with good rate function ${\cal I'}$ where
${\cal I'}={\cal I}$ on $C([0,T])$ and ${\cal I'}=\infty$ on 
$D([0,T])\setminus C([0,T])$. The basic step needed to extend the LDP of
Theorem \ver.6 to the present case is to establish that the measures 
$\widetilde{\cal P}_{\e,\phi_0}$ and $\widehat{\cal P}_{\e,\phi_0}$, both
defined on $\left(D([0,T]),{\cal S}\right)$, are exponentially equivalent. 
As will become clear in the next chapter (see Lemma 3.1), this property is 
very easily seen to hold.

Let us finally make some remarks on the large deviation principle we have 
obtained. The rate function \eqv(2.19) has the form of a classical action 
functional with $\ov\L^*(t,x,v)$ being a (in general time dependent) 
Lagrangian. Note however that in contrast to the setting of classical 
mechanics, the function space is one of absolutely continuous function,
rather than functions with absolutely continuous derivatives. Therefore the 
minimizers in the LDP need not be solutions of the corresponding 
Euler-Lagrange equations everywhere, but jumps between solutions can occur.
A particular feature, that is due to the discrete-time nature of the process
is the presence of a maximal velocity (i.e. a ``speed of light''), due to the
 fact that  the Lagrangian is infinite for $v\not \in D$. In that respect 
one can consider  the rate function as the action of a relativistic classical
mechanics.

\newpage



\chap{3. The basic large deviation estimates.}3

The aim of this short chapter is to bring into focus the basic large deviation 
estimates on which the proof of Theorem 2.6 relies. These estimates are 
established in a subset of the continuous paths space, chosen in
such a way as to retain the underlying geometrical properties of the paths of 
$Y$. Assuming these estimates we then proceed to give the proof of Theorem 2.6.

More precisely set:
$$
{\cal E}([0,T])=\Bigl\{
\phi\in C([0,T])\Big|\,
\frac{\phi(t)-\phi(t')}{t-t'}\in D\,\,
\forall t\in[0,T],\, \forall t'\in[0,T], t\neq t'
\Bigr\}
\Eq(3.1)
$$

\lemma{\ver.1} {\it 
$\wt\PP_{\e,\phi_0}({\cal E}([0,T]))=1$ for all $\e>0$.} 

\proof Assume that $t>t'$ and set $t=(i+\g)\e$, $t'=(j+\g')\e$ where
$i,j\in\N$, $\g,\g'\in[0,1)$. By \eqv(2.17bis),
$$
\frac{Y(t)-Y(t')}{t-t'}
=\frac{X(i)-X(j)+\g[X(i+1)-X(i)]-\g'[X(j+1)-X(j)]}{[(i+\g)-(j+\g')]\e}
\Eq(3.2)
$$
Using that all sample paths of $X$ have increments of the form 
$X(k+1)-X(k)=\e\d_{k+1}$ with $\d_{k}\in\D$, \eqv(3.2) yields
$$
\frac{Y(t)-Y(t')}{t-t'}
=\cases
\d_{i+1} & \hbox{if}\,\, i=j\cr
&\cr
\displaystyle
\frac{(1-\g')\d_{j+1}+ 
\left(\sum_{k=j+2}^{i}\d_k\1_{\{i>j+1\}}\right)
+\g\d_{i+1}}{(1-\g')+(i-j-1)+\g}
& \hbox{if}\,\, i\geq j+1
\endcases
\Eq(3.3)
$$
The last line in the r.h.s. of \eqv(3.3) is a convex combination of elements 
of $\D$. Thus, remembering that $D=\conv\D$, the lemma is proven.
\endproof

Being a subset of a metric space, ${\cal E}([0,T])$ is itself a metric space
with metric given by the supremum norm-derived metric, and thus, can be 
considered a topological space in it's own right in the metric topology. 
In addition, it inherits the topology induced by $C([0,T])$. But those two
topologies are easily seen to coincide, i.e.,
$\BB({\cal E}([0,T]))=\left\{{\cal A}\cap{\cal E}([0,T]) : 
{\cal A}\in\BB(C([0,T])\right\}$.
From this and Lemma \eqv(3.1) it follows that\break 
$\left({\cal E}([0,T]), \BB({\cal E}([0,T])), \wt\PP_{\e,\phi_0}\right)$
is a measure space.

Let $\BB_{\rho}(\phi)\in {\cal E}([0,T])$ denote the open ball of radius $\rho$
around $\phi$ and let $\bar\BB_{\rho}(\phi)$ be it's closure.
Our first result will be a pair of upper and lower bounds that hold under 
much weaker hypothesis than those of Theorem 2.6.

\proposition{\ver.2} {\it Assume that Hypothesis 2.2, 2.3 and 2.4 hold. Let  
$\left\{\wt\PP_{\e,\phi_0}, \e>0\right\}$ be defined on
$\left({\cal E}([0,T]), {\cal B}({\cal E}([0,T]))\right)$. 
Then, for any $\rho> 0$ and $\phi\in{\cal E}([0,T])$, 
$$
\limsup_{\e\rightarrow 0}\e\log\wt\PP_{\e,\phi_0}(\bar\BB_{\rho}(\phi))
\leq -\inf_{{\psi\in \bar\BB_{\rho}(\phi)\cap\ov\DD([0,T]):}\atop{
\psi(0)=\phi_0}}
\int_0^T\,dt \ov\LL^*(t,\psi(t),\dot\psi(t))
\Eq(3.4bis)
$$
$$
\liminf_{\e\rightarrow 0}\e\log\wt\PP_{\e,\phi_0}(\BB_{\rho}(\phi))
\geq -\inf_{{\psi\in \BB_{\rho}(\phi)\cap\DD^\circ([0,T]):}\atop{
\psi(0)=\phi_0}}
\int_0^T\,dt \LL^*(t,\psi(t),\dot\psi(t))
\Eq(3.5bis)
$$
}

In Section 4 we will prove the following lemma:

\lemma{\ver.3} {\it Under the same hypothesis as Proposition \ver.2,
for
all $\psi\in\DD^\circ([0,T])$, 
$$
\int_0^T\,dt \LL^*(t,\psi(t),\dot\psi(t))=\int_0^T\,dt \ov\LL^*(t,\psi(t),\dot\psi(t))
\Eq(3.10)
$$
}

This lemma together with hypothesis (H4) will in fact imply the stronger

\proposition{\ver.4} {\it If in  addition to the assumptions of Proposition 
\ver.2  condition (H4) is satisfied.  
Then, for any $\rho> 0$ and $\phi\in{\cal E}([0,T])$, 
$$
\limsup_{\e\rightarrow 0}\e\log\wt\PP_{\e,\phi_0}(\bar\BB_{\rho}(\phi))
\leq -\inf_{\psi\in \bar\BB_{\rho}(\phi)}{\cal J}(\psi)
\Eq(3.4)
$$
$$
\liminf_{\e\rightarrow 0}\e\log\wt\PP_{\e,\phi_0}(\BB_{\rho}(\phi))
\geq -\inf_{\psi\in \BB_{\rho}(\phi)}{\cal J}(\psi)
\Eq(3.5)
$$
where
${\cal J}:{\cal E}([0,T])\ni\psi\mapsto {\cal J}(\psi)\equiv{\cal I}(\psi)$ 
is the restriction of ${\cal I}$ to ${\cal E}([0,T])$.
}

\proof We prove the proposition assuming Proposition \ver.2 and Lemma \ver.3.
Using first  (H4) and then (ii) of Lemma \ver.3, we get 
$$
\eqalign{
\inf_{\psi\in \bar\BB_{\rho}(\phi)\cap\ov\DD([0,T])}
\int_0^T\,dt \ov\LL^*(t,\psi(t),\dot\psi(t))
&=\inf_{\psi\in \bar\BB_{\rho}(\phi)\cap\DD^\circ([0,T])}
\int_0^T\,dt \ov\LL^*(t,\psi(t),\dot\psi(t))\cr
=\inf_{\psi\in \bar\BB_{\rho}(\phi)\cap\DD^\circ([0,T])}
\int_0^T\,dt \LL^*(t,\psi(t),\dot\psi(t))
&\geq \inf_{\psi\in \bar\BB_{\rho}(\phi)\cap\ov\DD([0,T])}
\int_0^T\,dt \LL^*(t,\psi(t),\dot\psi(t))\cr
}\Eq(3.21)
$$
which implies \eqv(3.4). 

On the other hand, using first (ii) of Lemma \ver.3, then (H4), and finally the
fact that, since for any $r,\e>0$, $\LL_\e^{(r)*}(t,u,v^*)\leq 
\LL_\e^*(t,u,v^*)$, we have  
$
\ov\LL(t,u,v^*)\leq 
\LL^*(t,u,v^*)$,
we also get
$$\eqalign{
\inf_{\psi\in \BB_{\rho}(\phi)\cap\DD^\circ([0,T])}
\int_0^T\,dt \LL^*(t,\psi(t),\dot\psi(t))&=
\inf_{\psi\in \BB_{\rho}(\phi)\cap\DD^\circ([0,T])}
\int_0^T\,dt \ov\LL^*(t,\psi(t),\dot\psi(t))
\cr
=\inf_{\psi\in \BB_{\rho}(\phi)\cap\ov\DD([0,T])}
\int_0^T\,dt \ov\LL^*(t,\psi(t),\dot\psi(t))
&\leq \inf_{\psi\in \BB_{\rho}(\phi)\cap\ov\DD([0,T])}
\int_0^T\,dt \LL^*(t,\psi(t),\dot\psi(t))
}
\Eq(3.22)
$$
which implies \eqv(3.5). \endproof

The proof of Theorem 2.6, assuming Proposition
\ver.4 and Proposition 2.8, is now classical.

\proofof{Theorem 2.6} 
Assume Proposition
\ver.2 and Proposition 2.3 to hold. Then, 
on the one hand, since $C([0,T])$ is Polish, goodness of 
the rate function entails exponential tightness of the
family $\left\{\wt\PP_{\e,\phi_0}, \e>0\right\}$, which implies that the 
full LDP obtains whenever it's weak version obtains. 
On the other hand, since  ${\cal E}([0,T])$ is compact, 
it follows from Proposition \ver.2 that 
the family $\left\{\wt\PP_{\e,\phi_0}, \e>0\right\}$ on ${\cal E}([0,T])$ 
obeys a weak LDP with rate function ${\cal J}$.
The connection between these LDP's is made in through:

\lemma{\ver.5}([DZ], Lemma 4.1.5) {\it Let $S$ be a regular topological space
and $\{\mu_{\e}, \e\geq 0\}$ a family of probability measures on $S$.
Let ${\cal S}$ be a measurable subset of $S$ such that $\mu_{\e}({\cal S})=1$
for all $\e>0$. Assume ${\cal S}$ equipped with the topology induced by $S$.
\item{(i)} if ${\cal S}$ is a closed subset of $S$ and $\{\mu_{\e}\}$ satisfies
the LDP in ${\cal S}$ with rate function ${\cal J}$, then 
$\{\mu_{\e}\}$ satisfies the LDP in $S$ with rate function 
${\cal I}={\cal J}$ on ${\cal S}$ and ${\cal I}=\infty$ on 
$S\setminus{\cal S}$.
\item{(ii)} If $\{\mu_{\e}\}$ satisfies the LDP in $S$ with rate function 
${\cal I}$ and $\dom{\cal I}\subset{\cal S}$, then the same LDP holds in 
${\cal S}$.
}

\remark Lemma \ver.5 holds for the weak as well as the full LDP.
 
Theorem 2.6 now follows from  Lemma \ver.5 together with  Lemma \ver.1 and
the fact that  being compact,
 ${\cal E}([0,T])$ is closed in $C([0,T])$
\endproof

\newpage


\chap{4. Convexity related results}4

This rather lengthy chapter establishes most of the basic analytic
properties of the logarithmic 
moment generating functions and their Legendre transforms
that will be needed to prove the  upper and lower large deviation estimates 
in Section 5. We begin by fixing some notations.

Let ${\LL}_{\e}$ and ${\LL}_\e^*$ be the functions, mapping
$\R^+\times\R^d\times\R^d$ into $\R$, defined by: 
$$
{\LL}_{\e}({s},u,v)
=
\log\sum_{\d\in\D}e^{(v,\d)+f_{\e}^{(0)}({s},u,\d)}
\Eq(4.1)
$$
$$
{\LL}_{\e}^*({s},u,v^*)=
\sup_{v\in\R^d}\left\{(v,v^*)-{\LL}_{\e}({s},u,v)\right\}
\Eq(4.2)
$$
It plainly follows from Hypothesis 2.2 and (H0) of Hypothesis 2.3 that on 
$\R^+\times(\R^d\setminus\L)\times\R^d$,
${\LL}_{\e}=-\infty$, ${\LL}_{\e}^*=+\infty$, $\LL\equiv -\infty$
, and
$\LL^*=+\infty$. We will thus limit our 
attention to the behaviour of these functions on $\R^+\times\L\times\R^d$.

Let $\MM(\D)$ denote the set of all probability measures on $\D$. 
The support of a measure $\nu\in\MM(\D)$, denoted $\supp\,\nu$, is defined by 
$\supp\,\nu=\{\d\in\D\mid \nu(\d)>0\}$.
For any fixed $(s,u)\in\R^+\times\L$ and any $v\in\R^d$ let
$\nu^{v}_{\e,{s},u}$ be the probability measure on 
$\MM(\D)$ that assigns to $\d\in\D$ the density
$$
\nu^{v}_{\e,{s},u}(\d)=
\frac{e^{(v,\d)+f_{\e}^{(0)}({s},u,\d)}}
{\sum_{\d\in\D}e^{(v,\d)+f_{\e}^{(0)}({s},u,\d)}}
\Eq(4.7)
$$
Similarly $\nu^{v}_{{s},u}\in\MM(\D)$ is defined by \eqv(4.7) with
$f_{\e}^{(0)}({s},u,\d)$ replaced by $f({s},u,\d)$. 

Observe that if $u\in\L$ then either $u\in\inte_{\e}\L$ or 
$u\in\bd_{\e}\L$ and, according to the remark following Hypothesis
2.2, 
$$
\supp\,\nu^{0}_{\e,{s},u}=\D,\,\,\,\,\,\,\,\,\,\,
\forall (s,u)\in\R^+\times\inte_{\e}\L
\Eq(4.11)
$$
whereas
$$
\emptyset\neq\supp\,\nu^{0}_{\e,{s},u}\subset\D,\,\,\,\,\,\,\,\,\,\,
\forall (s,u)\in\R^+\times\bd_{\e}\L
\Eq(4.11bis)
$$
Moreover, for $\chi$ a random variable with law
$\nu^{v}_{\e,{s},u}$,
$$
{\LL}_{\e}({s},u,v)=\E_{\nu^{0}_{\e,{s},u}}e^{(v,\chi)}
+\log\sum_{\d\in\D}e^{f_{\e}^{(0)}({s},u,\d)}
\Eq(4.8)
$$
where $\E_{\nu^{v}_{\e,{s},u}}$ denotes the expectation w.r.t.
$\nu^{v}_{\e,{s},u}$. Thus, up to a small term (which goes to zero as
$\e\downarrow 0$) $\LL_{\e}$ is the logarithmic moment generating function
of $\nu^{v}_{\e,{s},u}$, ${\LL}_{\e}^*$ being termed it's conjugate.

With ${\LL}$ and ${\LL}^*$ given by \eqv(2.15) and \eqv(2.17), 
for fixed  $(s,u)\in\R^+\times\L$, we further define the functions, 
mapping $\R^d$ into $\R$:
$$
\eqalign{
&\Phi_{\e,{s},u}(v)={\LL}_{\e}({s},u,v)\cr
&\Phi^*_{\e,{s},u}(v^*)={\LL}_{\e}^*({s},u,v^*)\cr
&\Phi_{{s},u}(v)={\LL}({s},u,v)\cr
&\Phi^*_{{s},u}(v)={\LL}^*({s},u,v^*)\cr
}
\Eq(4.3)
$$

This chapter is divided into five subchapters. In the first 
subchapter we establish the
properties of the functions ${\Phi}_{\e}$, ${\Phi}$, and their 
conjugates ${\Phi}^*_{\e}$, ${\Phi}^*$. Although most of them are 
well know folklore of the theory of convex analysis, it is more convenient to 
state them at once rather then laboriously recall them from the literature 
when we need to put them in use. The proofs are merely compilations of 
references, chiefly taken from the books by Rockafellar [Ro] and
Ellis [E]. In the second 
subchapter we go back to the functions ${\LL}_{\e}$, ${\LL}^*$, and their 
limits, and establish their topological properties.
The third subchapter establishes some basic properties of semi-continuous
regularisations of our functions, and in particular provides an important 
result on the uniform convergence of the regularised functions as 
$\e\downarrow 0$.   In the forth subsection we present a result,
based on these topological properties, which shall be crucial in establishing 
the large deviation bounds, while the last subsection is devoted to the 
proof of Proposition 2.7.

Most of the results of this section will be established simultaneously 
for either the function ${\LL}_{\e}$ or ${\LL}_{\e}^*$ at fixed $\e$, and 
(what we shall see are their limits) ${\LL}$ or ${\LL}^*$. We stress here once 
for all that, according to the remark following Hypothesis 2.4, all results 
for  ${\LL}_{\e}$ or ${\LL}_{\e}^*$ directly infered from Hypothesis 2.2 and
2.3 obviously carry through to the limiting functions. As a rule we 
systematically skip the proofs of results for ${\LL}$ or ${\LL}^*$ whenever 
they are simple repetitions of those for ${\LL}_{\e}$ or ${\LL}_{\e}^*$.

\vskip.3truecm
\noindent{\bf \ver.1. The functions ${\Phi}_{\e}$, ${\Phi}$, and their 
conjugates.}
\vskip.3truecm

We begin with a short reminder of terminology and a few definitions.
Recall that $\dom g\equiv\{x\in X\mid g(x)<\infty\}$. All through 
this chapter we shall adopt the usual convention that consists in identifying 
a convex function $g$ on  $\dom g$ with the convex function defined 
throughout $\R^d$ by setting $g(x)=+\infty$ for $x\notin\dom g$.
A real valued function $g$ on a convex set $C$
is said to be {\it strictly}  convex on $C$ if
$$
g((1-\l)x+\l y)<(1-\l)g(x)+\l g(y),\,\,\,\,0<\l<1
\Eq(4.4)
$$
for any two different points $x$ and $y$ in $C$. It is called {\it proper} if 
$g(x)<\infty$ for at least one $x$ and $g(x)>-\infty$ for every $x$.
The {\it closure} of a convex function $g$ is defined to be the lower 
semi-continuous hull of $g$ if $g$ nowhere has value $-\infty$, whereas
the closure of $g$ is defined to be the convex function $-\infty$ if $g$ is
an improper convex function such that $g(x)=-\infty$ for some $x$. Either way
the closure of $g$ is another convex function and is denoted $\cl\, g$.
A convex function is said to be closed if $g=\cl\, g$.
For a proper convex function closedness is thus the same as lower 
semi-continuity.
A function $g$ on $R^d$ is said to be continuous {\it relative to}
a subset $\SS$ of $\R^d$ if the restriction of $d$ to $\SS$ is a continuous 
function.

For any set $C$ in $\R^d$ we denote by $\cl\, C$, $\inte\, C$ and
by $\bd\, C=(\cl\, C)\setminus(\inte\, C)$ 
the {\it closure, interior} and {\it boundary} of $C$.  
If $C$ is convex, we denote by $\ri\, C$ and 
$\rbd\, C=(\cl\, C)\setminus(\ri\, C)$
it's {\it relative interior} and {\it relative boundary}.

\definition{\ver.1}{\it
A proper convex function $g$ on $\R^d$ is called essentially smooth if it 
satisfies the following three conditions for $C=\inte(\dom g)$:
\item{(a)} $C$ is non empty;
\item{(b)} $g$ is differentiable throughout $C$;
\item{(c)} $\lim_{i\rightarrow\infty}|\nabla g(x_i)|=+\infty$ whenever 
$x_1,x_2,\dots$, is a sequence in $C$ converging to a boundary point $x$ 
of $C$.

Note that a smooth convex function on $\R^d$ is in particular essentially 
smooth (since (c) holds vacuously).
}

\definition{\ver.2}{\it The conjugate $g^*$ of an arbitrary function $g$ 
on $\R^d$ is defined by
$$
g^*(x^*)=\sup_{x\in\R^d}\left\{(x,x^*)-g(x)\right\}
\Eq(4.5)
$$
}

Note that both $\Phi_{\e}$,  $\Phi_{\e}^*$ and $\Phi$,  $\Phi^*$ are pairs 
of conjugate functions. 

\lemma{\ver.3}{\it ([Ro], Theorem 12.2)
Let $g$ be a convex function. The conjugate function $g^*$ is then a closed 
convex function, proper if and only if $g$ is proper. Moreover
$(\cl\, g)^*=g^*$ and $g^{**}=\cl\,g$.
}

Finally, for $g$ an extended-real-valued function on $\R^d$ which is
is finite and twice differentiable throughout $\R^d$, we denote by 
$
\nabla g(x)\equiv\left(
\sfrac{\partial g}{\partial x_1}(x),\dots,\sfrac{\partial g}{\partial x_d}(x)
\right)
$,
$
\nabla^2g(x)\equiv\left(
\sfrac{\partial g}{\partial x_i\partial x_j}(x)\right)_{i,j=1,\dots,d}
$
and 
$\D g(x)=\sum_{i=1}^d\sfrac{\partial^2 g}{\partial^2 x_i}(x)$, respectively,
the gradient, the Hessian, and the Laplacian of $g$ at $x$.

In order to unburden formulas the indices ${s}$ and $u$ in
\eqv(4.3) and \eqv(4.7) will systematically be dropped in the sequel. 
We start by listing some of the properties of $\Phi_{\e}$ and $\Phi$.

\lemma{\ver.4}{\it For all $\e>0$ the following conclusions hold.
For any fixed $(s,u)\in\R^+\times\L$,
\item{(i)} $|\Phi_{\e}(v)|<\infty$ for all $v\in\R^d$.
\item{(ii)} $\Phi_{\e}$ is a closed, convex, and 
continuous function on $\R^d$.
\item{(iii)} $\Phi_{\e}$ has mixed partial derivatives of all order which 
can be calculated by differentiation under the sum sign. In particular, for
all $v\in\R^d$,
if $\chi=(\chi_1,\dots,\chi_d)$ denotes a random vector with law
$\nu^{v}_{\e,{s},u}$,
$$
\nabla\Phi_{\e}(v)=\E_{\nu^{v}_{\e}}\chi
=\left(\E_{\nu^{v}_{\e}}\chi_i\right)_{i=1,\dots,d}
\Eq(4.12)
$$
$$
\nabla^2\Phi_{\e}(v)=\left(\E_{\nu^{v}_{\e}}\chi_i\chi_j
-\E_{\nu^{v}_{\e}}\chi_i\E_{\nu^{v}_{\e}}\chi_j\right)_{i,j=1,\dots,d}
\Eq(4.13)
$$
$$
\D\Phi_{\e}(v)
=\E_{\nu^{v}_{\e}}|\chi-\E_{\nu^{v}_{\e}}\chi|^2
=\sum_{i=1}^d\E_{\nu^{v}_{\e}}|\chi_i-\E_{\nu^{v}_{\e}}\chi_i|^2
\Eq(4.14)
$$

Moreover, for any fixed $(s,u)\in\R^+\times\inte_{\e}\L$, $\Phi_{\e}$
is a strictly convex function on $\R^d$.

All assertions above hold with $\Phi_{\e}$ replaced by $\Phi$
and $\nu^{v}_{\e}$ replaced by $\nu^{v}$.
}

\proof  If $u\in\L$ then, by Hypothesis 2.2,
$$
\textstyle\left|\log\sum_{\d\in\D}e^{f_{\e}^{(0)}({s},u,\d)}\right|<\infty
\Eq(4.15)
$$
Assertion (i) is then a consequence of \eqv(4.8).
Given assertion (i), 
assertions (ii) and (iii) are proven, e.g., in [E] (see  
pp230 for the former and Theorem VII.5.1 for the latter); 
formulae \eqv(4.12), \eqv(4.13) and \eqv(4.14) 
may be found in [BG]. Finally, a necessary and 
sufficient condition for $\Phi_{\e}$ to be strictly convex 
(see e.g. [E], Proposition VIII.4.2) is that the affine hull of
$\supp\,\nu^{0}_{\e}$ coincides with $\R^d$; but by Hypothesis 2.1 
this condition is fulfilled whenever $u\in\inte_{\e}\L$.
The lemma is proven.
\endproof

We next turn to the functions $\Phi_{\e}^*$ and $\Phi^*$. We first
 state an important relationship between the support of $\nu^{0}_{\e}$
and the effective their effective domains.

\lemma{\ver.5}{\it 
Let $d\geq 1$,  $\e>0$ and $(s,u)\in\R^+\times\L$. Then,
$$
\dom\Phi_{\e}^*=\conv(\supp\,\nu^{0}_{\e})
\Eq(4.16)
$$
In particular, if $(s,u)\in\R^+\times\inte_{\e}\L$,
$$
\dom\Phi_{\e}^*=\conv\,\D
\Eq(4.17)
$$
The same holds with $\Phi_{\e}$ replaced by $\Phi$
and $\inte_\e\L$ replaced by $\inte\L$.
}

\remark Since $\supp\, \nu^0_{\e,s,u}=\left\{\d\in\D\,\big|
\,f^{(0)}_\e(s,u,\d)>-\infty\right\}$, we have by the second remark following 
Hypothesis 2.2 and (H0) that  $\exists \e'=\e'(u)>0$ s.t.  
$\forall 0<\e\leq \e'$ 
$$
\supp\, \nu^0_{\e,s,u}=\left\{\d\in\D\,\big|\L^{(\d)}\ni u\right\}
\Eq(4.900)
$$
and therefore 
$$
\dom \Phi^*_{\e,s,u}=\dom\Phi^*_{s,u}
\Eq(4.902)
$$

\proof Obviously, if $\nu^{0}_{\e}$ is the unit mass at $\d^*$,
$\Phi_{\e}^*(v^*)=0$ if $v^*=\d^*$ whereas $\Phi_{\e}^*(v^*)=+\infty$ if 
$v^*\neq\d^*$, so that \eqv(4.16) and \eqv(4.17) hold true.
Assume now that $\nu^{0}_{\e}$ is non degenerate. The
starting point to prove the lemma under this assumption is a theorem by Ellis 
([E], Theorem VIII.4.3)
which, rephrased in our setting and putting 
$S\equiv\conv(\supp\,\nu^{0}_{\e})$, states that,
$$
\dom\Phi_{\e}^*\subseteq S \,\,\text{and} \inte(\dom\Phi_{\e}^*)=\inte\,S
\Eq(4.17bis)
$$
 From this \eqv(4.16) automatically follows if we can show that 
$\Phi_{\e}^*(v^*)<\infty$ for 
$v^*\in\bd\,S$.
The proof is built upon the fact 
that, since $\supp\,\nu^{0}_{\e}\subseteq\D$, the set $S$ is  
a polytope and hence is closed.
Let $\{a_1,\dots,a_{\kappa}\}$ be the subset of $\D$ generating 
$S$ that is, the smallest subset of $\D$
such that $\conv(\{a_1,\dots,a_{\kappa}\})=S$. 
Set $\kappa\equiv|\supp\,\nu^{0}_{\e}|$. By assumption
$\nu^{0}_{\e}$ is non degenerate so that $\kappa>1$.
All points 
$v$ of $\bd\,S$ can then be expressed in the form 
$v^*=\sum_{i=1}^{\kappa}\l_ia_i$ where
$\sum_{i=1}^{\kappa}\l_i=1$, $\l_i\geq 0$, the number of non zero
coefficients $\l_i$ being at most $\kappa-1$.

We now introduce a representation of
$\Phi_{\e}^*$ due to Donsker and Varadhan ([DV], p. 425). For 
$\mu\in\MM(\D)$ define the relative entropy of $\mu$ with respect to 
$\nu_{\e}^{0}$ by
$$
I(\mu)=\sum_{\d\in\D}
\mu(\d)\log\left(\sfrac{\mu(\d)}{\nu_{\e}^{0}(\d)}\right)
\Eq(4.18)
$$
Then
$$
\Phi_{\e}^*(v^*)=
\inf\left\{I(\mu)
\Big| \mu\in\MM(\D), 
{\textstyle\sum_{\d\in\D}}\d\mu(\d)=v^*\right\}
-{\textstyle\log\sum_{\d\in\D}e^{f_{\e}^{(0)}({s},u,\d)}}
\Eq(4.19)
$$
First, observe that for $v=a\in\{a_1,\dots,a_{\kappa}\}$ the 
set $\left\{\mu\in\MM(\D), 
{\textstyle\sum_{\d\in\D}}\d\mu(\d)=a\right\}$ reduces to the unit mass
at $a$, and, by 
\eqv(4.19) and \eqv(4.15), 
$$
I(\mu)=-\log(\nu_{\e}^{0}(a))
-{\textstyle\log\sum_{\d\in\D}e^{f_{\e}^{(0)}({s},u,\d)}}
<\infty
\Eq(4.20)
$$
Next, by Lemma \ver.3, $\Phi_{\e}^*$ is convex so
that
$$
\Phi_{\e}^*\biggl(
\sum_{i=1}^{\kappa}\l_ia_i
\biggr)
\leq\sum_{i=1}^{\kappa}\l_i\Phi_{\e}^*(a_i)
<\infty
\Eq(4.21)
$$
proving that $\bd\,S\subset\dom\Phi_{\e}^*$. The lemma is proven. \endproof

We now list some of the properties of $\Phi^*_{\e}$ and $\Phi^*$. 

\lemma{\ver.6}  {\it
For all $\e>0$ the following conclusions hold.
For any fixed $(s,u)\in\R^+\times\L$,
\item{(i)}  $\Phi^*_{\e}$ is a closed convex function on $\R^d$.
\item{(ii)} $\Phi^*_{\e}$ has compact level sets.
\item{(iii)} Let $v^*_0= \E_{\nu^{v}_{\e}}\chi|_{v=0}$.
Then for any $v^*\in\R^d$, $\Phi^*_{\e}(v^*)\geq 0$
and $\Phi^*_{\e}(v^*)= 0$ if and only if $v^*=v^*_0$.
\item{(iv)} For $d=1$, $\Phi^*_{\e}$ is strictly convex and for $d\geq 2$,
$\Phi^*_{\e}$ is strictly convex on $\ri(\dom\, \Phi^*_{\e})$.
\item{(v)} $\Phi_{\e}^*$ is continuous relative to $\dom\Phi_{\e}^*$.

Moreover, for any $(s,u)\in\R^+\times\inte_\e\L$,  $\Phi^*_{\e}$ is 
essentially smooth.

All assertions above hold with $\Phi_{\e}$ replaced by $\Phi$
and $\nu^{v}_{\e}$ replaced by $\nu^{v}$.
}

\proof Assertions (i) to (iv) are taken from  [E], Theorem VII.5.5.
Since by Lemma \ver.6 $\Phi^*_{\e}$ is closed, and since by Lemma \ver.5 
$\dom\Phi^*_{\e}$ is a polytope, then (v) is a special case of [Ro], 
Theorem 10.2.
Finally, the essential smoothness of  $\Phi^*_{\e}$ follows from
the fact that, by Lemma \ver.4 , $\Phi_{\e}$ is strictly convex for
$(s,u)\in\R^+\times\inte_\e\L$ together with Theorem 26.3 of [Ro], implying
that the conjugate of a proper and  strictly convex function having
effective domain $\R^d$ is essentially smooth. \endproof

The following lemma finally relates 
the functions $\Phi_{\e}$ and $\Phi$ to their conjugates.

\lemma{\ver.7}{\it  Let $(s,u)\in\R^+\times\L$, $\e>0$. 
 For any $v\in\R^d$, the following three conditions on $v^*$ are
equivalent to each other:
\item{(i)} $v^*=\nabla\Phi_{\e}(v)$;
\item{(ii)} $(v',v^*)-\Phi_{\e}(v')$ achieves it's supremum in $v'$ at 
$v'=v$;
\item{(iii)} $(v,v^*)-\Phi_{\e}(v)=\Phi^*_{\e}(v^*)$.

If $(s,u)\in\R^+\times\inte_{\e}\L$, two more conditions can be added to 
this list;
\item{(iv)} $v=\nabla\Phi^*_{\e}(v^*)$;
\item{(v)} $(v,v')-\Phi^*_{\e}(v')$ achieves it's supremum in $v'$ at 
$v'=v^*$.
  
The same holds when $\Phi_{\e}$ and $\Phi^*_{\e}$ are replaced by $\Phi$
and $\Phi^*$.
}

\proof
By lemma \ver.4 and the definition of essential smoothness,
$\Phi_{\e}$ and $\Phi$ are closed, proper, convex, essentially smooth functions
and are differentiable throughout
$\R^d$. By Lemma \ver.5 and Lemma \ver.6, for each 
$(s,u)\in\R^+\times\inte_{\e}\L$,
$\Phi^*_{\e}$ and $\Phi^*$ are closed, proper, convex,
essentially smooth functions with effective domain $\conv\D$; hence they
are differentiable on $\inte(\conv\D)$. Since for a closed, proper, convex, and
essentially smooth function $g$ on $\R^d$, the subgradient 
of $g$ at $x$, denoted by $\partial g(x)$, reduces to the gradient 
mapping $\nabla g(x)$
\note{that is, $\partial g(x)$ consists of the vector $\nabla g(x)$
alone when $x\in\inte(\dom\, g)$, while $\partial g(x)=\emptyset$ when
$x\notin\inte(\dom\, g)$.}
(see [Ro], Theorem 26.1), then Lemma \ver.5 is a 
special case of Theorem 23.5 of [Ro].
\endproof

\vskip.3truecm
\noindent{\bf \ver.2. Topological properties of the functions ${\LL}_{\e}$, 
${\LL}^*_{\e}$, and their limits.}
\vskip.3truecm

We have so far gathered information on the collections of convex functions
$v\mapsto{\LL}_{\e}({s},u,v)$, $v\mapsto{\LL}^*_{\e}({s},u,v^*)$, and their 
limits for $s\in\R^+$ and $u$ in either $\L$, $\inte_{\e}\L$ or $\inte\L$. We
saw in particular that ${\LL}_{\e}$ (respectively ${\LL}$) is continuous in
$v$ throughout $\R^d$ and that if $u\in\inte_{\e}\L$ 
(respectively $u\in\inte\L$) then ${\LL}^*_{\e}$ (respectively ${\LL}^*$) is 
continuous in $v^*$ relative to $\conv\D$. In order to complete this picture
we shall devote this subchapter to establishing the continuity properties of 
these functions in the variables $t$ and $u$.

\lemma{\ver.8}{\it For all $\e>0$,
\item{(i)} There exists a constant $0<\theta<\infty$ such that:
$$
\sup_{\scriptstyle{u\in\L}\atop\scriptstyle{v\in\R^d}}
\left|{\LL}_{\e}({s},u,v)-{\LL}_{\e}({s'},u,v)\right|
\leq\theta|{s}-{s}'|
,\,\,\,\,\,\,\,\,\forall {s}\in\R^+, \forall {s}'\in\R^+
\Eq(4.22)
$$
\item{(ii)} For any closed bounded subset ${\cal S}\subset\inte_{\e}\L$, there 
exists a positive constant $\vartheta\equiv \vartheta({\cal S})<\infty$ such 
that:
$$
\sup_{\scriptstyle{{s}\in\R^+}\atop\scriptstyle{v\in\R^d}}
\left|{\LL}_{\e}({s},u,v)-{\LL}_{\e}({s},u',v)\right|
\leq\vartheta|u-u'|
,\,\,\,\,\,\,\,\,\forall u\in{\cal S}, \forall u'\in{\cal S}
\Eq(4.23)
$$
\item{(iii)} The function 
${\LL}_{\e}({s},u,v)$ is jointly continuous in $s$, $u$ and $v$
relative to $\R^+\times \inte(\inte_{\e}\L)\times\R^d$.
 
Assertions (i)-(iii) hold with ${\LL}_{\e}$ replaced by  ${\LL}$
and $\inte_{\e}\L$ replaced by $\inte{\L}$.

\noindent In addition:
\item{(iv)}
For any $u\in\L$, $s\in\R^+$, the function 
${\LL}_{\e}({s},u,\cdot)$ converges uniformly to 
$ {\LL}({s},u,\cdot)$ on $\R^d$.
\item{(v)} For any closed bounded  $\SS\subset\inte{\L}$, 
    ${\LL}_{\e}$ converges uniformly to ${\LL}$ on $\R^+\times \SS\times\R^d$.
}

\proof By Lemma \ver.4, both
${\LL}_{\e}$ and ${\LL}$ are finite on
$\R^+\times\L\times\R^d$. Using Hypothesis 2.2 and (H2) of 
Hypothesis 2.3 we may write, 
for any $s\in\R^+$, $s'\in\R^+$, and any $(u,v)\in\L\times\R^d$,
$$
\left|{\LL}_{\e}({s},u,v)-{\LL}_{\e}({s'},u,v)\right|
\leq
\sup_{\scriptstyle{\d\in\D:}\atop\scriptstyle{\L^{(\d,\e)}\ni u}}
|f^{(0)}_{\e}({s},u,\d)-f^{(0)}_{\e}({s'},u,\d)|
\leq \theta|{s}-{s}'|
\Eq(4.25)
$$
This proves (i). Assertions (ii) and (iv) are likewise deduced from
$(H3)$ of Hypothesis 2.3 and Hypothesis 2.4. Knowing (i), (ii), and (ii) of
Lemma \ver.4, the proof of assertion (iii) is similar to that of Lemma 2.5. 
Assertion (iv) is an immediate consequence of Hypothesis (H4).

To prove (iv), by the second remark following Hypothesis 2.2, for any 
$(s,u)\in\R^+\times \L$,  there exists $\e'=\e'(u)>0$ such that for all 
$\e\leq\e'$ such that 
$$
{\LL}_{\e}({s},u,v)
=
\log\sum_{\d\in\D : \L^{(\d)}\ni u}
e^{(v,\d)+f_{\e}^{(0)}({s},u,\d)}
\Eq(4.101)
$$
This implies that 
$$
\left|\LL_\e(s,u,v)-\LL(s,u,v)\right|
 \leq\sup_{\d\in\D : \L^{(\d)}\ni u}|f_{\e}^{(0)}({s},u,\d)-f({s},u,\d)|
\Eq(4.100)
$$
where the right hand side is independent of $v$ and, by Hypothesis 2.4,
converges to zero. This yields (iv).

Finally, the prove of (v) is almost identical to that of (iv). We only need to 
observe that the $\e'(u)$ can be chosen uniform for $u\in \SS$ if $\SS$ is a 
compact subset of the interior of $\L$, and that as indicated in the 
remark following Hypothesis 2.4, the right hand side of \eqv(4.100)
converges to zero uniformly on $\R^+\times \SS$. \endproof

\lemma{\ver.9}{\it 
For all $\e>0$,
\item{(i)} There exists a constant $0<\theta<\infty$ such that:
$$
\sup_{\scriptstyle{u\in\L}\atop\scriptstyle{v^*\in\conv\,\D}}
\left|{\LL}_{\e}^*({s},u,v^*)-{\LL}_{\e}^*({s'},u,v^*)\right|
\leq\theta|{s}-{s}'|
,\,\,\,\,\,\,\,\,\forall {s}\in\R^+, \forall {s}'\in\R^+
\Eq(4.26)
$$
\item{(ii)} For any closed bounded subset ${\cal S}\subset\inte_{\e}\L$, there 
exists a positive constant $\vartheta\equiv \vartheta({\cal S})<\infty$ such 
that:
$$
\sup_{\scriptstyle{{s}\in\R^+}\atop\scriptstyle{v^*\in\conv\,\D}}
\left|{\LL}_{\e}^*({s},u,v^*)-{\LL}_{\e}^*({s},u',v^*)\right|
\leq\vartheta|u-u'|
,\,\,\,\,\,\,\,\,\forall u\in{\cal S}, \forall u'\in{\cal S}
\Eq(4.27)
$$
\item{(iii)} The function 
${\LL}^*_{\e}({s},u,v^*)$ is jointly continuous in $s$, $u$ and $v^*$
relative to $\R^+\times \inte(\inte_{\e}\L)\times\conv\D$. 

Moreover (i)-(iii) hold with ${\LL}_{\e}^*$ replaced by  
${\LL}^*$
and $\inte_{\e}\L$ replaced by $\inte\L$.

\noindent In addition:
\item{(iv)}
For each $({s},u,v^*)\in\R^+\times\L\times\conv\,\D$,
$$
\lim_{\e\rightarrow 0}{\LL}_{\e}^*({s},u,v^*)={\LL}^*({s},u,v^*)
\Eq(4.28)
$$
exists and is finite for all  $({s},u,v^*)$ such that 
$s\in\R^+$, $u\in\L$, $v^*\in\dom \Phi_{s,u}^*$.
\item{(v)} For every closed bounded set $\SS\subset \inte\L$, 
${\LL}_{\e}^*$ converges uniformly to ${\LL}^*$ on 
$\R^+\times \SS\times\conv\,\D$.
}

\proof By Lemma \ver.5, both ${\LL}_{\e}^*$ and ${\LL}^*$ 
are finite on $\R^+\times\inte\L\times\conv\,\D$. 
To prove (i) we write that
for any $s\in\R^+$, $s'\in\R^+$, and any $(u,v^*)\in\L\times\conv\D$,
$$
\eqalign{
{\LL}_{\e}^*({s},u,v^*)
&=\sup_{v\in\R^d}\left\{(v,v^*)-{\LL}_{\e}({s'},u,v)
  +({\LL}_{\e}({s'},u,v)-{\LL}_{\e}({s},u,v))\right\}\cr
&\leq \sup_{v\in\R^d}\left\{(v,v^*)-{\LL}_{\e}({s'},u,v)
  +\sup_{v\in\R^d}\left|{\LL}_{\e}({s'},u,v)-{\LL}_{\e}({s},u,v))\right|
\right\}\cr
&= {\LL}_{\e}^*({s}',u,v^*)+
\sup_{v\in\R^d}\left|{\LL}_{\e}({s'},u,v)-{\LL}_{\e}({s},u,v))\right|
\cr
}
\Eq(4.29)
$$
and by \eqv(4.22) of Lemma \ver.8, 
$$
{\LL}_{\e}^*({s},u,v^*)-{\LL}_{\e}^*({s}',u,v^*)
\leq \theta|{s}-{s}'|
\Eq(4.30)
$$
Similarly we can show that
$$
{\LL}_{\e}^*({s},u,v^*)-{\LL}_{\e}^*({s}',u,v^*)
\geq -\theta|{s}-{s}'|
\Eq(4.31)
$$
Thus (i) is proven. Assertions (ii) is  obtained in the 
same way
on the basis of assertion (ii) of Lemma \ver.8.
whereas (iii) is deduced from Lemma \ver.8, (iii), together with Lemma \ver.6,
(v).
To prove (iv), note that using the remark following Lemma \ver.5,
there exists $\e'=\e'(u)>0$ such that 
for $\e<\e'(u)$, for any $v^*\in\dom \Phi^*_{s,u}$
$$
\left|\LL^*_\e(s,u,v^*)-\LL^*(s,u,v^*)\right|\leq \sup_{v\in\R^d}
\left|\LL_\e(s,u,v)-\LL(s,u,v)\right|
\Eq(4.31bis)
$$
and the right hand side converges to zero by  Lemma \ver.8 (iv). Note that the 
convergence is even uniform in $v^*$. 
(v) now follows by the same arguments that were used in the proof of
(v) of Lemma \ver.8.
The proof is done. \endproof

\vskip.3truecm
\noindent{\bf \ver.3. Some properties of semi-continuous regularisations.}
\vskip.3truecm

The results established in the previous sub-section will be mainly used for 
the lower bounds. For these the use of the functions $\LL_\e$, $\LL^*_\e$,
defined in terms of the functions $f^{(0)}_\e$ will be convenient.
 The upper bounds will rely on the use of (upper-, resp. 
lower)
semi-continuous regularisations of the functions $L_\e$, resp. $L^*_\e$. 
Let us first note that all results of in \ver.2  that did 
not rely to the Lipshitz continuity of $f^{(0)}_\e$ are also valid for 
$L_\e$ and $L^*_\e$.

For $r>0$ we define:
$$
\LL_\e^{(r)}(s,u,v)\equiv \sup_{s':|s-s'|\leq r} \sup_{u':|u-u'|\leq r}
L_\e(s',u',v)
\Eq(43.1)
$$
Set $\L(r)\equiv\{u\in\R^d\,|\,\dist(u,\L)\leq r\}$. 
The following lemma establishes some simple properties of $\LL_\e^{(r)}$ 
we will 
need later.

\lemma{\ver.10} {\it 
\item{(i)} On $\R^+\times(\R^d\ba\L(r))\times \R^d$, $\LL_\e^{(r)}=-\infty$.
\item{(ii)} For all $(s,v)\in \R^+\times\R^d$, and all $e>0,r>0$ the function
$u\rightarrow\LL_\e^{(r)}(s,u,v)$ is upper semi-continuous (u.s.c.) at each
$u\in\L(r)$. 
\item{(iii)}   For all $(s,u)\in \R^+\times\L(r)$, the function
$\Phi_{\e,s,u}^{(r)}$ is convex and $\dom \Phi_{\e,s,u}^{(r)}=\R^d$.
}

\proof The proof is  trivial and is left to the reader.\endproof

The next Lemma relates the function $\LL^{(r)}_\e$ to the function 
$\LL_\e^{(r)}$ defined in \eqv(2.18). 

\lemma{\ver.11}{\it For any $(s,u,v^*)\in\R^+\times\R^d\times \R^d$,
$$
\left(\LL_\e^{(r)}\right)^*(s,u,v^*)=\LL_\e^{(r)*}(s,u,v^*)
\Eq(43.2)
$$
}

\proof 
We first prove that $\left(\LL_\e^{(r)}\right)^*(s,u,v^*)\geq
\LL_\e^{(r)*}(s,u,v^*)$.
For any $\wt v\in\R^d$,
$$
\eqalign{
\left(\LL_\e^{(r)}\right)^*(s,u,v^*)&
\geq (\wt v,v^*)-\LL_\e^{(r)}(s,u,\wt v)\cr
&= \inf_{s':|s-s'|\leq r} \inf_{u':|u-u'|\leq r}\left\{
 (\wt v,v^*)-L_\e(s',u',\wt v)\right\}\cr
}
\Eq(43.3)
$$
Now we choose for $\wt v$ the value s.t.
$$
\sup_{v\in\R^d}\left\{(v,v^*)-L_\e(s',u',v)\right\}
=(\wt v,v^*)-L_\e(s',u',\wt v)
\Eq(43.4)
$$
With this choice \eqv(43.3) becomes indeed
$$
\left(\LL_\e^{(r)}\right)^*(s,u,v^*)\geq  
\inf_{s':|s-s'|\leq r} \inf_{u':|u-u'|\leq r}L^*_\e(s',u',v^*)
=\LL_\e^{(r)*}(s,u,v^*)
\Eq(43.5)
$$
Next we show the converse inequality. Note that for any 
$\wt s, \wt u$ s.t. $|s-\wt s|\leq r,|u-\wt u|\leq r$, and any $v\in\R^d$,
$$
\sup_{s':|s-s'|\leq r} \sup_{u':|u-u'|\leq r}L_\e(s',u',v)
\geq L_\e(\wt s,\wt u,v)
\Eq(43.6)
$$
Hence
$$
\eqalign{&
\left(\LL_\e^{(r)}\right)^*(s,u,v^*)=
\sup_{v\in\R^d}\left\{
(v,v^*)-\sup_{s':|s-s'|\leq r} \sup_{u':|u-u'|\leq r}L_\e(s',u',v)\right\}
\cr
&\leq\sup_{v\in\R^d}\left\{
(v,v^*)- L_\e(\wt s,\wt u,v)\right\}=L_\e^*(\wt s,\wt u,v^*)
}
\Eq(43.7)
$$
Since \eqv(43.7) holds for all $\wt s,\wt u$ in the given sets,
it follows that 
$$
\left(\LL_\e^{(r)}\right)^*(s,u,v^*)\leq  
\inf_{\wt s:|\wt s-s|\leq r}
\inf_{\wt u:|\wt u-u|\leq r}L_\e^*(\wt s,\wt u,v^*)
=\LL_\e^{(r)*}(s,u,v^*)
\Eq(43.8)
$$
we obtain the desired inequality. The two inequalities imply \eqv(43.2).
\endproof

The previous Lemma allows to deduce the following analog of Lemma \ver.10:

\lemma{\ver.12} {\it 
\item{(i)} On $\R^+\times(\R^d\ba\L(r))\times \R^d$, $\LL_\e^{(r)*}=+\infty$.
\item{(ii)} For all $(s,v^*)\in \R^+\times\R^d$, and all $e>0,r>0$ the function
$u\rightarrow\LL_\e^{(r)*}(s,u,v^*)$ is lower semi-continuous (l.s.c.) at each
$u\in\L(r)$. 
\item{(iii)}   For all $(s,u)\in \R^+\times\L(r)$, the function
$\Phi_{\e,s,u}^{(r)*}$ is convex and for $(s,u)\in\R^+\times
\inte_\e\L(r)$, $\dom \Phi_{\e,s,u}^{(r)*}=\conv\D$.
}

Finally we come to the central result of this sub-section.

\lemma {\ver.13}{\it For any $r>0$ and for any 
closed bounded $\SS\subset\inte\L(r)$ the following holds:
\item{(i)} $\LL_\e^{(r)}$ converges uniformly to $\LL^{(r)}$ on 
$\R^+\times\SS\times\R^d$.
\item{(ii)} $\LL_\e^{(r)*}$ converges uniformly to $\LL^{(r)*}$ on 
$\R^+\times\SS\times\conv\D$.
}

\proof    Since (ii) follows from (i) in the same way as Lemma \ver.9 follows
from Lemma \ver.8, we concentrate on the proof of (i). Fix $r>0$. 
Define the sets
$$
\eqalign{
A_\e\equiv\bigl\{(s^*,u^*,v)&\in\R^+\times\L(r)\times\R^d\,\big|\,
\exists (s,u): |s-s^*|\leq r, |u-u^*|\leq r:
\cr&
\LL_\e(s^*,u^*,v)=\sup_{s':|s-s'|\leq r} \sup_{u':|u-u'|\leq r}L_\e(s',u',v)
\bigr\}
}
\Eq(43.9)
$$
and put
$$
A^{\e}\equiv\cup_{0\leq\e'\leq\e}A_{\e'}
\Eq(43.10)
$$

Define 
$$
\LL_{\e,\e_0}^{(r)}(s,u,v)\equiv
\lim_{\e\downarrow 0}\sup_{{s':|s-s'|\leq r, u':|u-u'|\leq r}\atop
{(s',u',v)\in A}} L_\e(s',u',v)
\Eq(43.11)
$$
Write
$$
\eqalign{
&\left| \LL_{\e}^{(r)}(s,u,v)-\LL^{(r)}(s,u,v)\right|\cr
&\leq\left|\LL_{\e}^{(r)}(s,u,v)-\LL_{\e,\e_0}^{(r)}(s,u,v)\right|+
\left|
\LL_{\e,\e_0}^{(r)}(s,u,v)\LL^{(r)}(s,u,v)\right|
}
\Eq(43.12)
$$
By definition of the set $\A^{\e}$, for $\e_0\geq \e$,
$$
\left|\LL_{\e}^{(r)}(s,u,v)-\LL_{\e,\e_0}^{(r)}(s,u,v)\right|=0
\Eq(43.13)
$$
On the other hand, for $(s^*,u^*,v)\in A^{\e_0}$, $\exists \e'\leq \e_0$       
and $(s,u)$ with $|s-s^*|\leq r, |u-u^*|\leq r$, such that for all
$(s',u')$ with   $|s-s'|\leq r, |u-u'|\leq r$,
$$
L_{\e'}(s^*,u^*,v)\geq L_{\e'}(s',u',v)
\Eq(43.14)
$$
Recalling the definition of $L_{\e'}$,
\eqv(43.14) implies that for any $\g>0$,
$$
\eqalign{
&\sum_{{\d\in\D}\atop{g_{\e'}(s^*,u^*,\d)\geq \g}}
e^{(\d,v)}g_{\e'}(s^*,u^*,\d)+
\sum_{{\d\in\D}\atop{g_{\e'}(s^*,u^*,\d)< \g}}
e^{(\d,v)}g_{\e'}(s^*,u^*,\d)
\cr\geq 
&\sum_{{\d\in\D}\atop{g_{\e'}(s^*,u^*,\d)\geq \g}}
e^{(\d,v)}g_{\e'}(s',u',\d)+
\sum_{{\d\in\D}\atop{g_{\e'}(s^*,u^*,\d)< \g}}
e^{(\d,v)}g_{\e'}(s',u',\d)
}
\Eq(43.15)
$$
The important point is now that since $\SS\subset \inte\L(r)$, no matter what
$u^*\in \SS$, there exists a $q=\dist(\SS,\L(r)^c)>0$, such that for some
$u'$ with $|u'-u|\leq r$. By Hypothesis 2.2, and the continuity assumptions 
of Hypothesis 2.3, one has that there exists a constant $c_q>0$ such that 
for all these points, and for all $\d\in\D$, $g_{\e'}(s',u',\d)>c_q$.
Choosing such $u'$ and $s'=s^*$, \eqv(43.15)
implies that
$$
(c_q-\g)\sum_{{\d\in\D}\atop{g_{\e'}(s^*,u^*,\d)< \g}}
e^{(\d,v)}
\leq
\sum_{{\d\in\D}\atop{g_{\e'}(s^*,u^*,\d)\geq \g}}
e^{(\d,v)}g_{\e'}(s^*,u^*,\d)
\Eq(43.16)
$$
By Hypothesis 2.4, $g_\e$ converges uniformly. Therefore, for any $\eta>0$,
there exists $\e_0>0$, such that for all $\e,\e'\leq \e_0$, 
and all $(s^*,u^*,\d)\in\R^+\times(\SS\cap\L)\times\R^d$,
$$
|g_{\e'}(s^*,u^*,\d)-g_{\e'}(s^*,u^*,\d)|\leq \eta
\Eq(43.17)
$$
Given $\eta$, let now $\e_0$ be such that \eqv(43.17) holds. Then \eqv(43.16)
implies that for all $\e\leq \e_0$,
$$
\eqalign{
(c_q-\g)\sum_{{\d\in\D}\atop{g_{\e}(s^*,u^*,\d)< \g-\eta}}
e^{(\d,v)}
\leq
&\sum_{{\d\in\D}\atop{g_{\e}(s^*,u^*,\d)\geq \g-\eta}}
e^{(\d,v)}(g_{\e}(s^*,u^*,\d)+\eta)
\cr&\leq(1+\frac \eta\g)
\sum_{{\d\in\D}\atop{g_{\e}(s^*,u^*,\d)\geq \g+\eta}}
e^{(\d,v)}g_{\e}(s^*,u^*,\d)
 }
\Eq(43.18)
$$
Therefore, for all $\e\leq \e_0$, and $(s^*,u^*,v)\in\A^{\e_0}$,
$$
\eqalign{
&\LL_{\e}(s^*,u^*,v)
\cr&=
\ln\left(\sum_{{\d\in\D}\atop{g_{\e}(s^*,u^*,\d)\geq \g}}
e^{(\d,v)}g_{\e}(s^*,u^*,\d)
\left(1+\frac{\sum_{{\d\in\D}\atop{g_{\e}(s^*,u^*,\d)< \g}}
e^{(\d,v)}g_{\e}(s^*,u^*,\d)}{\sum_{{\d\in\D}\atop{g_{\e}(s^*,u^*,\d)\geq \g}}
e^{(\d,v)}g_{\e}(s^*,u^*,\d)}\right)\right)\cr
&=\ln\sum_{{\d\in\D}\atop{g_{\e}(s^*,u^*,\d)\geq \g}}
e^{(\d,v)}g_{\e}(s^*,u^*,\d)\cr
&+\ln\left(1+\frac{\sum_{{\d\in\D}\atop{g_{\e}(s^*,u^*,\d)< \g}}
e^{(\d,v)}g_{\e}(s^*,u^*,\d)}{\sum_{{\d\in\D}\atop{g_{\e}(s^*,u^*,\d)\geq \g}}
e^{(\d,v)}g_{\e}(s^*,u^*,\d)}\right)
}
\Eq(43.19)
$$
The last term in \eqv(43.19) is bounded by 
$$
\ln\left(1+\frac {\g}{c_q-\g-\eta}\right)\leq \frac {\g}{c_q-\g-\eta}
\Eq(43.20)
$$
which will be made small by choosing $\g$ small enough.
On the other hand,
$$
\eqalign{
&\left|\ln\sum_{{\d\in\D}\atop{g_{\e}(s^*,u^*,\d)\geq \g}}
e^{(\d,v)}g_{\e}(s^*,u^*,\d)-
\ln\sum_{{\d\in\D}\atop{g_(s^*,u^*,\d)\geq \g}}
e^{(\d,v)}g(s^*,u^*,\d)\right|\leq \frac \eta\g
}
\Eq(43.21)
$$
Therefore, choosing $\g=\sqrt {c_q\eta}$, we see that for all
$\e\leq\e_0$,
$$
\left|
\LL_{\e,\e_0}^{(r)}(s,u,v)-\LL^{(r)}(s,u,v)\right|\leq 3\sqrt{ \eta/c_q}
\Eq(43.23)
$$
Combining both observations, we see that with $\e=\e_0$, we get in fact that
$$
\left| \LL_{\e_0}^{(r)}(s,u,v)-\LL^{(r)}(s,u,v)\right|\leq 3\sqrt{\eta c_q}
\Eq(43.24)
$$
which implies the desired uniform convergence and proves (i).
(ii) follows easily in the same way as the convergence result in Lemma
\ver.9 (v) follows from Lemma \ver.8 (v).\endproof

\proofof{Lemma 3.3} By definition, for any $(s,u,v^*)\in \R^+\times\L\times
\conv\D$,
$$
\ov{\LL}^*(s,u,v^*)=\liminf_{{s'\rightarrow s}\atop{u'\rightarrow u}}
{\LL}^{*}(s',u',v^*)
\Eq(43.25)
$$
But by Lemma \ver.11, the function ${\LL}^{*}(s,u,v^*)$ is jointly continuous 
in the variables $s,u$ at any   $(s,u,v^*)\in \R^+\times\inte\L\times
\conv\D$ so that on this set the right hand side of \eqv(43.25) 
coincides with ${\LL}^{*}(s,u,v^*)$. This proves Lemma 3.3.\endproof

\vskip.3truecm
\noindent{\bf \ver.4. A continuity derived result.}
\vskip.3truecm

We shall here be interested in the case $u\in\inte\L$ only. As seen in 
Lemma 4.7 the
conjugacy correspondence between   
$\Phi$ and $\Phi^*$
is closely connected to their differentiability properties. To this we may add:

\lemma{\ver.14}{\it Let $(s,u)\in\R^+\times\inte\L$. Then 
$\nabla\Phi^*(v^*)$ is bounded if and only if $v^*\in\ri(\dom\D)$.
}

\proof We know from Lemma \ver.5 and Lemma \ver.6
that for each $(s,u)\in\R^+\times\inte_{}\L$,  
$\Phi^*$ is a proper, closed, and strictly convex function 
having effective domain $\conv\D$. Moreover, we saw in the proof of 
lemma \ver.5 that the subgradient of $\Phi^*$ reduces to the
gradient mapping. Finally, invoking Theorem 23.4 of [Ro], the subgradient of 
$\Phi^*$ at $v^*$ is a non empty and 
bounded set if and only if $v^*\in\ri(\dom\D)$. The lemma is proven. \endproof

Now boundedness of $\nabla\Phi^*$ turns out to be an essential ingredient
of the proof of the large deviations estimates of Chapter 5. The 
particular place where it is needed appears in the context of
the minimisation problem of 
Lemma \ver.15 below. There, we shall see that the continuity property of 
$\Phi^*$, which in 
contrast with it's differentiability properties hold up to $\rbd(\dom \D)$, 
enables us to restrict ourselves to situations where $\nabla\Phi^*$ is
bounded. 

\lemma{\ver.15}{\it Let $\FF\subset\DD^{\circ}([0,T])$ be a convex subset of 
$\DD^{\circ}([0,T])$ and set 
$$
\GG\equiv
\left\{\psi\in{\cal F}\,\Big|\, \dot\psi(t)\in\ri(\conv\D),\,\,0\leq t\leq T
\right\}
\Eq(4.33)
$$
Then,
$$
\inf_{\psi\in{\cal F}}
\int_0^T dt{\cal L}^*(t,\psi(t),\dot\psi(t))
=
\inf_{\psi\in{\cal G}}
\int_0^T dt{\cal L}^*(t,\psi(t),\dot\psi(t))
\Eq(4.34)
$$
}

\proof With $\DD^{\circ}([0,T])$ defined in \eqv(2.34) recall that
$\Phi^*{t,\psi(t)}(\cdot)= {\cal L}^*(t,\psi(t),\cdot)$.
As seen in the proof of Lemma \ver.14, for $\psi\in\DD^{\circ}([0,T])$, 
$\Phi^*_{t,\psi(t)}$
is a proper, closed, strictly convex, and positive function 
having effective domain $\conv\D$. This in particular ensures that both sides 
of \eqv(4.34) are finite. 
 Since ${\cal F}\supseteq{\cal G}$,
$$
\inf_{\psi\in{\cal F}}
\int_0^T dt{\cal L}^*(t,\psi(t),\dot\psi(t))
\leq
\inf_{\psi\in{\cal G}}
\int_0^T dt{\cal L}^*(t,\psi(t),\dot\psi(t))
\Eq(4.38)
$$
and we only have to prove the reverse inequality. To do so we will use that 
for any $\psi_1\in\GG$ and any $\psi_2\in\FF$ the path $\a\psi_1+(1-\a)\psi_2$
belongs to $\GG$ for each $0<\a\leq 1$: obviously, by the convexity assumption 
on $\FF$,  $\a\psi_1+(1-\a)\psi_2\in\FF$;
but since for each $t\leq 0\leq 1$ $\dot\psi_1(t)$ is a point in 
$\ri(\conv\D)$ and
$\dot\psi_2(t)$ a point in $\conv\D$, the point 
$\a\dot\psi_1(t) +(1-\a)\dot\psi_2(t)$ lies in  
$\ri(\conv\D)$ for each $0<\a\leq 1$ (see [Ro], Theorem 6.1) so that
$\a\psi_1+(1-\a)\psi_2$ lies in $\GG$.
Thus, given $\psi_1\in\GG$ and $\psi_2\in\FF$ we have, for each $0<\a\leq 1$,
$$
\eqalign{
&\inf_{\psi\in{\cal G}}\int_0^T dt{\cal L}^*(t,\psi(t),\dot\psi(t))
\cr
&\leq
\int_0^T dt{\cal L}^*
(t,\psi_2(t)+\a[\psi_1(t)-\psi_2(t)],\dot\psi_2(t)+
\a[\dot\psi_1(t)-\dot\psi_2(t)])
}
\Eq(4.39)
$$
where the integrand in the last 
line is positive and bounded for each $0<\a\leq 1$.
Thus, taking the limit $\a\downarrow 0$, we may write, 
using Lebesgue's dominated 
convergence Theorem,
$$
\eqalign{
&\inf_{\psi\in{\cal G}}\int_0^T dt{\cal L}^*(t,\psi(t),\dot\psi(t))
\cr
&\leq
\lim_{\a\downarrow 0}\int_0^T dt{\cal L}^*
(t,\psi_2(t)+\a[\psi_1(t)-\psi_2(t)],\dot\psi_2(t)+
\a[\dot\psi_1(t)-\dot\psi_2(t)])
\cr
&=\int_0^T dt\lim_{\a\downarrow 0}{\cal L}^*
(t,\psi_2(t)+\a[\psi_1(t)-\psi_2(t)],\dot\psi_2(t)+
\a[\dot\psi_1(t)-\dot\psi_2(t)])
\cr
&=\int_0^T dt{\cal L}^*
(t,\psi_2(t),\dot\psi_2(t))
}
\Eq(4.40)
$$
where in the last line we used that  ${\cal L}^*(s,u,v^*)$ is jointly 
continuous in the variables $s,u$, and $v^*$ relative to $\DD^{\circ}([0,T])$
(see Lemma \ver.9, last line and  assertion (iii)). 
Finally,
since \eqv(4.40) is true for any $\psi_2\in\FF$, 
$$
\inf_{\psi\in\GG}\int_0^T dt{\cal L}^*(t,\psi(t),\dot\psi(t))
\leq 
\inf_{\psi\in\FF}\int_0^T dt{\cal L}^*(t,\psi(t),\dot\psi(t))
\Eq(4.41)
$$
which concludes the proof of the lemma.\endproof

\vskip.3truecm
\noindent{\bf \ver.5. Proof of Proposition 2.7.}
\vskip.3truecm

The proof of Proposition 2.7 goes along the same lines as that of Lemma
\ver.14. 

Let $\psi_1$ be any path in $\AA\cap\DD^{\circ}([0,T])$ and let 
$\psi_2$ be any path in $\AA\cap\ov\DD([0,T])$.
It follows from the convexity of $\AA$ together with the definitions of 
$\DD^{\circ}([0,T])$ and $\ov\DD([0,T])$ that the path 
$\a\dot\psi_1(t)+(1-\a)\dot\psi_2(t)$
lies in $\AA\cap\DD^{\circ}([0,T])$ for each $0<\a\leq 1$.
Hence, for each such $\a$,
$$
\eqalign{
&\inf_{\psi\subset\AA\cap \DD^\circ([0,T])}\displaystyle{\int_0^T 
\overline{\cal L}^*(t,\psi(t),\dot\psi(t))dt}
\cr
&\leq \int_0^T \overline{\cal L}^*(t,
\a\psi_1(t)+(1-\a)\psi_2(t),\a\dot\psi_1(t)+(1-\a)\dot\psi_2(t))dt
\cr
&\leq
\int_0^T \overline{\cal L}^*(t,\a\psi_1(t)+(1-\a)\psi_2(t),\dot\psi_2(t))dt\cr
&+\a\Bigl\{
\int_0^T \overline{\cal L}^*(t,
\a\psi_1(t)+(1-\a)\psi_2(t),\dot\psi_1(t))dt
\cr&\quad-
\int_0^T \overline{\cal L}^*(t,
\a\psi_1(t)+(1-\a)\psi_2(t),\dot\psi_2(t))dt
\Bigr\}
\cr}
\Eq(4.42)
$$
Now condition (i) implies that
$$
\lim_{\a\downarrow 0}\int_0^T 
\overline{\cal L}^*(t,\a\psi_1(t)+(1-\a)\psi_2(t),\dot\psi_2(t))dt\leq \int_0^T \overline{\cal L}^*(t,\psi_2(t),\dot\psi_2(t))dt
\Eq(4.43)
$$
while condition (ii) guarantees that 
$$
\eqalign{
\lim_{\a\downarrow 0}
&\Bigl\{
\int_0^T \overline{\cal L}^*(t,
\a\psi_1(t)+(1-\a)\psi_2(t),\dot\psi_1(t))dt
\cr&\quad-
\int_0^T \overline{\cal L}^*(t,
\a\psi_1(t)+(1-\a)\psi_2(t),\dot\psi_2(t))dt
\Bigr\}=0
}
\Eq(4.44)
$$
Since this is true for all $\psi_2\in\AA\cap\ov\DD([0,T])$, we have
$$
\inf_{\psi\subset\AA\cap \DD^\circ([0,T])}\displaystyle{\int_0^T 
\overline{\cal L}^*(t,\psi(t),\dot\psi(t))dt}
\leq\inf_{\psi\subset\AA\cap \ov\DD([0,T])}\displaystyle{\int_0^T 
\overline{\cal L}^*(t,\psi(t),\dot\psi(t))dt}
\Eq(4.45)
$$
As the reverse inequality trivially holds, the proposition is proven.
\endproof

\newpage


\chap{5. Proof of Proposition 3.2}5

We are now ready to prove the main estimates of the paper. Basically, the 
idea of the proof is simple and consist of exploiting the ``almost-independence''
of consecutive jumps over length scales large compared to $1$ but small 
compared to $1/\e$, as in Wentzell's work. 
The source of this almost independence are of course the 
regularity properties of the transition probabilities. On the basis of this
independence, we bring to bear classical Cram\'er type-techniques. The 
main difficulties arise from the non-uniformity of our regularity assumptions 
near the  boundaries.

The chapter is divided in three subchapters. We will first get
equipped with some preparatory tools. Armed with these, the basic
upper and lower bounds are next derived. Lastly, using results from
Chapter 4, the proof is brought to a close.
 From now on the letter $t$ will be used exclusively for time parameters taking
value in $[0, T]$ (that is, on `macroscopic scale' 1) while $k$ will be 
reserved 
for discrete time parameters (on `microscopic scale' $\e^{-1}$).

\vskip.3truecm
\noindent{\bf 5.1: Preparatory steps.}
\vskip.3truecm

Lemma \ver.1 below provides a covering of the ball $\BB_{\rho}(\phi)$ into 
basic `tubes'. 

$\L^c$ denotes the complement of $\L$ in $\R^d$. 
For $x\in\R^d$ and $A\subset\R^d$, $\dist(x,A)\equiv\inf_{y\in A}|x-y|$.
Recall that given $\rho> 0$ and $\phi\in{\cal E}([0,T])$,
$
\BB_{\rho}(\phi)=\left\{\psi\in{\cal E}([0,T])\,\Big|\, 
\max_{0\leq t\leq T}|\psi(t)-\phi(t)|<\rho\right\}
$.


\lemma{\ver.1}{\it Let $0=t_0<t_1<\dots<t_n=T$ be any 
partition of $[0,T]$ into $n$ intervals and set
$$
\tau\equiv\max_{0\leq i\leq n}\e^{-1}|t_{i+1}-t_i|
\Eq(5.1)
$$
For $\eta>0$ and 
$\psi\in{\cal E}([0,T])$ define,
$$
{\cal A}_{\eta}(\psi)=\left\{\psi'\in{\cal E}([0,T])\,\Big|\,
\max_{0\leq i\leq n}|\psi'(t_i)-\psi(t_i)|\leq 2\eta
\right\}
\Eq(5.3)
$$
and for $\g\geq 0$ define,
$$
\eqalign{
\BB_{\rho,\g}(\phi)&=\left\{\psi'\in\BB_{\rho}(\phi)\,\Big|\,
\inf_{0\leq t\leq T}\dist(\psi'(t),\L^c)\geq\g\right\}\cr
\bar\BB_{\rho,\g}(\phi)&=\left\{\psi'\in\bar\BB_{\rho}(\phi)\,\Big|\,
\inf_{0\leq t\leq T}\dist(\psi'(t),\L^c)\geq\g\right\}\cr
}
\Eq(5.2)
$$
the restrictions of $B_\rho(x)$ and its closure to the $\g$-interior of $\L$. 
\item{(i)}
For any $\g\geq 0$ and $\eta>0$ such that  $\rho>2\eta$, 
there exists a subset ${\cal R}_{\rho,\eta,\g}(\phi)$
of ${\cal E}([0,T])$ such that:
$$
{\cal R}_{\rho,\eta,\g}(\phi)\subset\bar\BB_{\rho,\g}(\phi)
\subset\bigcup_{\psi\in{\cal R}_{\rho,\eta,\g}(\phi)}{\cal A}_{\eta}(\psi)
\Eq(5.4)
$$
$$
\left|{\cal R}_{\rho,\eta,\g}(\phi)\right|\leq
e^{dn\left(\log\left(\frac{\rho}{\eta}\right)+2\right)},
\,\,\,\,\,\forall \g\geq 0
\Eq(5.6)
$$
\item{(ii)}
For any $\g\geq 0$ and $\eta>0$ such that $\rho>2(\eta+\e\tau\diam\D)$,
$$
\bigcup_{\psi\in
\BB_{\rho-2(\eta+\e\tau\diam\D),\g}(\phi)}{\cal A}_{\eta}(\psi)
\subset\BB_{\rho}(\phi)
\Eq(5.5)
$$
}

\proof 
The proof of \eqv(5.4) relies on the following construction.
Given $\eta>0$ let ${\cal W}_{\eta}$ be the Cartesian lattice in $\R^d$
with spacing $\frac{\eta}{\sqrt{d}}$. For $y\in\R^d$ set
${\cal W}_{\rho,\eta}(y)={\cal W}_{\eta}\cap
\{y'\in\R^d\mid |y'-y|\leq\rho\}$
and for $\phi\in{\cal E}([0,T])$, 
${\cal V}_{\rho,\eta}(\phi)=\times_{i=0}^n{\cal W}_{\rho,\eta}(\phi(t_i))$.
Next, for
$x=(x_0,\dots,x_n)\in {\cal V}_{\rho,\eta}(\phi)$, define
$$
A_{\rho,\eta,\g}(x)=
\left\{\psi'\in\bar\BB_{\rho,\g}(\phi)\,\Big|\,
\max_{0\leq i\leq n}|\psi'(t_i)-x_i|\leq \eta,
\right\}
\Eq(5.8)
$$
Thus $A_{\rho,\eta,\g}(x)$ is the set of paths in $\bar\BB_{\rho,\g}(\phi)$ 
which at time
$t_i$ are within a distance $\eta$ of the lattice point $x_i$. Obviously,
the collection of all (not necessarily disjoint and possibly empty sets)
$A_{\rho,\eta,\g}(x)$ form a covering of $\bar\BB_{\rho,\g}(\phi)$: 
$$
\bar\BB_{\rho,\g}(\phi)=
\bigcup_{x\in{\cal V}_{\rho,\eta}(\phi)}
A_{\rho,\eta,\g}(x)
\Eq(5.9)
$$
In each of those sets $A_{\rho,\eta,\g}(x)$ that are non empty
pick one element arbitrarily and label it $\psi_{x}$. Clearly
$\psi_{x}\in\bar\BB_{\rho,\g}(\phi)$. Moreover for all 
$\psi'\in A_{\rho,\eta,\g}(x)$,
$|\psi'(t_i)-\psi_{x}(t_i)|\leq 2\eta$ for all $i=0,\dots,n$, and hence
$A_{\rho,\eta,\g}(x)\subset {\cal A}_{\eta}(\psi_{x})$.
Putting these information together with \eqv(5.9) and 
taking
${\cal R}_{\rho,\eta,\g}(\phi)=\left\{\psi_{x}\mid x\in{\cal V}_{\rho,\eta}(\phi)
\right\}$ yields \eqv(5.4). 
Finally \eqv(5.6) follows from the bound
$
\left|{\cal R}_{\rho,\eta,\g}(\phi)\right|
\leq \left|{\cal V}_{\rho,\eta}(\phi)\right|
\leq\left(\max_i \left|{\cal W}_{\rho,\eta}(\phi(t_i))\right|\right)^n
$
together with
the estimate $\left|{\cal W}_{\rho,\eta}(y)\right|\leq
\exp\left\{d\left(\log\left(\frac{\rho}{\eta}\right)+2\right)\right\}$, 
$y\in\R^d$, whose (simple) proof can be found e.g. in [BG5].

We now prove \eqv(5.5). 
Set $\bar\rho\equiv\rho-2(\eta+\e\tau\diam\D)$.
Let 
$\psi'\in\bigcup_{\psi\in\BB_{\bar\rho,\g}(\phi)}{\cal A}_{\eta}(\psi)$.
Then
$\psi'\in{\cal A}_{\eta}(\psi)$ for some 
$\psi\in\BB_{\bar\rho,\g}(\phi)$. Hence,
$$
\eqalign{
\max_{0\leq t\leq T}|\psi(t)-\phi(t)|
&\leq \max_{0\leq t\leq T}(|\psi'(t)-\psi(t)|+|\psi(t)-\phi(t)|)\cr
&< \max_{0\leq t\leq T}|\psi'(t)-\psi(t)|+\bar\rho\cr
&< \max_{0\leq i\leq n}\max_{t_i\leq t\leq t_{i+1}}
(|\psi'(t)-\psi'(t_i)|+|\psi'(t_i)-\psi(t_i)|+|\psi(t)-\psi(t_i)|)+
\bar\rho\cr
&<  \max_{0\leq i\leq n}\max_{t_i\leq t\leq t_{i+1}}
(|\psi'(t)-\psi'(t_i)|+|\psi(t)-\psi(t_i)|)
+2\eta+\bar\rho\cr
}\Eq(5.14)
$$
Thus, using that for $\psi''\in{\cal E}([0,T])$,
$$
\eqalign{
\max_{0\leq i\leq n}\max_{t_i\leq t\leq t_{i+1}}|\psi''(t)-\psi''(t_i)|
\leq\max_{0\leq i\leq n}|t_{i+1}-t_i|\,\diam\D\leq\e\tau\,\diam\D
}
\Eq(5.15)
$$
\eqv(5.14) entails
$\psi'\in\BB_{\bar\rho+2(\eta+\e\tau\diam\D)}(\phi)$, proving \eqv(5.5).
Lemma \ver.1 is proven.\endproof

\remark  Note that in general $B_{\rho,0}(x)\neq B_\rho(x)$. However,
due to Lemma 3.1, it is true that 
$$
\wt\PP_{\e,\phi_0}\left(B_\rho(\phi)\right)=
\wt\PP_{\e,\phi_0}\left(B_{\rho,\g}(\phi)\right)
\Eq(5.-1)
$$
and the same holds true for the closed balls. Thus it will suffice to 
get upper and lower bounds for the set $B_{\rho,\g}$, for all $\g\geq 0$.
Therefore the following Lemma will be a sufficient starting point.

Lemma \ver.1 allows us to control the probabilities in path space 
by the probabilities of some discrete time observations of the chain. 
This is the content of the next lemma.

\lemma{\ver.2} {\it With the notation of
Lemma \ver.1, the following holds for any $0=t_0<t_1<\dots<t_n=
 T$, $t_i\in\R$, $n\in\N$. 
\item{(i)} For any $\g\geq 0$ and $\eta>0$ such that $\rho>2\eta$,
$$
\eqalign{
\log\widetilde{\cal P}_{\e,\phi_0}(\bar\BB_{\rho,\g}(\phi))
\leq &\sup_{\psi\in\bar\BB_{\rho,\g}(\phi)}
\log{\cal P}_{\e,\phi_0}
\left(\max_{0\leq i\leq n}\left|X(\left[\sfrac{t_i}{\e}\right])
-\psi(\e\left[\sfrac {t_i}\e\right])\right|\leq 2\eta+2\e\,\diam\D\right)\cr
&+dn\left(\log\left(\frac{\rho}{\eta}\right)+2\right)\cr
}
\Eq(5.16)
$$
\item{(ii)} For any $\g\geq 0$, any $\eta$ such that $\eta>\e\diam\D$ 
and $\rho>2(\eta+\e\tau\diam\D)$, and any 
$\psi\in\BB_{\rho-2(\eta+\e\tau\diam\D),\g}(\phi)$, 
$$
\log\widetilde{\cal P}_{\e,\phi_0}(\BB_{\rho}(\phi))\geq
\log{\cal P}_{\e,\phi_0}
\left(
\max_{0\leq i\leq n}\left|X(\left[\sfrac{t_i}{\e}\right])-\psi(\e\left[\sfrac {t_i}\e\right])\right|
<2\eta-2\e\,\diam\D
\right)
\Eq(5.17)
$$

}

\proof
We first prove assertion (i). Assume that $\eta,\rho$ and $\g$ satisfy
the conditions of Lemma \ver.1, (i). Then, by \eqv(5.4),
$$
\eqalign{
\widetilde{\cal P}_{\e,\phi_0}(\bar\BB_{\rho,\g}(\phi))
&\leq
\left|{\cal R}_{\rho,\eta,\g}(\phi)\right|
\exp\left\{\sup_{\psi\in{\cal R}_{\rho,\eta,\g}(\phi)}
\log\widetilde{\cal P}_{\e,\phi_0}({\cal A}_{\eta}(\psi))
\right\}
\cr
&\leq
\left|{\cal R}_{\rho,\eta,\g}(\phi)\right|
\exp\left\{\sup_{\psi\in \bar\BB_{\rho,\g}(\phi)}
\log\widetilde{\cal P}_{\e,\phi_0}({\cal A}_{\eta}(\psi))
\right\}
}
\Eq(5.18)
$$
Now
$$
\eqalign{
\widetilde{\cal P}_{\e,\phi_0}({\cal A}_{\eta}(\psi))
= &\widetilde{\cal P}_{\e,\phi_0}
\left(\max_{0\leq i\leq n}\left|Z(t_i)-\psi(t_i)\right|\leq2\eta\right)
\cr
\leq &{\cal P}_{\e,\phi_0}\left(
\max_{0\leq i\leq n}\left|X(\left[\sfrac{t_i}{\e}\right])
-\psi(\e\left[\sfrac {t_i}\e\right])\right|\leq2\eta+\e\,\diam\D\right)
\cr
}
\Eq(5.19)
$$
where we used that, (for $Z\in{\cal E}([0,T])$),
$$
\eqalign{
|Z(t)-\psi(t)|
&\geq\left|Z(\e\left[\sfrac{t}{\e}\right])-\psi(\e\left[\sfrac t\e\right])
\right|
-\left|Z(t)-Z(\e\left[\sfrac{t}{\e}\right])\right|-\left|\psi(t)-\psi(\e\left[
\sfrac t\e\right])\right|
\cr
&\geq
\left|X(\left[\sfrac{t}{\e}\right])-\psi(t)\right|
-2\left|t-\e\left[\sfrac{t}{\e}\right]\right|\diam\D
\cr
&\geq \left|X(\left[\sfrac{t}{\e}\right])-\psi(t)\right|
-2\e\,\diam\D\cr
}
\Eq(5.20)
$$
Inserting \eqv(5.6) and \eqv(5.19) in \eqv(5.18) gives \eqv(5.16).
Similarly we derive assertion (ii) of Lemma \ver.2 from 
assertion (ii) of Lemma \ver.1, writing first that by \eqv(5.5),
for any 
$\psi\in\BB_{\rho-2(\eta+\e\tau\diam\D),\g}(\phi)$, 
$$
\log\widetilde{\cal P}_{\e,\phi_0}(\BB_{\rho}(\phi))\geq
\log\widetilde{\cal P}_{\e,\phi_0}({\cal A}_{\eta}(\psi))
\Eq(5.21)
$$
and using next that, since $Z\in{\cal E}([0,T])$, analogous to \eqv(5.20),
$$
\eqalign{
|Z(t)-\psi(t)|
&\leq 2\e\,\diam\D+ \left|X(\left[\sfrac{t}{\e}\right])-
\psi(\e\left[\sfrac t\e\right])\right|\cr
}
\Eq(5.22)
$$
so that 
$$
\eqalign{
\widetilde{\cal P}_{\e,\phi_0}({\cal A}_{\eta}(\psi))
\geq
{\cal P}_{\e,\psi_0}\left(
\max_{0\leq i\leq n}\left|X(\left[\sfrac{t_i}{\e}\right])-
\psi(\e\left[\sfrac{t_i}\e\right])\right|
\leq 2\eta-\e\diam\D
\right)
\cr
}
\Eq(5.23)
$$
This concludes the proof of Lemma \ver.2. \endproof

\remark We could arrange to use Lemma \ver.2 with $t_i$ that are multiples
of $\e$ only, except that $t_n=T$ has to be allowed to be what it wants to 
be. Thus we prefer to write the more homogeneous form above.

In view of Lemma \ver.2 the problem is reduced to estimating the probability 
 that the chain $X(t)$ be pinned in a small neighbourhood of a 
prescribed point $\psi(t_i)$ at each time $t_i$. As explained earlier
we will do this by comparing the chain in each time
interval $[t_{i-1}, t_{i})$ with a random walk whose steps, on microscopic 
time scale, take value in $\D$ and are distributed according to 
$p_{\e}(\left[{t_{i-1}}/{\e}\right],\psi(t_{i-1}),\cdot)$.
Let $P_{\e,k}=(p_{\e}(k,x,y))_{y\in\G_{\e},x\in\G_{\e}}$ denote the transition
matrix of the chain at time $k$ and, for $\ell\geq 1$, let
$P_{\e,k}^{(k, k+\ell)}=
\left(P_{\e,k}^{(k, k+\ell)}(x,y)\right)_{y\in\G_{\e},x\in\G_{\e}}$
denote the matrix product
$$
P_{\e,k}^{(k, k+\ell)}=\prod_{l=1}^{k}P_{\e,k+l-1}
\Eq(5.01)
$$
Set $k_i\equiv \left[\sfrac{t_i}{\e}\right]$. By the Markov property, 
for $\zeta>0$,
$$
\eqalign{
&{\cal P}_{\e,\phi_0}
\left(\max_{0\leq i\leq n}\left|X(\left[\sfrac{t_i}{\e}\right])
-\psi(t_i)\right|\leq\zeta\right)\cr
= &\sum_{x(k_0)\in\G_{\e}}
\pi_{\e,\phi_0}(x(k_0))
\1_{\{|x(k_0)-\psi(t_0)|\leq\zeta\}}
\sum_{x(k_1)\in\G_{\e}}
\1_{\{|x(k_1)-\psi(t_1)|\leq\zeta\}}
P_{\e,k_0}^{(k_0, k_1)}(x(k_0),x(k_1))
\dots
\cr
&\,\,\,\,\,\,\,\,\,\,\,\,\,\,\,\,\,\,\,\,\,\,\,\,\,\,\,\,\,\,\,\,
\,\,\,\,\,\,\,\,\,\,\,\,\,\,\,\,\,\,\,\,\,
\dots
\sum_{x(k_{n})\in\G_{\e}}
\1_{\{|x(k_n)-\psi(t_n)|\leq\zeta\}}
P_{\e,k_{n-1}}^{(k_{n-1}, k_{n})}(x(k_{n-1}),x(k_{n}))
\cr
}
\Eq(5.27)
$$
The following lemma provides estimates for terms of the form
$P_{\e,k_{i-1}}^{(k_{i-1}, k_{i})}(x(k_{i-1}),x(k_{i}))$.

\lemma{\ver.3}{\it
Let ${\cal S}$ be any closed bounded subset of $\inte\L$. 
Let ${\cal S'}$ be
an open subset of ${\cal S}$ and, for $\ell$ an integer, assume that the 
following condition is satisfied: for each $\ell\geq 1$ and $\e>0$ small
enough,
$$
\inf_{x\in{\cal S'}}\dist(x, {\cal S}^c)>\e\ell\diam\D
\Eq(5.02)
$$
For $r\geq 0$ set
$$
q(\ell,r)
=\e\sfrac{\ell^2}{2}(\theta +{\vartheta({\cal S})}{\diam}\D)
+\ell(r+2\e K({\cal S}))
\Eq(5.09)
$$
with $\theta$, $\vartheta({\cal S})$ and $K({\cal S})$ as in Hypothesis 2.3. 
Then, for any $x\in{\cal S'}$ and any $z\in{\cal S'}$,
$$
P_{\e,k}^{(k, k+\ell)}(x,y)\lg e^{\pm q(\ell,|x-z|)}
\sum_{\d(1)\in\D}\dots\sum_{\d(\ell-1)\in\D}
\prod_{l=1}^{\ell}e^{f^{(0)}_{\e}\left(\e k,z,\d(l)\right)}
\1_{\{
\e^{-1}(y-x)-\sum_{m=1}^{\ell-1}\d(m)
\in\D\}}
\Eq(5.010)
$$
}

\proof
First note that if $y$ is such that $P_{\e,k}^{(k, k+\ell)}(x,y)=0$ then
$\1_{\{\e^{-1}(y-x)-\sum_{m=1}^{\ell-1}\d(m)\in\D\}}=0$ 
for all sequences 
$(\d(1),\dots,\d(\ell-1))\in\times_{l=1}^{\ell-1}\D$, and hence 
\eqv(5.010) holds true. Assume that $y$ is such that  
$P_{\e,k}^{(k, k+\ell)}(x,y)\neq 0$ and set 
$x(k)\equiv x$, $x(k+\ell)\equiv y$, and
$$
\eqalign{
\d(0)&\equiv 0 \cr
\d(\ell)&\equiv \e^{-1}(x(k+\ell)-x(k))-\sum_{m=1}^{\ell-1}\d(m)\cr
}
\Eq(5.010bis)
$$
(We slightly abuse the notation in that $\d(0)$ and $\d(\ell)$ do not 
necessarily belong to $\D$).
By \eqv(5.01), 
$$
\eqalign{
&P_{\e,k}^{(k, k+\ell)}(x(k),x(k+\ell))\cr
=&
\sum_{x(k+1)\in\G_{\e}}\dots\sum_{x(k+\ell-1)\in\G_{\e}}
\prod_{l=1}^{\ell}p_{\e}(k+l-1,x(k+l-1),x(k+l))\cr
=&
\sum_{\d(1)\in\D}\dots\sum_{\d(\ell-1)\in\D}
\prod_{l=1}^{\ell}
p_{\e}\left(k+l-1,x(k)+\e\sum_{m=0}^{l-1}\d(m),x(k)+\e\sum_{m=1}^{l}\d(m)
\right)
\1_{\{\d(\ell)\in\D\}}
\cr
}
\Eq(5.013)
$$
Note that since 
$$
\e\sup|\sum_{m=1}^{\ell}\d(m)|\leq \e \ell\diam\D
\Eq(5.012ter)
$$
it follows from \eqv(5.02) that 
$$
\inf_{x\in{\cal S'}}\dist(x, {\cal S}^c)>\e\sup|\sum_{m=1}^{\ell}\d(m)|
\Eq(5.012bis)
$$
so that the chain 
starting at $x(k)\in{\cal S}'$ at time $k$ cannot reach the boundary of 
${\cal S}$ by time $k+\ell$. This in particular implies that 
for each $x(k)\in{\cal S}'$,
each sequence
$(\d(1),\dots,\d(\ell-1))\in\times_{l=1}^{\ell-1}\D$, and each 
$l=1,\dots,\ell-1$,
$$
x(k)+\e\sum_{m=1}^{l}\d(m)\in\inte_{\e}{\cal S}\subset\inte_{\e}\L
\Eq(5.012)
$$
Thus by \eqv(2.1) and Hypothesis 2.2 (see e.g. \eqv(2.8)), each of the 
probabilities in the last line of \eqv(5.013) is strictly positive.
In addition, under our assumption on $z$, by $(H0)$ of Hypothesis 2.3,
$e^{f^{(0)}_{\e}\left(\e k,z,\d(l)\right)}>0$. We may thus write 
$$
P_{\e,k}^{(k, k+\ell)}(x(k),x(k+\ell))
=
\sum_{\d(1)\in\D}\dots\sum_{\d(\ell-1)\in\D}
\prod_{l'=1}^{\ell}R_{l'}
\prod_{l=1}^{\ell}e^{f^{(0)}_{\e}\left(\e k,z,\d(l)\right)}
\1_{\{\d(\ell)\in\D\}}
\Eq(5.015)
$$
where
$$
R_l\equiv
p_{\e}\left(k+l-1,x(k)+\e\sum_{m=0}^{l-1}\d(m),x(k)+\e\sum_{m=1}^{l}\d(m)
\right)
e^{-f^{(0)}_{\e}\left(\e k,z,\d(l)\right)},
\,\,\,\,\,\,\forall l=1,\dots,\ell
\Eq(5.016)
$$
Setting $k'=k+l-1$ and $x'=x(k)+\e\sum_{m=0}^{l-1}\d(m)$ and 
using \eqv(2.1) and \eqv(2.10), we have
$$
\eqalign{
\left|\log R_l\right|
&=\left|f_{\e}(\e k',x',\d(l))-f^{(0)}_{\e}(\e k,z,\d(l))\right|\cr
&\leq
\e\left|f^{(1)}_{\e}(\e k',x',\d)\right|
+\left|f^{(0)}_{\e}(\e k',x',\d)-f^{(1)}_{\e}(\e k,z,\d(l))\right|\cr
}
\Eq(5.ex46)
$$
where by $(H1)$ of Hypothesis 2.3, 
$\left|f^{(1)}_{\e}(\e k',x',\d(l))\right|\leq K({\cal S})$
and by $(H2)$ and $(H3)$ of Hypothesis 2.3,
$$
\eqalign{
&\left|f^{(0)}_{\e}(\e k',x',\d(l))-f^{(0)}_{\e}(\e k,z,\d(l))\right|\cr
\leq &
\left|f^{(0)}_{\e}(\e k',x',\d(l))-f^{(0)}_{\e}(\e k,x',\d(l))\right|
+\left|f^{(0)}_{\e}(\e k,x',\d(l))-f^{(0)}_{\e}(\e k,z,\d(l))\right|\cr
\leq &\e\theta|{k}-{k}'|+\vartheta({\cal S})|z-x'|
}
\Eq(5.ex47)
$$
Thus 
$$
\left|\log R_l\right|\leq
\e\theta l+\vartheta({\cal S})
\left|(x(k)-z)+\e{\textstyle\sum_{m=1}^{l-1}}\d(m)\right|+\e K({\cal S})
\Eq(5.ex48)
$$
and for $\d(\ell)\in\D$, 
we have
$$
\eqalign{
\left|\log\left(\prod_{l=1}^{\ell}R_{l}\right)\right|&\leq
\sum_{l=1}^{\ell}
\left(\e\theta l+\vartheta({\cal S})
\left|(x(k)-z)+\e{\textstyle\sum_{m=1}^{l-1}}\d(m)\right|+\e K({\cal S})\right)
\cr
&\leq
\e\theta \sfrac{\ell(\ell-1)}{2}
+\e{\vartheta({\cal S})}{\diam\D}{\sfrac{\ell(\ell-1)}{2}}
+\vartheta({\cal S})\ell|x(k)-z|+\e\ell K({\cal S})
\cr
}
\Eq(5.017)
$$
Inserting the bound \eqv(5.017) in \eqv(5.015) yields \eqv(5.010).
This concludes the proof of the lemma.
\endproof

\vskip.3truecm
\noindent{\bf 5.2: Basic upper and lower large deviation estimates.}
\vskip.3truecm
We define the following sets: 
$$
\bar\L_{\rho,\gamma}(\phi)=
\left\{
x\in\L\mid \exists\psi\in\bar\BB_{\rho,\g}(\phi), \,
\exists t\in[0,T]\,\, s.t.\,\, \psi(t)=x
\right\}
\Eq(5.a1)
$$
$$
\bar{\cal S}_{\rho,r}(\phi)=
\cl \left(\left\{
x\in\L\mid \dist\left(x,\bar\L_{\rho,\gamma}(\phi)\right)\leq r
\right\}\right),\,\,\,\,\,\,\,r\geq 0
\Eq(5.a2)
$$
Observe that for $r<\gamma$, $\bar{\cal S}_{\rho,r}(\phi)$ is a closed 
bounded subset of $\inte\L$.
$$
{\cal T}(\phi_0)=\phi_0+\left[-(T+\e\sqrt{d})\diam\D, 
(T+\e\sqrt{d})\diam\D\right]^d
\Eq(5.a3)
$$
(this definition has to do with the fact that the initial condition 
$\pi_{\e,\phi_0}$ of the chain has support in
$\{x\in\G_{\e}\mid |x-\phi_0|\leq \e\sqrt{d}\}$). Finally, 
$$
{\cal S}_{\gamma/2}(\phi_0)=\cl\left(\left\{
x\in\L\mid\dist\left(x, ({\cal T}(\phi_0)\cap\L)^c\right)\geq \gamma/2
\right\}\right)
\Eq(5.a4)
$$

The upper bound we will prove is analogous to that of [DEW].

\lemma{\ver.4} {\it Let 
$0=t_0<t_1<\dots<t_n=\e\left[\sfrac T\e\right]$
be such that for all
$0\leq i\leq n-1$,
$
t_i=\e\left[\sfrac{t_i}{\e}\right]\equiv\e k_i, \,\,k_i\in\N$.
Assume that the conditions of Lemma \ver.2, (i), are verified 
and set $\bar\zeta= 2\eta+2\e\,\diam\D$. For any fixed $r>0$ assume that 
$\eta$, $\e$ and $\tau$ are such that
$$
r>2\bar\zeta+\e\tau\diam\D
\Eq(5.24bis)
$$
Then the following conclusions hold
for any $\psi$ in $\bar\BB_{\rho,0}(\phi)$.

\item{(i)} If $|\psi(t_0)-\phi_0|\leq \bar\zeta+\e\sqrt{d}$ then,
$$
\eqalign{
 &\e\log{\cal P}_{\e,\phi_0}
\left(\max_{0\leq i\leq n}\left|X(\sfrac{t_i}{\e})
-\psi(t_i)\right|\leq\bar\zeta\right)
\cr&
\leq \sup_{\psi':\forall_{i=0}^n|\psi'(t_i)-\psi(t_i)|\leq \bar\z}
\left(-\sum_{i=1}^n 
(t_i-t_{i-1})
\LL^{(r)*}_\e
\left(t_{i-1}, \psi'(t_i),\sfrac{\psi'(t_i)-\psi'(t_{i-1})}{t_i-t_{i-1}}\right)
\right)
}
\Eq(5.25)
$$
\item{(ii)} If 
$|\psi(t_0)-\phi_0|>\bar\zeta+\e\sqrt{d}$ then,
$$
\e\log{\cal P}_{\e,\phi_0}
\left(\max_{0\leq i\leq n}\left|X(\sfrac{t_i}{\e})
-\psi(t_i)\right|\leq\bar\zeta\right)
=-\infty 
\Eq(5.25a)
$$
}

\proof 
 The proof starts from equation \eqv(5.27), replacing $\z$ by $\bar\z$. 
We follow the procedure 
used by Varadhan [Va] for the multidimensional Cram\'er theorem\note{This 
allows us to avoid Wentzell's
 assumptions of boundedness of the derivatives of the
Lagrangian function $\LL^*$ with respect to the velocities.} 
 and write
$$
\eqalign{
&\prod_{i=0}^n\1_{\{|x(k_i)-\psi(t_i)|\leq \bar\z\}}
\leq \inf_{\bar\l_1,\dots,\bar\l_n\in\R^d}\sup_{{\psi'(t_1),\dots,\psi'(t_n)}
\atop{\forall_{i} |\psi'(t_i)-\psi(t_i)|\leq\bar\z}}
e^{\sum_{i=1}^n\left(\bar\l_i, x(k_i)-\psi'(t_i)\right)}
\cr &\quad\times \prod_{i=0}^n\1_{\{|x(k_i)-\psi(t_i)|\leq \bar\z\}} \cr
&=\inf_{\bar\l_1,\dots,\bar\l_n\in\R^d}\sup_{{\psi'(t_1),\dots,\psi'(t_n)}
\atop{\forall_{i} |\psi'(t_i)-\psi(t_i)|\leq\bar\z}}
\prod_{i=0}^n\1_{\{|x(k_i)-\psi(t_i)|\leq \bar\z\}}\cr
&\times e^{\sum_{i=1}^n\left(\bigl(\sum_{j=i}^n\bar\l_j\bigr),
\bigl(x(k_i)-x(k_{i-1})- \psi'(t_i)+\psi'(t_{i-1})\bigr)\right)}
\cr
&\times  e^{\left(\bigl(\sum_{j=1}^n\bar\l_j\bigr),x(k_0)-\psi(t_0)\right)}
\cr&\leq 
\inf_{\l_1,\dots,\l_n\in\R^d}\sup_{{\psi'(t_1),\dots,\psi'(t_n)}
\atop{\forall_{i} |\psi'(t_i)-\psi(t_i)|\leq\bar\z}}
\prod_{i=0}^n\1_{\{|x(k_i)-\psi(t_i)|\leq \bar\z\}}\cr
&\times e^{\sum_{i=2}^n\left(\l_i,x(k_i)-x(k_{i-1})\right)-
\left(\l_i, \psi'(t_i)-\psi'(t_{i-1})\right)} 
\cr
&\times  e^{\left(\l_1,x(k_1)-x(k_0)\right)-
\left(\l_1, \psi'(t_1)-\psi(t_{0})+x_0-\psi(t_0)\right)}
}
\Eq(50.1)
$$
We now insert \eqv(50.1) into \eqv(5.27). Relaxing  
all constraints on the endpoints of summations 
(this is reasonable since we already  assume 
that $\psi(t) $ remains in $\L$) we obtain, using \eqv(5.013),
$$
\eqalign{
&{\cal P}_{\e,\phi_0}
\left(\max_{0\leq i\leq n}\left|X(\left[\sfrac{t_i}{\e}\right])
-\psi(t_i)\right|\leq\bar\zeta\right)
\cr
&\leq \sum_{x(k_0)\in\G_\e}\pi_{\e,\phi_0}(x(k_0)) 
\1_{\{|x(k_0)-\psi(t_0)|\leq \bar\z\}}
\cr&\times\inf_{\l_1,\dots,\l_n\in\R^d}\sup_{{\psi'(t_1),\dots,\psi'(t_n)}
\atop{\forall_{i} |\psi'(t_i)-\psi(t_i)|\leq\bar\z}}
\prod_{i=2}^n\Biggl( 
e^{-\left(\l_i,\psi'(t_i)-\psi'(t_{i-1})\right)}\cr
&\times\sup_{x(k_i)\in\G_\e:|x(k_{i-1})-\psi(t_{i-1})|\leq \bar \z}
\sum_{\d(1),\dots,\d(\ell_i)}
\prod_{l=1}^{\ell_i} e^{f_\e\left(t_{i-1}+l-1,x(k_{i-1})+\e
\sum_{k=1}^{l-1} \d(k),
\d(l)\right)}
e^{(\e\l_i,\d(l))} \Biggr)\cr
&\times
\left(e^{-\left(\l_1,\psi'(t_1)-x(k_0)\right)}
\sum_{\d(1),\dots,\d(\ell_1)}
\prod_{l=1}^{\ell_1} e^{f_\e\left(t_{0}+l-1,x(k_{0})+
\e\sum_{k=1}^{l-1} \d(k),
\d(l)\right)}
e^{(\e\l_1,\d(l))}\right) 
}
\Eq(50.13)
$$
where $\ell_i\equiv k_{i+1}-k_i$. Taking into account the constraints on
the suprema over
the $x(k_i)$ and the $\psi(t_i)$, we see that all terms
$x(k_i) +\e\sum_{k=1}^{l-1} \d(k)$ appearing satisfy
$|x(k_i) +\e\sum_{k=1}^{l-1} \d(k)-\psi(t_i)|\leq 2\zeta +\e\t\diam\D$.
Therefore, for $r> 2\bar\zeta+\e\t\diam\D$,
$$
\eqalign{\sup_{x(k_i)\in\G_\e:|x(k_{i-1})-\psi(t_{i-1})|\leq \bar \z}
\sum_{\d(l)}
&e^{f_\e(t_{i-1}+\e(l-1),x(k_{i-1})+
\e\sum_{k=1}^{l-1} \d(k),
\d(l))}
e^{(\e\l_i,\d(l))}\cr
&\leq  
\sup_{t':|t'-t_{i-1}|\leq r}\sup_{u:|u-\psi(t_{i-1})|\leq r}
\sum_{\d(l)}
e^{f_\e(t',u,
\d(l))}e^{(\e\l_i,\d(l))}
} 
\Eq(50.14)
$$
to bound all the summations over the $\d(l)$ successively. This leads with the
above notation to the bound
$$
\eqalign{
&{\cal P}_{\e,\phi_0}
\left(\max_{0\leq i\leq n}\left|X(\sfrac{t_i}{\e})
-\psi(t_i)\right|\leq\bar\zeta\right)
\cr
&\leq \sum_{x(k_0)\in\G_\e}\pi_{\e,\phi_0}(x(k_0))
\1_{\{|x(k_0)-\psi(t_0)|\leq \bar\z\}}
\cr&\times\inf_{\l_1,\dots,\l_n\in\R^d}\sup_{{\psi'(t_1),\dots,\psi'(t_n)}
\atop{\forall_{i} |\psi'(t_i)-\psi(t_i)|\leq\bar\z}}
\prod_{i=2}^n
e^{-\left(\l_i,\psi'(t_i)-\psi'(t_{i-1})\right)+\ell_i  
\LL_\e^{(r)}(t_{i-1},\psi(t_{i-1}),\e\l_i)}\cr
&\times
e^{-\left(\l_1,\psi'(t_1)-x(k_0)\right)+{\ell_1} 
\LL_\e^{(r)}(t_{0},\psi(t_{0}),\e\l_1)}
}
\Eq(50.15)
$$
Using that for $|\psi-\psi'|\leq \bar\z$,
$\sup_{u:|u-\psi|\leq r} L_\e(t,u,v)\leq \sup_{u:|u-\psi'|\leq r+\bar\z} 
L_\e(t,u,v)$, we can replace $\psi(t_{i-1})$ by $\psi'(t_{i-1})$
in the second argument of $\LL^{(r)}_\e$ at the expense of increasing
$r$ by $\bar\z$ (which will lead to the condition
$r>2\bar\zeta +\e\t\diam\D$). 
The argument in the $\inf\sup$ is convex in the variables $\l_i$ and
concave (since linear) in the $\psi'(t_i) $ and verifies the 
assumptions of the minimax theorem (see [Ro], Section 37 Corollary 37.3.1.)
so that we may interchange the order in which they are taken. Thus we obtain
$$
\eqalign{
&{\cal P}_{\e,\phi_0}
\left(\max_{0\leq i\leq n}\left|X(\sfrac{t_i}{\e})
-\psi(t_i)\right|\leq\bar\zeta\right)
\leq \sum_{x(k_0)\in\G_\e}\pi_{\e,\phi_0}(x(k_0))
\1_{\{|x(k_0)-\psi(t_0)|\leq \bar\z\}}
\cr
&\times\sup_{{\psi'(t_0),\dots,\psi'(t_n)}
\atop{\forall_{i} |\psi'(t_i)-\psi(t_i)|\leq\bar\z}}
\exp\left(-\e^{-1}\sum_{i=1}^n (t_i-t_{i-1})  
\LL_\e^{(r)}{}^*\left(t_{i-1},\psi'(t_{i-1}),\sfrac {\psi'(t_i)-\psi'(t_{i-1})}
{t_i-t_{i-1}}\right)\right)
}
\Eq(50.16)
$$
The first factor in the last line is always less than one which implies (i) and
is zero 
if $|\psi(t_0)-\phi(0)|> \bar\z+\e\sqrt d$.
This implies (ii).\endproof

We now turn to the lower bound. Recall from  \eqv(4.3) that
$\Phi_{\e,t_{i-1},\psi(t_{i-1})}(\,\cdot\,)=
{\cal L}_{\e}(t_{i-1},\psi(t_{i-1}),\,\cdot\,)$.

\lemma{\ver.5}{\it The notation is the same as in Lemma \ver.4.
Assume that the conditions of Lemma \ver.2, (ii), are verified
and set $\zeta\equiv 2\eta-2\e\diam\D$. 
Define the set
$$
{\cal E}^{\circ}([0,T])=
\left\{\psi\in{\cal E}([0,T])\,\Big|\,
\sfrac{\psi(t)-\psi(t')}{t-t'}\in\ri(\conv\D)\,\,\forall t\in[0,T],
\forall t'\in[0,T], t\neq t'
\right\}
\Eq(5.25riter)
$$
Then, for any $\psi$ in
$$
\BB_{\rho-2(\eta+\e\tau\diam\D),\g}(\phi)\cap
{\cal E}^{\circ}([0,T])
\Eq(5.25ribis)
$$
there exist positive constants $c_0\equiv c_0(\psi)<\infty$  
such that, if
$\eta$, $\e$, and $\tau$ are such that
$$
\sfrac{\g}{2}\geq\zeta+\e\tau\diam\D
\text{and}
\sqrt{2\e T}\diam \D+{\e\sqrt{d}}<\zeta,
\Eq(5.25riribis)
$$
the following holds:
$$
\eqalign{
&\e\log{\cal P}_{\e,\phi_0}
\left(\max_{0\leq i\leq n}\left|X(\sfrac{t_i}{\e})
-\psi(t_i)\right|\leq\zeta\right)
\cr
\geq &
\cases
-{\displaystyle{\sum_{i=1}^n}} (t_i-t_{i-1}){\cal L}_{\e}^*
\left(t_{i-1},\psi(t_{i-1}),\frac{\psi(t_i)-\psi(t_{i-1})}{t_i-t_{i-1}}\right)
-Q\left(
\bar{\cal S}_{\rho,\gamma/2}(\phi),
\zeta, c_0\right)
& \hbox{if}\,\, |\psi(t_0)-\phi_0|\leq \e\sqrt{d}\cr
&\cr
-\infty &\hbox{ otherwise}\cr
\endcases
\cr
}
\Eq(5.26bis)
$$
where
$$
Q({\cal S},\zeta, c_0)
\equiv
3n(\e\tau)^2(\theta+\vartheta({\cal S})\diam\D)
+3 T\left(\zeta+2\e K({\cal S})\right)
+4n\zeta c_0
+\e\log(8d^2+4)
\Eq(5.25qua)
$$

}

\proof Obviously, for any $\varrho\leq\zeta$,
$$
{\cal P}_{\e,\phi_0}
\left(\max_{0\leq i\leq n}\left|X(\sfrac{t_i}{\e})
-\psi(t_i)\right|\leq\zeta\right)
\geq {\cal P}_{\e,\phi_0}
\left(\max_{0\leq i\leq n}\left|X(\sfrac{t_i}{\e})
-\psi(t_i)\right|\leq\varrho\right)
\Eq(5.40)
$$
As will turn out, the generic term for which we shall want a lower bound
is  of the form:
$$
{\cal T'}_i\equiv
\1_{\{|x(k_{i-1})-\psi(t_{i-1})|\leq\varrho\}}
\sum_{x(k_{i})\in\G_{\e}}
\prod_{j=i}^n
\1_{\{|(x(k_i)-\psi(t_i))+a_{i,j}|\leq\varrho\}}
P_{\e,k_{i-1}}^{(k_{i-1}, k_{i})}(x(k_{i-1}),x(k_{i}))
\Eq(5.41)
$$
where, for each $j=i,\dots,n$, $a_{i,j}\in\R^d$ is independent of
$\{x(k_{j})\}_{i\leq j\leq n}$. We shall however only treat the term
$$
{\cal T}_i\equiv
\1_{\{|x(k_{i-1})-\psi(t_{i-1})|\leq\varrho\}}
\sum_{x(k_{i})\in\G_{\e}}
\1_{\{|(x(k_i)-\psi(t_i))+a|\leq\varrho\}}
P_{\e,k_{i-1}}^{(k_{i-1}, k_{i})}(x(k_{i-1}),x(k_{i}))
\Eq(5.41bis)
$$
for $a\in\R^d$ an arbitrary constant, the extension of the resulting bound
to ${\cal T'}_i$ being straightforward. Naturally our bound on ${\cal T}_i$ 
will  be derived by means of Lemma \ver.3.
Let ${\cal G}$ denote the set \eqv(5.25ribis).
Since $\psi$ belongs to ${\cal G}$ it belongs in particular
to $\BB_{\rho-2(\eta+\e\tau\diam\D),\g}$ 
and hence to $\bar\BB_{\rho,\g}$. Thus, 
under the assumptions \eqv(5.25riribis), we may apply
Lemma \ver.3 with 
$\ell\equiv\tau$,
${\cal S}\equiv \bar{\cal S}_{\rho,\gamma/2}(\phi)$,
${\cal S}'\equiv \bar{\cal S}_{\rho,\zeta}(\phi)$, and,
in each time interval $(k_{i-1},k_i)$, choose 
$z\equiv\psi(t_{i-1})$ in \eqv(5.010).

Following the classical pattern of Cramer's type techniques, the lower bound 
will come from `centering the variables' (i.e. introducing a Radon-Nikodym 
factor). 
For a given $\psi\in{\cal G}$  let
$\l_i^*\equiv\l_i^*\left(\sfrac{\psi(t_i)-\psi(t_{i-1})}{t_i-t_{i-1}}\right)$,
$1\leq i\leq n$,
be defined through:
$$
\left(\e\l_i^*,\sfrac{\psi(t_i)-\psi(t_{i-1})}{t_i-t_{i-1}}\right)
-\displaystyle{\cal L}_{\e}(t_{i-1},\psi(t_{i-1}),\e\l_i^*)
= {\cal L}_{\e}^*
\left(t_{i-1},\psi(t_{i-1}),\sfrac{\psi(t_i)-\psi(t_{i-1})}{t_i-t_{i-1}}\right)
\Eq(5.25ter)
$$
Obviously the conditions in 
\eqv(5.25riribis) imply
that $\psi(t_i)\in\inte(\inte_{\e}\L)$
for all $1\leq i\leq n$. The point is that from this, Corollary 4.10, and 
the equivalence $(ii)\Leftrightarrow (iv)$ 
of Lemma 4.7 we can conclude that there exists a positive constant 
$c_0\equiv c_0(\psi)<\infty$ such that:
$$
\max_{1\leq i\leq n}
|\l^*_i|<c_0
\Eq(5.44new)
$$
We then rewrite ${\cal T}_i$ as
$$
{\cal T}_i ={\cal T}_{i,1}{\cal T}_{i,2}
\Eq(5.44)
$$
where
$$
{\cal T}_{i,1}
\equiv
\1_{\{|x(k_{i-1})-\psi(t_{i-1})|\leq\varrho\}}
\sum_{x(k_{i})\in\G_{\e}}
e^{(\l_i^*,x(k_i)-\psi(t_i))}
P_{\e,k_{i-1}}^{(k_{i-1}, k_{i})}(x(k_{i-1}),x(k_{i}))
\Eq(5.45)
$$
$$
\eqalign{
&{\cal T}_{i,2}
\equiv 
\1_{\{|x(k_{i-1})-\psi(t_{i-1})|\leq\varrho\}}
\cr
&
\times\sum_{x(k_{i})\in\G_{\e}}
\frac{e^{(\l_i^*,x(k_i)-\psi(t_i))}
P_{\e,k_{i-1}}^{(k_{i-1}, k_{i})}(x(k_{i-1}),x(k_{i}))
\1_{\{|(x(k_i)-\psi(t_i))+a|\leq\varrho\}}
}
{\sum_{x(k_{i})\in\G_{\e}}
e^{(\l_i^*,x(k_i)-\psi(t_i))}
P_{\e,k_{i-1}}^{(k_{i-1}, k_{i})}(x(k_{i-1}),x(k_{i}))}
e^{-(\l_i^*,x(k_i)-\psi(t_i))}
\cr
}
\Eq(5.46)
$$

We first prove a lower bound for the term
$$
\eqalign{
{\cal T}_{i,3}
\equiv&
\1_{\{|x(k_{i-1})-\psi(t_{i-1})|\leq\varrho\}}
\cr
&\times\sum_{x(k_{i})\in\G_{\e}}
\1_{\{|(x(k_i)-\psi(t_i))+a|\leq\varrho\}}
e^{(\l_i^*,x(k_i)-\psi(t_i))}
P_{\e,k_{i-1}}^{(k_{i-1}, k_{i})}(x(k_{i-1}),x(k_{i}))
}\Eq(5.46bis)
$$
Setting $\ell_i\equiv k_{i}-k_{i-1}$ 
and using \eqv(5.010),
$$
\eqalign{
{\cal T}_{i,3}\geq 
&
e^{-q(\ell_i,|x(k_{i-1})-\psi(t_{i-1})|)}
\1_{\{|x(k_{i-1})-\psi(t_{i-1})|\leq\varrho\}}
\cr
&
\times\sum_{x(k_{i})\in\G_{\e}}
\1_{\{|(x(k_i)-\psi(t_i))+a|\leq\varrho\}}
e^{(\l_i^*,x(k_i)-\psi(t_i))}
\cr
&\,
\times\sum_{\d(1)\in\D}\dots\sum_{\d(\ell_i-1)\in\D}
\prod_{l=1}^{\ell_i}
e^{f^{(0)}_{\e}\left(t_{i-1},\psi(t_{i-1}),\d(l)\right)}
\1_{\{\d(\ell_i)\in\D\}}
\1_{\left\{x(k_i)-x(k_{i-1})=\e\sum_{m=1}^{\ell_i}\d(m)\right\}}
\cr
}
\Eq(5.48)
$$
We have,
$$
\eqalign{
&\1_{\{|x(k_{i-1})-\psi(t_{i-1})|\leq\varrho\}}
\1_{\left\{x(k_i)-x(k_{i-1})=\e\sum_{m=1}^{\ell_i}\d(m)\right\}}
e^{(\l_i^*,x(k_i)-\psi(t_i))}
\cr
\geq &
\1_{\{|x(k_{i-1})-\psi(t_{i-1})|\leq\varrho\}}
\1_{\left\{x(k_i)-x(k_{i-1})=\e\sum_{m=1}^{\ell_i}\d(m)\right\}}
e^{-\varrho|\l_i^*|
+(\e\l_i^*,\sum_{m=1}^{\ell_i}\d(m))
-(\l_i^*,\psi(t_i)-\psi(t_{i-1}))}
\cr
}
\Eq(5.49)
$$
Consequently,
$$
\eqalign{
{\cal T}_{i,3}\geq 
&
e^{-q(\ell_i,\varrho)-\varrho|\l_i^*|}
\1_{\{\left|x(k_{i-1})-\psi(t_{i-1})\right|\leq\varrho\}}
\cr
&\,\,\,\,\,\,
\times\sum_{\d(1)\in\D}\dots\sum_{\d(\ell_i)\in\D}
\prod_{l=1}^{\ell_i}e^{(\e\l_i^*,\d(l))}
e^{f^{(0)}_{\e}\left(t_{i-1},\psi(t_{i-1}),\d(l)\right)}
\cr
&\,\,\,\,\,\,\,\,\,\,\,
\times\1_{\left\{\left|
\e\sum_{m=1}^{\ell_i}\d(m)
-(\psi(t_i)-\psi(t_{i-1}))+(x(k_{i-1})-\psi(t_{i-1}))+a
\right|\leq\varrho\right\}}
\cr
}
\Eq(5.49bis)
$$
The same arguments applied to ${\cal T}_{i,1}$ give
$$
\eqalign{
{\cal T}_{i,1}\geq 
&
e^{-q(\ell_i,\varrho)-\varrho|\l_i^*|}
\1_{\{|x(k_{i-1})-\psi(t_{i-1})|\leq\varrho\}}\cr
&\,\,\,\,\,\,\,\,\,\,\,\,\,\,\,\,\,\,
\times e^{ -(\l_i^*,\psi(t_i)-\psi(t_{i-1}))}
\prod_{l=1}^{\ell_i}
\sum_{\d(l)\in\D}e^{(\e\l_i^*,\d(l))}
e^{f^{(0)}_{\e}\left(t_{i-1},\psi(t_{i-1}),\d(l)\right)}
\cr
= &
e^{-q(\ell_i,\varrho)-\varrho|\l_i^*|}
\1_{\{|x(k_{i-1})-\psi(t_{i-1})|\leq\varrho\}}
e^{-\ell_i\left\{
\left(\e\l_i^*,\sfrac{\psi(t_i)-\psi(t_{i-1})}{t_i-t_{i-1}}\right)
-\displaystyle{\cal L}_{\e}(t_{i-1},\psi(t_{i-1}),\e\l_i^*)
\right\}}
\cr
}
\Eq(5.50)
$$
and, by definition of $\l_i^*$,
$$
{\cal T}_{i,1}\geq 
e^{-q(\ell_i,\varrho)-\varrho|\l_i^*|}
\1_{\{|x(k_{i-1})-\psi(t_{i-1})|\leq\varrho\}}
e^{-\e^{-1}(t_i-t_{i-1})\displaystyle{\cal L}_{\e}^*
\left(t_{i-1},\psi(t_{i-1}),\sfrac{\psi(t_i)-\psi(t_{i-1})}{t_i-t_{i-1}}
\right)
}
\Eq(5.51)
$$
which is precisely the form of the bound we need.

We now turn to the term ${\cal T}_{i,2}$ and first write
$$
\eqalign{
{\cal T}_{i,2}
\geq &
\1_{\{|x(k_{i-1})-\psi(t_{i-1})|\leq\varrho\}}e^{-\varrho|\l_i^*|}
\cr
&
\times\frac{
\sum_{x(k_{i})\in\G_{\e}}
e^{(\l_i^*,x(k_i)-\psi(t_i))}
P_{\e,k_{i-1}}^{(k_{i-1}, k_{i})}(x(k_{i-1}),x(k_{i}))
\1_{\{|(x(k_i)-\psi(t_i))+a|\leq\varrho\}}
}
{\sum_{x(k_{i})\in\G_{\e}}
e^{(\l_i^*,x(k_i)-\psi(t_i))}
P_{\e,k_{i-1}}^{(k_{i-1}, k_{i})}(x(k_{i-1}),x(k_{i}))
}
\cr
}
\Eq(5.52)
$$
\eqv(5.49bis) allows to bound the numerator in \eqv(5.52) from above.
Virtually the same arguments allow to bound the 
denominator from above:
$$
\eqalign{
&{\sum_{x(k_{i})\in\G_{\e}}
e^{(\l_i^*,x(k_i)-\psi(t_i))}
P_{\e,k_{i-1}}^{(k_{i-1}, k_{i})}(x(k_{i-1}),x(k_{i}))
}
\cr&\leq
e^{\{q(\ell_i,\varrho)+
\varrho|\l_i^*|\}}
\prod_{l=1}^{\ell_i}{\sum_{\d(l)\in\D}
e^{(\e\l_i^*,\d(l))+f^{(0)}_{\e}(t_{i-1},\psi(t_{i-1}),\d(l))}}
}\Eq(5.34) 
$$
Combining these yields
$$
\eqalign{
{\cal T}_{i,2}\geq 
&
e^{-\{2q(\ell_i,\varrho)
+3\varrho|\l_i^*|\}}
\1_{\{\left|x(k_{i-1})-\psi(t_{i-1})\right|\leq\varrho\}}
\cr
&\,\,\,\,\,\,\,\,\,\,\,\,\,\,
\times\sum_{\d(1)\in\D}\dots\sum_{\d(\ell_i)\in\D}
\prod_{l=1}^{\ell_i}
\frac
{e^{(\e\l_i^*,\d(l))+f^{(0)}_{\e}(t_{i-1},\psi(t_{i-1}),\d(l))}}
{\sum_{\d(l)\in\D}
e^{(\e\l_i^*,\d(l))+f^{(0)}_{\e}(t_{i-1},\psi(t_{i-1}),\d(l))}}
\cr
&\,\,\,\,\,\,\,\,\,\,\,\,\,\,\,\,\,\,\,\,\,\,\,\,\,\,\,\,\,\,\,\,
\times\1_{\left\{\left|
\e\sum_{m=1}^{\ell_i}\d(m)
-(\psi(t_i)-\psi(t_{i-1}))+(x(k_{i-1})-\psi(t_{i-1}))+a
\right|\leq\varrho\right\}}
\cr
}
\Eq(5.54)
$$
At this point \eqv(5.54) may be recast in the following form:
let $\left\{\chi_{m,i}\right\}_{1\leq m\leq\ell_i}$ be a family of 
i.i.d. r.v.'s taking values in $\D$ with law, $\nu_i$, defined through
(see \eqv(4.7))
$$
\nu_i(\d)
\equiv
\nu_{\e,t_{i-1},\psi(t_{i-1})}^{\l_i^*}(\d)
=
\frac
{e^{(\e\l_i^*,\d)+f^{(0)}_{\e}(t_{i-1},\psi(t_{i-1}),\d)}}
{\sum_{\d\in\D}
e^{(\e\l_i^*,\d)+f^{(0)}_{\e}(t_{i-1},\psi(t_{i-1}),\d)}},
\,\,\,\,\,\,\,\forall\d\in\D
\Eq(5.55)
$$
Set
$$
\overline\chi_{m,i}=\chi_{m,i}-\sfrac{\psi(t_i)-\psi(t_{i-1})}{t_i-t_{i-1}}
\Eq(5.56)
$$
$$
S_i=\sum_{m=1}^{\ell_i}\overline\chi_{m,i}
\Eq(5.57)
$$
and let $\E_{\{\nu_i\}}$ denote the expectation w.r.t. 
$\left\{\chi_{m,i}\right\}$. Then \eqv(5.54) reads,
$$
\eqalign{
{\cal T}_{i,2}\geq &
e^{-\{2q(\ell_i,\varrho)
+3\varrho|\l_i^*|\}}
\1_{\{\left|x(k_{i-1})-\psi(t_{i-1})\right|\leq\varrho\}}
\E_{\{\nu_i\}}
\1_{\left\{\left|
\e S_i +(x(k_{i-1})-\psi(t_{i-1}))+a
\right|\leq\varrho\right\}}
\cr
}
\Eq(5.58)
$$
Collecting \eqv(5.44), \eqv(5.51) and \eqv(5.58) we thus obtain
$$
\eqalign{
{\cal T}_{i}&\geq 
e^{-\varsigma_i-\e^{-1}(t_i-t_{i-1})\displaystyle{\cal L}_{\e}^*
\left(t_{i-1},\psi(t_{i-1}),\sfrac{\psi(t_i)-\psi(t_{i-1})}{t_i-t_{i-1}}
\right)
}
\cr
&\,\,\,\,\,\,\,\,\,
\times\1_{\{|x(k_{i-1})-\psi(t_{i-1})|\leq\varrho\}}
\E_{\{\nu_i\}}
\1_{\left\{\left|
\e S_i +(x(k_{i-1})-\psi(t_{i-1}))+a
\right|\leq\varrho\right\}}
\cr
}
\Eq(5.59)
$$
where
$$
\varsigma_i\equiv 3q(\ell_i,\varrho)
+4\varrho|\l_i^*|
\Eq(5.60)
$$
We are now in a position to deal with the r.h.s. of \eqv(5.27). Applying
\eqv(5.59) to ${\cal T}_{n}$ gives rise to a term of the form
${\cal T'}_{n-1}$ (see definition \eqv(5.41)) with 
$a_{n-1,n-1}=0$ and $a_{n-1,n}=\e S_n$. The second iteration step thus yields
$$
\eqalign{
&\1_{\{|x(k_{n-2})-\psi(t_{n-2})|\leq\varrho\}}
\sum_{x(k_{n-1})\in\G_{\e}}
\1_{\{|x(k_{n-1})-\psi(t_{n-1})|\leq\varrho\}}
P_{\e,k_{n-2}}^{(k_{n-2}, k_{n-1})}(x(k_{n-2}),x(k_{n-1}))\cr
&\,\,\,\,\,\,\,\,
\times\sum_{x(k_{n})\in\G_{\e}}
\1_{\{|x(k_n)-\psi(t_n)|\leq\varrho\}}
P_{\e,k_{n-1}}^{(k_{n-1}, k_{n})}(x(k_{n-1}),x(k_{n}))
\cr
\geq &
e^{-(\varsigma_{n}+\varsigma_{n-1})
-\e^{-1}\sum_{i=n-1}^n (t_i-t_{i-1})\displaystyle{\cal L}_{\e}^*
\left(t_{i-1},\psi(t_{i-1}),\sfrac{\psi(t_i)-\psi(t_{i-1})}{t_i-t_{i-1}}
\right)
}
\1_{\{|x(k_{n-2})-\psi(t_{n-2})|\leq\varrho\}}\cr
&\,\,\times\E_{\{\nu_{n-1}\}}
\1_{\left\{\left|
\e S_{n-1} +(x(k_{n-2})-\psi(t_{n-2}))
\right|\leq\varrho\right\}}
\E_{\{\nu_{n}\}}\1_{\left\{\left|
\e(S_{n-1}+S_{n}) +(x(k_{n-2})-\psi(t_{n-2}))
\right|\leq\varrho\right\}}
\cr
}
\Eq(5.61)
$$
and gradually, setting
$$
a_{i,j}=
\cases
0 & \hbox{if} \,\, j=i\cr
\e(S_{j+1}+\dots+S_{n})& \hbox{if}\,\, i+1\leq j\leq n\cr
\endcases
\Eq(5.61bis)
$$
in \eqv(5.41) at step $i$, we obtain,
$$
\eqalign{
&{\cal P}_{\e,\phi_0}
\left(\max_{0\leq i\leq n}\left|X(\sfrac{t_i}{\e})
-\psi(t_i)\right|\leq\varrho\right)
\cr&\geq 
e^{-\sfrac{1}{\e}\widetilde Q-\sfrac{1}{\e}\sum_{i=1}^n(t_i-t_{i-1})
\displaystyle{\cal L}_{\e}^*
\left(t_{i-1},\psi(t_{i-1}),\sfrac{\psi(t_i)-\psi(t_{i-1})}{t_i-t_{i-1}}
\right)}
\cr
&\,\,\,\,\times\sum_{x(t_0)\in\G_{\e}}
\pi_{\e,\phi_0}(x(k_0))
\1_{\{|x(k_0)-\psi(t_0)|\leq\varrho\}}
\cr
&\,\,\,\,\times\E_{\{\nu_{1}\}}
\1_{\left\{\left|
\e S_{1} +(x(k_{0})-\psi(t_{0}))
\right|\leq\varrho\right\}}
\dots
\E_{\{\nu_{n}\}}\1_{\left\{\left|
\e(S_{1}+\dots+S_{n}) +(x(k_{0})-\psi(t_{0}))
\right|\leq\varrho\right\}}
\cr
&=
{\cal R}e^{-\sfrac{1}{\e}\widetilde Q-\sfrac{1}{\e}\sum_{i=1}^n(t_i-t_{i-1})
\displaystyle{\cal L}_{\e}^*
\left(t_{i-1},\psi(t_{i-1}),\sfrac{\psi(t_i)-\psi(t_{i-1})}{t_i-t_{i-1}}
\right)}
\cr
}
\Eq(5.62)
$$
where
$$
\widetilde Q\equiv \e\sum_{i=1}^n\varsigma_i
\Eq(5.63)
$$
$$
{\cal R}\equiv\sum_{x(t_0)\in\G_{\e}}
\pi_{\e,\phi_0}(x(k_0))
\1_{\{|x(k_0)-\psi(t_0)|\leq\varrho\}}
\E_{\{\nu\}}
\1_{\left\{\left|\e(S_{1}+\dots+S_{n}) +(x(k_{0})-\psi(t_{0}))
\right|\leq\varrho\right\}}
\Eq(5.64)
$$
and $\E_{\{\nu\}}$ denotes the expectation w.r.t. the joint law
of $\{S_i\}_{1\leq i\leq n}$.
We are left to estimate ${\cal R}$. Assume that
$\varrho\geq\e\sqrt{d}$. Then
$$
\eqalign{
{\cal R}
&\geq
\sum_{x(t_0)\in\G_{\e}}
\pi_{\e,\phi_0}(x(k_0))
\1_{\{|x(k_0)-\psi(t_0)|\leq\e\sqrt{d}\}}
\E_{\{\nu\}}
\1_{\left\{\left|S_{1}+\dots+S_{n}\right|
\leq\e^{-1}(\varrho-\e\sqrt{d})\right\}}
\cr
&=\sfrac{1}{|\{x\in\G_{\e}\mid |x-\phi_0|\leq \e\sqrt{d}\}|}
\sum_{x(t_0):|x(k_0)-\phi_0|\leq \e\sqrt{d}}
\1_{\{|x(k_0)-\psi(t_0)|\leq\e\sqrt{d}\}}
\E_{\{\nu\}}
\1_{\left\{\left|S_{1}+\dots+S_{n}
\right|\leq\e^{-1}(\varrho-\e\sqrt{d})\right\}}
\cr
&\geq
\sfrac{1}{4d^2+1}
\1_{\{|x(k_0)-\psi(t_0)|\leq\e\sqrt{d}\}}
\E_{\{\nu\}}
\1_{\left\{\left|S_{1}+\dots+S_{n}
\right|\leq\e^{-1}(\varrho-\e\sqrt{d})\right\}}
\cr
}
\Eq(5.65)
$$
for any $x(k_0)\in\{x\in\G_{\e}\mid |x-\phi_0|\leq \e\sqrt{d}\}$. Since 
$$
\bigcup_{x(k_0)\in\{x\in\G_{\e}\mid |x-\phi_0|\leq \e\sqrt{d}\}}
\{y\in\R^d\mid |y-x(k_0)|\leq \e\sqrt{d}\}\supset 
\{y\in\R^d\mid |y-\phi_0|\leq \e\sqrt{d} \}
\Eq(5.66)
$$
then
$$
{\cal R}
\geq
\cases
\sfrac{1}{4d^2+1}
\E_{\{\nu\}}
\1_{\left\{\left|S_{1}+\dots+S_{n}
\right|\leq\e^{-1}(\varrho-\e\sqrt{d})\right\}}
& \hbox{if} \,\, |\psi(t_0)-\phi_0|\leq \e\sqrt{d}\cr
& \cr
0 & \hbox{otherwise}\cr
\endcases
\Eq(5.67)
$$
and it remains to estimate the expectation. But this is immediate once
observed that, recalling \eqv(5.25ter) and 
combining Lemma 4.4, (iii), together with the equivalence
$(i)\Leftrightarrow(iii)$ of Lemma 4.7 we have, for all $1\leq m\leq \ell_i$,
$$
\E_{\nu_i}\chi_{m,i}=
\nabla
\Phi_{\e,t_{i-1},\psi(t_{i-1})}(v)|_{v=\e\l_i^*}
=\sfrac{\psi(t_i)-\psi(t_{i-1})}{t_i-t_{i-1}}
\Eq(5.68)
$$
and
$$
\eqalign{
&\E_{\nu_i}\overline\chi_{m,i}=0\cr
&\E_{\nu_i}\left|\overline\chi_{m,i}\right|^2=
\D\Phi_{\e,t_{i-1},\psi(t_{i-1})}(v)|_{v=\e\l_i^*}
\cr
}
\Eq(5.69)
$$
Defining
$$
\s^2\equiv\sigma^2(\{\psi(t_{i})\},\{\l^*_i\})
=T\max_{1\leq i\leq n}
\D
\Phi_{\e,t_{i-1},\psi(t_{i-1})}(v)|_{v=\e\l_i^*}
\Eq(5.26)
$$
Moreover,
$$
\s^2\leq T(\diam\D)^2
\Eq(5.68bis)
$$
Hence, by independence and Chebyshev's inequality
$$
\eqalign{
\E_{\{\nu\}}
\1_{\left\{\left|S_{1}+\dots+S_{n}
\right|\leq\e^{-1}(\varrho-\e\sqrt{d})\right\}}
&=
1-\E_{\{\nu\}}
\1_{\left\{\left|S_{1}+\dots+S_{n}
\right|>\e^{-1}(\varrho-\e\sqrt{d})\right\}}
\cr
&\geq 1-\left(\e(\varrho-\e\sqrt{d})^{-1}\right)^2
\E_{\{\nu\}}(S_{1}+\dots+S_{n})^2
\cr
&\geq 1-\left(\e(\varrho-\e\sqrt{d})^{-1}\right)^2
\sum_{i=1}^n \ell_i 
\D\Phi_{\e,t_{i-1},\psi(t_{i-1})}(v)|_{v=\e\l_i^*}
\cr
&\geq 1-\e\left(\varrho-\e\sqrt{d}\right)^{-2} 
\sigma^2(\{\psi(t_{i})\},\{\l^*_i\})
\cr
&\geq 1-\e T(\diam\D)^2\left(\varrho-\e\sqrt{d}\right)^{-2} 
\cr
&\geq\frac{1}{2}
}
\Eq(5.70)
$$
whenever
$\varrho\geq \sqrt{2\e T}\diam\D +\e\sqrt{d}$.
For such a $\varrho$, inserting \eqv(5.70) in \eqv(5.67) and 
combining with \eqv(5.62) proves Lemma \ver.5 since 
$\widetilde Q\leq 
Q\left(\bar{\cal S}_{\rho,\gamma/2}(\phi),\zeta,
\max_{1\leq i\leq n}|\e\l_i^*|\right)$
and since by \eqv(5.44new),
$$
\sup_{\psi\in{\cal G}}Q\left(
\bar{\cal S}_{\rho,\gamma/2}(\phi),
\zeta, \max_{1\leq i\leq n}|\e\l_i^*|\right)
\leq  Q(\bar{\cal S}_{\rho,\gamma/2}(\phi),\zeta,c_0)
\Eq(5.70bis)
$$
(see definitions \eqv(5.1), \eqv(5.09), and \eqv(5.44new) as well as 
\eqv(5.60) and \eqv(5.63) for the first of the last two inequalities).
\endproof

\vskip.3truecm
\noindent{\bf 5.3: Proof of Proposition 3.2 (concluded).}
\vskip.3truecm

To conclude the proofs of the upper and lower bounds, we need the following two
lemmata that will permit to replace the sums over $t_i$ by integrals.

\lemma{\ver.6}{\it Recall that $D=\conv\D$ and define the sets
$$
\eqalign{
{\cal K}([0,T])&=\left\{
\psi\in W([0,T])\, \Big|\,
\dot\phi(t)\in D,\, \hbox{for Lebesgue a.e.}\,\, t\in[0,T]
\right\}
\cr
{\cal K}^{\circ}([0,T])&=
\left\{
\psi\in W([0,T])\, \Big|\,
\dot\phi(t)\in \ri\,D,\, \hbox{for Lebesgue a.e.}\,\, t\in[0,T]
\right\}
\cr
}
\Eq(5.80)
$$
With ${\cal E}([0,T])$ and ${\cal E}^{\circ}([0,T])$ defined respectively in
\eqv(3.1) and \eqv(5.25riter) we have:
$$
\eqalign{
{\cal K}([0,T])&={\cal E}([0,T])\cr
{\cal K}^{\circ}([0,T])&\subset {\cal E}^{\circ}([0,T])\cr
}
\Eq(5.81)
$$
}

\proof The proof is elementary. 
Recall that by assumption  $D$ is a bounded  closed and convex subset of
$\R^d$. For any bounded convex subset $A$ in $\R^d$ and any 
$\psi\in{\cal C}([0,T])$ consider the following three conditions:
{\obeylines{
{$(i)$} $\psi\in L^1([0,T])$ and $\dot\psi(t)\in A$ for Lebesgue a.e. $t\in[0,T]$.
{$(ii)$} $\psi\in L^1([0,T])$ and $\sfrac{1}{t-t'}\int_{t'}^{t}ds\dot\psi(s)\in A$ $\forall t\in[0,T]$, $\forall t'\in[0,T]$, $t\neq t'$.
{$(iii)$} $\sfrac{\psi(t)-\psi(t')}{t-t'}\in A$ $\forall t\in[0,T]$, $\forall t'\in[0,T]$, $t\neq t'$.
}}
Then the following conclusions hold:
{\obeylines{
{$(iv)$} If $A=D$ or if $A=\ri\,D$ then $(ii)\Leftrightarrow (iii)$
{$(v)$}  If $A=D$ or if $A=\ri\,D$ then $(i)\Rightarrow (ii)$
{$(vi)$} If $A=D$ then $(ii)\Leftrightarrow (i)$
}}

We first prove {$(iv)$}: that $(ii)\Rightarrow (iii)$ is immediate whereas
since $A$ is bounded $\psi$ is Lipshitz and, in particular, absolutely 
continuous, yielding $(iii)\Rightarrow (ii)$. Whenever $A$ is 
a closed or opened set, the implication
$(i)\Rightarrow (ii)$ results from it's convexity and the integrability of 
$\dot\psi$: this proves $(v)$. If in addition $A$ is closed then,
by a standard result of real analysis,
$(ii)\Rightarrow (i)$ (see e.g. [Ru], Theorem 1.40);
this together with  $(v)$ yields $(vi)$.
Now $(iv)$ together with $(vi)$ implies the first relation in \eqv(5.81) 
while {$(iv)$} together with {$(v)$} implies the second.
The proof is done.\endproof

\lemma{\ver.7} {\it  Let $\SS$ be any closed bounded subset of 
$\inte(\inte_{\e}\L)$, and let $t_i$, $i=1,\dots,n$ be as in Lemma \ver.4.
\item{(i)} If  $\psi$ is in
$$ 
\left\{
\psi\in{\cal E}([0,T])\,\Big|\, \psi(t)\in{\cal S},\,\,\,\,\forall
t\in[0,T]\right\}
\Eq(5.100)
$$
then, for each $\varepsilon_0>0$
there corresponds $\varepsilon_1>0$ such that if $\e\tau<\varepsilon_1$,
$$
\eqalign{
&
\left|{\sum_{i=1}^n} (t_i-t_{i-1}){\cal L}_{\e}^*
\left(t_{i-1},\psi(t_{i-1}),\sfrac{\psi(t_i)-\psi(t_{i-1})}{t_i-t_{i-1}}\right)
-
\int_0^T dt {\cal L}_{\e}^*(t,\psi(t),\dot\psi(t))\right|
\cr
\leq &
\varepsilon_0 T+(\theta+\vartheta({\cal S})\diam\D)n\sfrac{(\e\tau)^2}{2}
}
\Eq(5.101)
$$
\item{(ii)} Let $t_i$, $i=0,\dots,n$, $n$, $\bar\z$ and $r$ be given  as in 
Lemma \ver.4. 
Assume that $\psi'(t_i)\in\R^d$ are such that 
$$
|\psi'(t_i)-\psi'(t_{i-1})|\leq |t_i-t_{i-1}| C, \quad \forall {i=1,\dots,n}
\Eq(5.200)
$$ 
for some constant $0<C<\infty$ and 
$$
\dist\left(\psi'(t_i),\L\right)\leq \bar\zeta
\Eq(5.2001)
$$
Let $\wt\psi(t)$, $t\in [0,T]$ be the linear interpolation of the 
points $\psi'(t_i)$.   
Then, for each $\varepsilon_0>0$
there exists $\varepsilon_1>0$ (depending on $r$ and $C$) such that,
 if $\e\tau<\varepsilon_1$,
$$
{\sum_{i=1}^n} (t_i-t_{i-1}){\cal L}_{\e}^{(r)*}
\left(t_{i-1},\psi'(t_{i-1}),
\sfrac{\psi'(t_i)-\psi'(t_{i-1})}{t_i-t_{i-1}}\right)
-
\int_0^T dt {\cal L}_{\e}^{(r)*}(t,\psi'(t),\dot\psi'(t))
\geq
-3\varepsilon_0 T
\Eq(5.201)
$$
}

\proof We first prove (i). Recall that
 $\Phi^*_{\e,t_{i-1},\psi(t_{i-1})}(\cdot)=
{\cal L}^*_{\e}(t_{i-1},\psi(t_{i-1}),\cdot)$
and
$\tau\equiv\max_{0\leq i\leq n}\e^{-1}|t_{i+1}-t_i|$ as defined in
\eqv(4.3) and \eqv(5.1). Let us write:
$$
\eqalign{
&
(t_i-t_{i-1}){\cal L}_{\e}^*
\left(t_{i-1},\psi(t_{i-1}),\sfrac{\psi(t_i)-\psi(t_{i-1})}{t_i-t_{i-1}}\right)
\cr
= &
\int_{t_{i-1}}^{t_i}ds
 {\cal L}_{\e}^*(s,\psi(s),\dot\psi(s))
\cr
& +
\left[\int_{t_{i-1}}^{t_i}ds
\left(
\Phi^*_{\e,t_{i-1},\psi(t_{i-1})}
\left(
{\textstyle{
\sfrac{1}{t_i-t_{i-1}}\int_{t_{i-1}}^{t_i}ds'\dot\psi(s')
}}
\right)
-
\Phi^*_{\e,t_{i-1},\psi(t_{i-1})}\left(\dot\psi(s)\right)
\right)\right]
\cr
& +
\left\{\int_{t_{i-1}}^{t_i}ds
\left(
{\cal L}_{\e}^*\left(t_{i-1},\psi(t_{i-1}),\dot\psi(s)\right)
-
{\cal L}_{\e}^*\left(s,\psi(s),\dot\psi(s)\right)\right)
\right\}\cr
= &
\int_{t_{i-1}}^{t_i}ds
{\cal L}_{\e}^*(s,\psi(s),\dot\psi(s))
+\left[I_i\right]+\left\{J_i\right\}
\cr
}
\Eq(5.103)
$$
where the last line defines the terms $I_i$ and $J_i$. 
In order to bound $J_i$ we use the decomposition
$$
\eqalign{
&
{\cal L}_{\e}^*\left(t_{i-1},\psi(t_{i-1}),\dot\psi(s)\right)
-
{\cal L}_{\e}^*\left(s,\psi(s),\dot\psi(s)\right)
\cr
=
&
\left[
{\cal L}_{\e}^*\left(t_{i-1},\psi(t_{i-1}),\dot\psi(s)\right)
-
{\cal L}_{\e}^*\left(s,\psi(t_{i-1}),\dot\psi(s)\right)
\right]
\cr+&
\left[
{\cal L}_{\e}^*\left(s,\psi(t_{i-1}),\dot\psi(s)\right)
-
{\cal L}_{\e}^*\left(s,\psi(s),\dot\psi(s)\right)
\right]
\cr
}
\Eq(5.104bis)
$$
and, applying Lemma 4.9, obtain 
$$
\eqalign{
\left|
{\cal L}_{\e}^*\left(t_{i-1},\psi(t_{i-1}),\dot\psi(s)\right)
-
{\cal L}_{\e}^*\left(s,\psi(s),\dot\psi(s)\right)\right|
\leq &
\theta|s-t_{i-1}|+\vartheta|\psi(s)-\psi(t_{i-1})|
\cr
\leq &
(\theta+\vartheta\diam\D)|s-t_{i-1}|
\cr
}
\Eq(5.104)
$$ 
where $\vartheta\equiv\vartheta({\cal S})$. Thus,
$$
|J_i|\leq (\theta+\vartheta\diam\D)\int_{t_{i-1}}^{t_i}ds|s-t_{i-1}|
=(\theta+\vartheta\diam\D)\sfrac{(t_i-t_{i-1})^2}{2}
\leq (\theta+\vartheta\diam\D)\sfrac{(\e\tau)^2}{2}
\Eq(5.105)
$$
We now bound $I_i$.
By Lemma 4.6, (i), $\Phi^*_{\e,t_{i-1},\psi(t_{i-1})}$ is convex and
lower semi-continuous. Convexity implies $I_i\leq 0$. For an upper 
bound note first that
by Lebesgue's Theorem: to each $\varepsilon_2>0$ there 
corresponds $\varepsilon_1>0$ such that, for Lebesgue almost every 
$s\in[t',t]$,
$$
\left|\int_{t'}^{t}ds'\dot\psi(s')-\dot\psi(s)\right|
<\varepsilon_2 |t'-t|
\Eq(5.106)
$$
for all $[t',t]\subset[0,T]$ verifying $s\in[t',t]$ and $|t-t'|<\varepsilon_1$
\note{the set of $s$'s for which \eqv(5.104) holds is usually called the 
Lebesgue set of $\psi$.}. Next, by definition of lower semi-continuity,
for any $x\in\R^d$ we have: to each 
$\varepsilon_0>0$ there corresponds $\varepsilon_2>0$ such that
if $|x-y|<\varepsilon_2$, then
$\Phi^*_{\e,t_{i-1},\psi(t_{i-1})}(x)\geq
\Phi^*_{\e,t_{i-1},\psi(t_{i-1})}(y)-\varepsilon_0$.
Thus, for each $\varepsilon_0>0$,
if $\e$ is sufficiently small so that $\e\tau<\varepsilon_1$ we have,
 on the Lebesgue set of $\psi$:
$$
\Phi^*_{\e,t_{i-1},\psi(t_{i-1})}
\left(
{\textstyle{
\sfrac{1}{t_i-t_{i-1}}\int_{t_{i-1}}^{t_i}ds'\dot\psi(s')
}}
\right)
\geq 
\Phi^*_{\e,t_{i-1},\psi(t_{i-1})}\left(\dot\psi(s)\right)
-\varepsilon_0
\Eq(5.107)
$$
and
$$
I_i\geq -(t_i-t_{i-1})\varepsilon_0
\Eq(5.108)
$$
Inserting our bounds on $I_i$ and $J_i$ in \eqv(5.103) and adding up  yields 
$$
\eqalign{
&
\left|{\sum_{i=1}^n} (t_i-t_{i-1}){\cal L}_{\e}^*
\left(t_{i-1},\psi(t_{i-1}),\sfrac{\psi(t_i)-\psi(t_{i-1})}{t_i-t_{i-1}}\right)
-
\int_0^{\e\left[\frac T\e\right]} 
dt {\cal L}_{\e}^*(t,\psi(t),\dot\psi(t))\right|
\cr
\leq &
\varepsilon_0 T+(\theta+\vartheta({\cal S})\diam\D)n\sfrac{(\e\tau)^2}{2}
}
\Eq(5.101fast)
$$
But 
$
\left|\int_{\e\left[\frac T\e\right]}^T 
dt {\cal L}_{\e}^*(t,\psi(t),\dot\psi(t))\right|
\leq \e const(\SS)
$
so that 
\eqv(5.101) obtains upon minor modification of $\varepsilon_0$.

To prove  (ii) 
we note that since $\wt\psi$ is linear between the points $t_i$,
in the analogue of \eqv(5.103) the term corresponding to $[I_i]$ is absent,
 i.e. we have
$$
\eqalign{
&
(t_i-t_{i-1}){\cal L}_{\e}^{(r)*}
\left(t_{i-1},\psi'(t_{i-1}),\sfrac{\psi'(t_i)-\psi'(t_{i-1})}{t_i-t_{i-1}}
\right)
\cr
= &
\int_{t_{i-1}}^{t_i}ds
 {\cal L}_{\e}^{(r)*}(s,\wt\psi(s),\dot{\wt\psi}(s))
\cr
& +
\int_{t_{i-1}}^{t_i}ds
\left(
{\cal L}_{\e}^{(r)*}\left(t_{i-1},\wt\psi(t_{i-1}),\dot{\wt\psi}(s)\right)
-
{\cal L}_{\e}^*\left(s,\wt\psi(s),\dot{\wt\psi}(s)\right)\right)
\cr
}
\Eq(5.103bis)
$$
To bound the second term in \eqv(5.103bis) we use the same decomposition as in 
\eqv(5.104bis). However, instead of the Lipshitz bounds \eqv(5.104) we use the
lower semi-continuity property 
of ${\cal L}_{\e}^{(r)*}$ (see Lemma 4.12) 
together with the  fact that  $\wt\psi$ is Lipshitz by \eqv(5.200),
it follows from the decomposition
\eqv(5.104bis) that: 
for each $\varepsilon_0$
there corresponds $\varepsilon_1'>0$ such that if $\e\tau<\varepsilon_1'$,
$$
{\cal L}_{\e}^{(r)*}\left(t_{i-1},\psi(t_{i-1}),\dot\psi(s)\right)
-
{\cal L}_{\e}^{(r)*}\left(s,\psi(s),\dot\psi(s)\right)
\geq -2\varepsilon_0
\Eq(5.104sup)
$$
The lemma is proven.
\endproof

\proofof{the lower bound \eqv(3.5bis)}: Given any $\g>0$ we may choose
$\zeta$ and $\tau$ depending 
on $\e$ in such a way that firstly, both $\zeta\downarrow 0$ and 
$\e\tau\downarrow 0$ as $\e\downarrow 0$ (hence $\eta\downarrow 0$ as 
$\e\downarrow 0$), and secondly, that the conditions \eqv(5.25riribis) 
of Lemma \ver.5 as well as those of Lemma \ver.2, (ii), are satisfied.  
It then easily follows from the first relation of Lemma 
\ver.6 that
$$
\bigcup_{\g>0}\bigcup_{\e>0}
\BB_{\rho-2(\eta+\e\tau\diam\D),\g}(\phi)
=\BB_{\rho}(\phi)\cap{\DD}^{\circ}([0,T])
\Eq(500.1)
$$
Setting 
$$
\eqalign{
\widetilde\GG&\equiv
\BB_{\rho}(\phi)\cap{\DD}^{\circ}([0,T])\cap
{\cal E}^{\circ}([0,T])\cr
\GG&\equiv\BB_{\rho}(\phi)\cap{\DD}^{\circ}([0,T])\cap{\cal K}^{\circ}([0,T])
\cr
}
\Eq(500.2)
$$
and using now the second relation of Lemma \ver.6, we moreover have
$\GG\subset\widetilde\GG$.
Let $\psi$ be any path in $\widetilde\GG$. Then obviously, 
$\exists\g_0>0$ s.t. $\forall 0<\g<\g_0$ 
$\exists 0<\e_0$ s.t. $\forall \e<\e_0$, 
$\psi\in\BB_{\rho-2(\eta+\e\tau\diam\D),\g}(\phi)
\cap{\cal E}^{\circ}([0,T])$. Thus, given $\g<\g_0$ and $\e<\e_0$ we may 
combine the bound \eqv(5.26bis) of 
 Lemma \ver.5 and Lemma \ver.7, (i), to write, 
under the assumptions of Lemma \ver.7, (i), and choosing 
$\SS\equiv\bar{\cal S}_{\rho,\gamma/2}(\phi)$ therein,
$$
\e\log{\cal P}_{\e,\phi_0}
\left(\max_{0\leq i\leq n}\left|X(\sfrac{t_i}{\e})
-\psi(t_i)\right|\leq\zeta\right)
\geq
-\int_0^T dt {\cal L}_{\e}^*(t,\psi(t),\dot\psi(t))
-\widetilde Q
\left(\varepsilon_0,\bar{\cal S}_{\rho,\gamma/2}(\phi),\zeta, c_0\right)
\Eq(500.3)
$$
where 
$$
\widetilde Q
\left(\varepsilon_0,\bar{\cal S}_{\rho,\gamma/2}(\phi),\zeta, c_0\right)
\equiv Q\left(\bar{\cal S}_{\rho,\gamma/2}(\phi),\zeta, c_0\right)
+\varepsilon_0 T+(\theta
+\vartheta(\bar{\cal S}_{\rho,\gamma/2}(\phi))\diam\D)n\sfrac{(\e\tau)^2}{2}
\Eq(500.4)
$$
Making use of  Lemma \ver.2, (ii), \eqv(500.3) entails
$$
\e\log\widetilde{\cal P}_{\e,\phi_0}(\BB_{\rho}(\phi))\geq
-\int_0^T dt {\cal L}_{\e}^*(t,\psi(t),\dot\psi(t))
-\widetilde Q
\left(\varepsilon_0,\bar{\cal S}_{\rho,\gamma/2}(\phi),\zeta, c_0\right)
\Eq(500.5)
$$

The next step consists in taking the limit as $\e\downarrow 0$. This will be 
done with the help of the following two observations. On the one hand,
by Lemma 4.5, $\LL^*_{\e}$ is positive and bounded on
$\R^+\times\inte_{\e}\L\times(\conv\D)$. Since, for all 
$\e$ sufficiently small, $\psi(t)$ is contained for all $0\leq t\leq T$
in  a compact subset of $\inte(\inte_\e\L)$, we have, by Lemma 4.9 (v) that
$\LL^*_{\e}(t,\psi(t),\dot\psi(t))$ converges uniformly in $t\in [0,T]$.
Hence, (for each $0<\g<\g_0$) ,
$$
\lim_{\e\rightarrow 0}\int_0^T dt {\cal L}_{\e}^*(t,\psi(t),\dot\psi(t))
=\int_0^T dt\lim_{\e\rightarrow 0}{\cal L}_{\e}^*(t,\psi(t),\dot\psi(t))
=\int_0^T dt{\cal L}^*(t,\psi(t),\dot\psi(t))
\Eq(500.6)
$$ 
On the other hand, for any $\psi\in\widetilde\GG$ and any 
$\g<\g_0$, $c_1\equiv c_1(\psi)<\infty$ and 
$\vartheta(\bar{\cal S}_{\rho,\gamma/2}(\phi))<\infty$.
Thus, given our choice of the parameters $\zeta$ and $\tau$,
$\widetilde Q
\left(\varepsilon_0,\bar{\cal S}_{\rho,\gamma/2}(\phi),\zeta, c_0\right)$
converges to zero when taking the limit $\e\downarrow 0$ first and the limit
$\varepsilon_0\downarrow 0$ next.

Combining the previous two observations and passing to the
limit $\e\downarrow 0$ in \eqv(500.5) we obtain that 
$$
\liminf_{\e\rightarrow 0}
\e\log\widetilde{\cal P}_{\e,\phi_0}(\BB_{\rho}(\phi))\geq
\cases
-{\displaystyle{\int_0^T dt {\cal L}^*(t,\psi(t),\dot\psi(t))}}
& \hbox{if}\,\, \psi(t_0)=\phi_0\cr
&\cr
-\infty & \hbox{otherwise}
\cr
\endcases
\Eq(500.7)
$$
and since this is true for any $\psi\in\widetilde\GG$,
$$
\eqalign{
\liminf_{\e\rightarrow 0}
\e\log\widetilde{\cal P}_{\e,\phi_0}(\BB_{\rho}(\phi))
\geq &
-\inf_{\scriptstyle{\psi\in\widetilde\GG:}
\atop\scriptstyle{\psi(t_0)=\phi_0}}
{\displaystyle{\int_0^T dt {\cal L}^*(t,\psi(t),\dot\psi(t))}}
\cr
\geq &
-\inf_{\scriptstyle{\psi\in\GG:}
\atop\scriptstyle{\psi(t_0)=\phi_0}}
{\displaystyle{\int_0^T dt {\cal L}^*(t,\psi(t),\dot\psi(t))}}
\cr
}
\Eq(500.8)
$$
where we used that $\GG\subset\widetilde\GG$ in the last line and
where the infimum is $+\infty$ vacuously. But by Lemma 4.15, taking
$\FF=\BB_{\rho}(\phi)\cap{\DD}^{\circ}([0,T])$ therein,
$$
\inf_{\scriptstyle{\psi\in\GG:}
\atop\scriptstyle{\psi(t_0)=\phi_0}}
{\displaystyle{\int_0^T dt {\cal L}^*(t,\psi(t),\dot\psi(t))}}
=
\inf_{\scriptstyle{\psi\in\FF:}
\atop\scriptstyle{\psi(t_0)=\phi_0}}
{\displaystyle{\int_0^T dt {\cal L}^*(t,\psi(t),\dot\psi(t))}}
\Eq(500.9)
$$
and so
$$
\liminf_{\e\rightarrow 0}
\e\log\widetilde{\cal P}_{\e,\phi_0}(\BB_{\rho}(\phi))
\geq 
-\inf_{\scriptstyle{\psi\in \BB_{\rho}(\phi)\cap{\DD}^{\circ}([0,T]):}
\atop\scriptstyle{\psi(t_0)=\phi_0}}
{\displaystyle{\int_0^T dt {\cal L}^*(t,\psi(t),\dot\psi(t))}}
\Eq(500.10)
$$
The lower bound is proven.
\endproof

\proofof{the upper bound \eqv(3.4bis)} To prove the upper bound we first 
combine Lemmata \ver.2 and \ver.4. to get
(with the notation of Lemma \ver.4)
$$
\eqalign{
& \e\log\widetilde{\cal P}_{\e,\phi_0}
\left(
\bar\BB_{\rho}(\phi)
\right)
\cr
\leq &
-
\inf_{{\psi\in\bar\BB_{\rho,0}(\phi):}
\atop{|\psi(t_0)-\phi_0|\leq\bar\zeta+\e\sqrt{d}}}
\inf_{ \psi'(t):\forall_{i=0}^n |\psi'(t_i)- \psi(t_i)|\leq \bar\z}
\sum_{i=1}^n (t_i-t_{i-1}){\cal L}_{\e}^{(r)*}
\left(t_{i-1},\psi'(t_{i-1}),
\sfrac{\psi'(t_i)-\psi'(t_{i-1})}{t_i-t_{i-1}}\right)
}
\Eq(5000.1)
$$
Next we want to use Lemma \ver.7 (ii) to replace sum in 
the right hand side 
by an integral. Before doing this, we observe, however, that the
second  infimum in \eqv(5000.1) will always be realized for 
$\psi'(t_i)$'s for which $\sfrac{\psi'(t_i)-\psi'(t_{i-1})}{t_i-t_{i-1}}\in
D$ (otherwise the infimum takes the value $+\infty$). Thus not only 
can we use Lemma \ver.7 (ii) with $C=\diam\D$, but we
actually have that  $\wt\psi\in\EE([0,T])$. Therefore we may first use
\eqv(5.201) and then replace the infimum over the values $\psi(t_i)$
by an infimum over functions $\wt\psi(t)\in \EE([0,T])$ that are piecewise
linear (p.l.) between the times $t_i$ , i.e. if $\e\t<
\varepsilon_1$,
$$
\eqalign{
& \e\log\widetilde{\cal P}_{\e,\phi_0}
\left(
\bar\BB_{\rho}(\phi)
\right)
\cr
\leq &
-
\inf_{{\psi\in\bar\BB_{\rho,0}(\phi):}
\atop{|\psi(t_0)-\phi_0|\leq\bar\zeta+\e\sqrt{d}}}
\inf_{{\wt \psi(t)\in\EE([0,T]), {\roman p.l.}
}\atop{\forall_{i=0}^n  |\wt\psi(t_i)-\psi(t_i)|\leq \bar\z}}
\int_{0}^T dt  \LL^{(r)*}_{\e}
\left(t,\wt\psi(t),\dot{\wt\psi}(t)\right) 
-3\varepsilon_0T\cr
}
\Eq(5000.2)
$$
Finally (using convexity arguments), 
the two infima can be combined to a single infimum 
over a slightly enlarged set: 
$$
\eqalign{
& \e\log\widetilde{\cal P}_{\e,\phi_0}
\left(
\bar\BB_{\rho}(\phi)
\right)
\leq
-
\inf_{{\psi\in\bar\BB_{\rho+\bar\z}(\phi):}
\atop{{|\psi(t_0)-\phi_0|\leq\bar\zeta+\e\sqrt{d}}\atop
{\forall t\in[0,T], \dist(\psi(t),\L)\leq \bar\z}}}
\int_{0}^T dt  \LL_{\e}^{(r)*}
\left(t,\psi(t),\dot{\psi}(t)\right) 
-3\varepsilon_0T\cr
}
\Eq(50.12)
$$

To conclude the proof of the upper bound what is left to do is to 
pass to the limits $\e\downarrow 0$, $\varepsilon_0\downarrow 0$,
 and $r\downarrow 0$ in \eqv(50.12). 
Note that by Lemma 4.12, for all $r>0$, the function 
$\LL_{\e}^{(r)*} (t,u,v^*)$ is uniformly bounded for all
$t\in\R^+, v^*\in D$, and $u$ such that $\dist(u,\L)\leq r/2$.
Moreover, on the same set it converges uniformly to $\LL^{(r)*}(t,u,v^*)$.
Thus we can use that 
$$
\eqalign{
&\inf_{{\psi\in\bar\BB_{\rho+\bar\z}(\phi):}
\atop{{|\psi(t_0)-\phi_0|\leq\bar\zeta+\e\sqrt{d}}\atop
{\forall t\in[0,T], \dist(\psi(t),\L)\leq \bar\z}}}
\int_{0}^T dt  \LL^{(r)*}_{\e}
\left(t,\psi(t),\dot\psi(t)\right) 
\geq \inf_{{\psi\in\bar\BB_{\rho+\bar\z}(\phi):}
\atop{{|\psi(t_0)-\phi_0|\leq\bar\zeta+\e\sqrt{d}}\atop
{\forall t\in[0,T], \dist(\psi(t),\L)\leq \bar\z}}}
\int_{0}^T dt  \LL^{(r)*}
\left(t,\psi(t),\dot\psi(t)\right) 
\cr&
-\sup_{{\psi\in\bar\BB_{\rho+\bar\z}(\phi):}
\atop{{|\psi(t_0)-\phi_0|\leq\bar\zeta+\e\sqrt{d}}\atop
{\forall t\in[0,T], \dist(\psi(t),\L)\leq \bar\z}}}
\int_{0}^T dt \left[ \LL^{(r)*}_\e
\left(t,\psi(t),\dot\psi(t)\right)-  \LL^{(r)*}
\left(t,\psi(t),\dot\psi(t)\right)\right]
\cr
}
\Eq(50.17)
$$
But
$$
\eqalign{
&\sup_{{\psi\in\bar\BB_{\rho+\bar\z}(\phi):}
\atop{{|\psi(t_0)-\phi_0|\leq\bar\zeta+\e\sqrt{d}}\atop
{\forall t\in[0,T], \dist(\psi(t),\L)\leq \bar\z}}}
\int_{0}^T dt \left[ \LL^{(r)*}_\e{}
\left(t,\psi(t),\dot\psi(t)\right)-  \LL^{(r)*}
\left(t,\psi(t),\dot\psi(t)\right)\right]
\cr&
\leq \sup_{{\psi\in\bar\BB_{\rho+r/2}(\phi):}
\atop{{|\psi(t_0)-\phi_0|\leq r/2}\atop
{\forall t\in[0,T], \dist(\psi(t),\L)\leq r/2}}}
\int_{0}^T dt \left[ \LL^{(r)*}_\e
\left(t,\psi(t),\dot\psi(t)\right)-  \LL^{(r)*}
\left(t,\psi(t),\dot\psi(t)\right)\right]
}
\Eq(50.18)
$$
By Lemma 4.13
 and  dominated convergence, 
the last integral in \eqv(50.18) converges to zero
as $\e\downarrow 0$ uniformly for any $\psi\in\bar\BB_{\rho+r/2}(\phi)$, 
and so \eqv(50.18) converges to zero.
Recall from the proof of the lower bound that $\eta$ and $\t$ were chosen such
that both $\e\t\downarrow 0$ and
$\eta\downarrow 0$ as  $\e\downarrow 0$. Hence $\bar\zeta\downarrow 0$
as  $\e\downarrow 0$. Taking the limit $\e\downarrow 0$ first
and $\varepsilon_0\downarrow 0$ in \eqv(50.12) 
yields that, for any $r>0$,
$$
\limsup_{\e\downarrow 0} \e\log\widetilde{\cal P}_{\e,\phi_0}
\left(
\bar\BB_{\rho}(\phi)
\right)
\leq -
\inf_{{\psi\in\bar\BB_{\rho}(\phi):}
\atop{{\psi(t_0)=\phi_0}\atop {\forall t\in[0,T],\psi(t)\in \L}}}
\int_{0}^T dt  \LL^{(r)*}
\left(t,\psi(t),\dot\psi(t)\right) 
\Eq(50.19)
$$
Finally we must pass to the limit as $r\downarrow 0$. Here the argument is
identical to the one given in [DEW]. It basically relies on Theorem 3.3 in 
[WF] which states that if ${\cal I}$ is a rate function with compact level 
sets $K(s)\equiv\{\psi:{\cal I}(\psi)\leq s\}$, than an  upper bound of the 
form \eqv(50.19)
with rate function ${\cal I}$ is equivalent to the statement that for 
any $c,c'>0$,
there is $\e_0>0$ such that for all $\e\leq \e_0$,
$$
{\cal P}_{\e,\phi_0}\left(\dist(\psi, K(s)\right)
\leq e^{-\frac 1\e (s-c')}
\Eq(50.20)
$$
Therefore, it is enough to show that with 
$K^{(r)}(\psi) \equiv \int_{0}^T dt  \LL^{(r)*}
\left(t,\psi(t),\dot\psi(t)\right) $, and $\bar K(\psi)
\equiv \int_{0}^T dt  {\bar\LL}_{}{}^*
\left(t,\psi(t),\dot\psi(t)\right) $,
for any $s,c,c'>0$, there exists $r>0$ such that
$$
K^{(r)}(s-c)\subset \left\{\psi:\dist(\psi,K(s))\leq c'\right\}
\Eq(50.21)
$$
which is established in Proposition 2.10 of [DEW]. This gives 
the upper bound of Proposition 3.2.\endproof

\newpage
{\headline={\ifodd\pageno\rightheadline \else \leftheadline \fi}}
\def\rightheadline{\it  {Sample path LDP}\hfil\tenrm\folio}
\def\leftheadline{\tenrm \folio \hfil\it  {References}}

\chap{References}0
\item{[AD]} R. Atar and P. Dupuis, ``Large deviations and queueing networks: 
methods for rate function identification, preprint 1998.
\item{[Az]} R. Azencott, ``Petites perturbations al\'eatoires des syst\`emes 
dynamiques: d\'eveloppements asymptotiques'', Bull. Sc. Math. {\bf 109},
253-308 (1985).
\item{[BEW]} M. Bou\'e, P. Dupuis, and R. Ellis, ``Large deviations 
for diffusions with discontinuous statistics'', preprint 1997.
\item{[BEGK]} A. Bovier, M. Eckhoff, V. Gayrard, and M. Klein,
``Metastability in stochastic dynamics of disordered mean field models'',
WIAS-preprint 1998
\item{[BG]} A. Bovier and V. Gayrard, ``Hopfield models as generalized 
             random mean field models'',  in 
            ``Mathematical aspects of spin glasses and neural networks'',  
             A. Bovier and P. Picco (eds.), Progress in Probability {\bf 41},
             Birkh\"auser, Boston, 1998.
\item{[CS]} Chiang Tzuu-Shuh and Sheu Shuenn-Jyi, ``Large deviation of diffusion 
processes and their occupation times with discontinuous limit'', preprint
Academia Sinica, Taipeh, 1998
\item{[DE]} P. Dupuis and  R.S. Ellis, ``A weak convergence approach to the 
theory of large deviations'', Wiley, New York, 1995.
\item{[DE2]} P. Dupuis and R.S. Ellis,  ``The large deviation principle for 
a general class of queueing systems. I.'', Trans.
Amer. Math. Soc. {\bf 347}, 2689-2751 (1995). 
\item{[DEW]} P. Dupuis, R.S. Ellis, and A. Weiss, ``Large deviations for
Markov processes with discontinuous statistics, I: General upper bounds'',
Ann. Probab. {\bf 19}, 1280-1297 (1991).
\item{[DR]} P. Dupuis and K. Ramanan, ``A Skorokhod problem formulation
and large deviation analysis of a processor sharing model'', Queueing 
Systems Theory Appl. {\bf 28}, 109-124 (1998).
\item{[DV]} M.D. Donsker and S.R.S. Varadhan, 
``Asymptotic evaluation of certain
Markov process expectations for large time. III'', Comm. Pure Appl. Math. 
{\bf 29}, 389-461 (1976).
\item{[DZ]} A. Dembo and O. Zeitouni, ``Large deviations techniques and 
applications'', Second edition. Applications of Mathematics {\bf 38},
 Springer, New York, 1998.
\item{[E]} R.S. Ellis, ``Entropy, large deviations, and statistical 
mechanics'', Springer, Berlin-Heidelberg-New York, 1985.
\item{[FW]} M.I. Freidlin and A.D. Wentzell, ``Random perturbations of 
dynamical systems'', Springer, Berlin-Heidelberg-Ney York, 1984.
\item{[vK]} N.G. van Kampen, ``Stochastic processes in physics and
chemistry'', North-Holland, Amsterdam, 1981 (reprinted in 1990).
\item{[IT]} A.D. Ioffe and V.M. Tihomirov, ``Theory of
extremal problems'', Studies in mathematics and its applications {\bf 6},
North-Holland, Amsterdam, 1979. 
\item{[Ki3]} Y. Kifer, ``Random perturbations of 
dynamical systems'', Progress in Probability and Statistics 16, Birkh\"auser,
Boston-Basel, 1988.
\item{[Ki4]} Y. Kifer, ``A discrete time version of the Wentzell-Freidlin 
theory'', Ann. Probab. {\bf 18}, 1676-1692 (1990).
\item{[Ku1]}  T. G. Kurtz, ``Solutions of ordinary differential equations 
as limits of pure jump Markov processes'', J. Appl. Probab. {\bf 7},
49-58 (1970).
\item{[Ku2]} T.G. Kurtz, ``Limit theorems for sequences of jump Markov 
processes approximating ordinary differential processes'', 
J. Appl. Probab. {\bf 8},  344-356 (1971). 
\item{[Mo]} A.A. Mogulskii, ``Large deviations for trajectories of 
multi-dimensional random walks'', Theor. Probab. Appl. {\bf 21}, 300-315 
(1976).
\item{[R]} R.T. Rockafeller, ``Convex analysis'', Princeton University Press,
Princeton, 1970.
\item{[SW]} A. Shwartz and A. Weiss, ``Large deviations for performance 
analysis'', Chapman and Hall, London, 1995.
\item{[W1]} A.D. Wentzell, ``Rough limit theorems on large deviations for 
Markov stochastic processes I.'', 
Theor. Probab. Appl. {\bf 21}, 227-242 (1976). 
\item{[W2]} A.D. Wentzell, ``Rough limit theorems on large deviations for 
Markov stochastic processes II.'', 
Theor. Probab. Appl. {\bf 21}, 499-512 (1976).
\item{[W3]} A.D. Wentzell, ``Rough limit theorems on large deviations for 
Markov stochastic processes III.'', 
Theor. Probab. Appl. {\bf 24}, 675-692 (1979).
\item{[W4]} A.D. Wentzell, ``Rough limit theorems on large deviations for 
Markov stochastic processes IV.'', 
Theor. Probab. Appl. {\bf 27}, 215-234 (1982).

\end